\documentclass{amsart}
\usepackage[all,2cell]{xy}
\UseAllTwocells
\usepackage{amssymb,amsmath,stmaryrd,mathrsfs,amsthm}
\usepackage{eucal}
\usepackage{wasysym}            

\title{Framed bicategories and monoidal fibrations}
\author{Michael Shulman}
\address{Department of Mathematics, University of Chicago\\
  5734 S. University Ave, Chicago IL 60615}
\email{shulman@math.uchicago.edu}

\newcommand{\MonCat}{\ensuremath{\mathcal{M}\mathit{on}\mathcal{C}\mathit{at}}}
\newcommand{\MonCatl}{\ensuremath{\mathcal{M}\mathit{on}\mathcal{C}\mathit{at}_{\ell}}}
\newcommand{\BrMonCat}{\ensuremath{\mathcal{B}\mathit{r}\mathcal{M}\mathit{on}\mathcal{C}\mathit{at}}}
\newcommand{\SymMonCat}{\ensuremath{\mathcal{S}\mathit{ym}\mathcal{M}\mathit{on}\mathcal{C}\mathit{at}}}
\newcommand{\Cat}{\ensuremath{\mathcal{C}\mathit{at}}}

\newcommand{\FibB}{\ensuremath{\mathcal{F}\mathit{ib}_{\sB}}}

\newcommand{\FibcB}{\ensuremath{\mathcal{F}\mathit{ib}_{\mathit{op}\ell,\sB}}}
\newcommand{\calMod}{\ensuremath{\mathcal{M}\mathit{od}}}
\newcommand{\calMat}{\ensuremath{\mathcal{M}\mathit{at}}}
\newcommand{\calDist}{\ensuremath{\mathcal{D}\mathit{ist}}}
\newcommand{\calEx}{\ensuremath{\mathcal{E}\mathit{x}}}
\newcommand{\Cart}{\ensuremath{\mathcal{C}\mathit{art}}}
\newcommand{\biCartpsc}{\ensuremath{\mathit{bi}\mathcal{C}\mathit{art}^{\mathrm{psc}}}}

\newcommand{\Bicatc}{\ensuremath{\mathfrak{Bicat}_{\mathit{op}\ell}}}

\newcommand{\bfBicat}{\ensuremath{\mathbf{Bicat}}}

\newcommand{\Dbl}{\ensuremath{\mathcal{D}\mathit{bl}}}
\newcommand{\Dbll}{\ensuremath{\mathcal{D}\mathit{bl}_{\ell}}}
\newcommand{\Dblc}{\ensuremath{\mathcal{D}\mathit{bl}_{\mathit{op}\ell}}}
\newcommand{\FrBi}{\ensuremath{\mathcal{F}\mathit{r}\mathcal{B}\mathit{i}}}
\newcommand{\FrBil}{\ensuremath{\mathcal{F}\mathit{r}\mathcal{B}\mathit{i}_{\ell}}}
\newcommand{\FrBilq}{\ensuremath{\mathcal{F}\mathit{r}\mathcal{B}\mathit{i}_{\ell}^q}}
\newcommand{\FrBilnq}{\ensuremath{\mathcal{F}\mathit{r}\mathcal{B}\mathit{i}_{\ell,n}^q}}
\newcommand{\FrBiq}{\ensuremath{\mathcal{F}\mathit{r}\mathcal{B}\mathit{i}^q}}
\newcommand{\FrBic}{\ensuremath{\mathcal{F}\mathit{r}\mathcal{B}\mathit{i}_{\mathit{op}\ell}}}
\newcommand{\MF}{\ensuremath{\mathcal{MF}}} 
\newcommand{\MFB}{\ensuremath{\mathcal{MF}_\sB}} 
\newcommand{\MFc}{\ensuremath{\mathcal{MF}_{\mathit{op}\ell}}}
\newcommand{\MFl}{\ensuremath{\mathcal{MF}_\ell}}
\newcommand{\MbF}{\ensuremath{\mathcal{M}\mathit{bi}\mathcal{F}}} 
\newcommand{\MbFl}{\ensuremath{\mathcal{M}\mathit{bi}\mathcal{F}_\ell}} 
\newcommand{\BMF}{\ensuremath{\mathcal{BMF}}} 
\newcommand{\SMF}{\ensuremath{\mathcal{SMF}}} 
\newcommand{\BMFB}{\ensuremath{\mathcal{BMF}_\sB}} 
\newcommand{\SMFB}{\ensuremath{\mathcal{SMF}_\sB}} 
\newcommand{\MFfr}{\ensuremath{\MF^{\mathrm{fr}}}} 
\newcommand{\MFfrc}{\ensuremath{\MFc^{\mathrm{fr}}}}
\newcommand{\MFfrl}{\ensuremath{\MFl^{\mathrm{fr}}}}

\newcommand{\Mod}{\ensuremath{\mathbb{M}\mathbf{od}}}
\newcommand{\CMod}{\ensuremath{\mathbb{CM}\mathbf{od}}}
\newcommand{\Comod}{\ensuremath{\mathbb{C}\mathbf{omod}}}
\newcommand{\Span}{\ensuremath{\mathbb{S}\mathbf{pan}}}
\newcommand{\Dist}{\ensuremath{\mathbb{D}\mathbf{ist}}}
\newcommand{\Mat}{\ensuremath{\mathbb{M}\mathbf{at}}}
\newcommand{\nCob}{\ensuremath{n\mathbb{C}\mathbf{ob}}}
\newcommand{\Ex}{\ensuremath{\mathbb{E}\mathbf{x}}}
\newcommand{\Sp}{\ensuremath{\mathbb{S}\mathbf{p}}}
\newcommand{\Adj}{\ensuremath{\mathbb{A}\mathbf{dj}}}

\newcommand{\bbFr}{\ensuremath{\mathbb{F}\mathbf{r}}}

\newcommand{\ttFam}{\ensuremath{\mathtt{Fam}}}
\newcommand{\ttMod}{\ensuremath{\mathtt{Mod}}}
\newcommand{\ttCMod}{\ensuremath{\mathtt{CMod}}}
\newcommand{\ttArr}{\ensuremath{\mathtt{Arr}}}
\newcommand{\ttRetr}{\ensuremath{\mathtt{Retr}}}
\newcommand{\ttSp}{\ensuremath{\mathtt{Sp}}}

\newcommand{\bbDz}{{\ensuremath{\mathbb{D}_0}}}
\newcommand{\bbDo}{{\ensuremath{\mathbb{D}_1}}}
\newcommand{\bbEz}{{\ensuremath{\mathbb{E}_0}}}
\newcommand{\bbPz}{{\ensuremath{\mathbb{P}_0}}}
\newcommand{\bbPo}{{\ensuremath{\mathbb{P}_1}}}

\newcommand{\sA}{\ensuremath{\mathscr{A}}}
\newcommand{\sB}{\ensuremath{\mathscr{B}}}
\newcommand{\sC}{\ensuremath{\mathscr{C}}}
\newcommand{\sD}{\ensuremath{\mathscr{D}}}
\newcommand{\sE}{\ensuremath{\mathscr{E}}}
\newcommand{\sF}{\ensuremath{\mathscr{F}}}
\newcommand{\sK}{\ensuremath{\mathscr{K}}}
\newcommand{\sS}{\ensuremath{\mathscr{S}}}
\newcommand{\sU}{\ensuremath{\mathscr{U}}}
\newcommand{\sV}{\ensuremath{\mathscr{V}}}
\newcommand{\V}{\ensuremath{\mathscr{V}}}

\newcommand{\calB}{\ensuremath{\mathcal{B}}}
\newcommand{\calD}{\ensuremath{\mathcal{D}}}
\newcommand{\calE}{\ensuremath{\mathcal{E}}}
\newcommand{\calF}{\ensuremath{\mathcal{F}}}
\newcommand{\calK}{\ensuremath{\mathcal{K}}}
\newcommand{\calM}{\ensuremath{\mathcal{M}}}

\newcommand{\bbD}{\ensuremath{\mathbb{D}}}
\newcommand{\bbE}{\ensuremath{\mathbb{E}}}

\newcommand{\bbP}{\ensuremath{\mathbb{P}}}
\newcommand{\bbZ}{\ensuremath{\mathbb{Z}}}

\newcommand{\fa}{\ensuremath{\mathfrak{a}}}
\newcommand{\fl}{\ensuremath{\mathfrak{l}}}

\newcommand{\fr}{\ensuremath{\mathfrak{r}}}
\newcommand{\fs}{\ensuremath{\mathfrak{s}}}

\newcommand{\fu}{\ensuremath{\mathfrak{u}}}
\newcommand{\fx}{\ensuremath{\mathfrak{x}}}

\newcommand{\fatil}{\ensuremath{\widetilde{\mathfrak{a}}}}

\newcommand{\fbar}{\ensuremath{\overline{f}}}
\newcommand{\gbar}{\ensuremath{\overline{g}}}
\newcommand{\kbar}{\ensuremath{\overline{k}}}

\newcommand{\ftil}{\ensuremath{\widetilde{f}}}

\newcommand{\altil}{\ensuremath{\widetilde{\alpha}}}

\newcommand{\ep}{\ensuremath{\varepsilon}}
\newcommand{\pr}{\ensuremath{\pi}}

\newcommand{\Top}{\ensuremath{\mathbf{Top}}}
\newcommand{\Fam}{\ensuremath{\mathbf{Fam}}}
\newcommand{\Set}{\ensuremath{\mathbf{Set}}}
\newcommand{\Ab}{\ensuremath{\mathbf{Ab}}}

\newcommand{\Id}{\ensuremath{\operatorname{Id}}}
\newcommand{\Hom}{\ensuremath{\operatorname{Hom}}}
\newcommand{\Ho}{\ensuremath{\operatorname{Ho}}}

\newcommand{\dn}{\downarrow}
\newcommand{\op}{\ensuremath{^{\mathit{op}}}}

\newcommand{\adj}{\dashv}
\newdir{ >}{{}*!/-10pt/@{>}}
\newdir{ |}{{}*!/-10pt/@{|}}
\newdir{> }{!/10pt/@{>}*{}}
\newcommand{\iso}{\cong}
\newcommand{\eqv}{\simeq}
\newcommand{\too}[1][]{\ensuremath{\overset{#1}{\longrightarrow}}}
\newcommand{\oot}[1][]{\ensuremath{\overset{#1}{\longleftarrow}}}
\renewcommand{\to}{\ensuremath{\rightarrow}}
\newcommand{\toto}{\ensuremath{\rightrightarrows}}
\newcommand{\toot}{\ensuremath{\rightleftarrows}}
\newcommand{\into}{\ensuremath{\hookrightarrow}}
\newcommand{\hto}{\ensuremath{\,\mathaccent\shortmid\rightarrow\,}}
\newcommand{\vto}{\to}
\newcommand{\sto}[2]{\ensuremath{\underset{#1}{\overset{#2}{\Longrightarrow}}}}
\newcommand{\hop}{\ensuremath{^{\mathit{h\cdot{}op}}}}
\newcommand{\maps}{\colon}
\newcommand{\spam}{\,:\!}

\newcommand{\del}{\ensuremath{\partial}}
\newcommand{\sm}{\wedge}
\newcommand{\exsm}{\ensuremath{\barwedge}}
\newcommand{\ten}{\ensuremath{\otimes}}
\newcommand{\xrhd}{\;\overline{\rhd}\;}
\newcommand{\xlhd}{\;\overline{\lhd}\;}

\makeatletter
\ifx\SK@label\undefined\let\SK@label\label\fi
\let\your@thm\@thm
\def\@thm#1#2#3{\gdef\currthmtype{#3}\your@thm{#1}{#2}{#3}}
\def\xlabel#1{{\let\your@currentlabel\@currentlabel\def\@currentlabel
{\currthmtype~\your@currentlabel}
\SK@label{#1@}}\label{#1}}
\def\xref#1{\ref{#1@}}
\makeatother

\newtheorem{thm}{Theorem}[section]
\newtheorem{cor}[thm]{Corollary}
\newtheorem{prop}[thm]{Proposition}
\newtheorem{lem}[thm]{Lemma}
\theoremstyle{definition}
\newtheorem{defn}[thm]{Definition}
\newtheorem{notn}[thm]{Notation}
\newtheorem{notns}[thm]{Notations}
\newtheorem{rmk}[thm]{Remark}
\newtheorem{eg}[thm]{Example}
\newtheorem{egs}[thm]{Examples}

\newenvironment{proofsk}{\begin{proof}[Sketch of Proof]}{\end{proof}}

\makeatletter
\let\c@equation\c@thm
\numberwithin{equation}{section}
\makeatother

\begin{document}

\maketitle

\begin{abstract}
  In some bicategories, the 1-cells are `morphisms' between the
  0-cells, such as functors between categories, but in others they are
  `objects' over the 0-cells, such as bimodules, spans, distributors,
  or parametrized spectra.  Many bicategorical notions do not work
  well in these cases, because the `morphisms between 0-cells', such
  as ring homomorphisms, are missing.  We can include them by using a
  pseudo double category, but usually these morphisms also induce base
  change functors acting on the 1-cells.  We avoid complicated
  coherence problems by describing base change `nonalgebraically',
  using categorical fibrations.  The resulting `framed bicategories'
  assemble into 2-categories, with attendant notions of equivalence,
  adjunction, and so on which are more appropriate for our examples
  than are the usual bicategorical ones.

  We then describe two ways to construct framed bicategories.  One is
  an analogue of rings and bimodules which starts from one framed
  bicategory and builds another.  The other starts from a `monoidal
  fibration', meaning a parametrized family of monoidal categories,
  and produces an analogue of the framed bicategory of spans.
  Combining the two, we obtain a construction which includes both
  enriched and internal categories as special cases.
\end{abstract}

\tableofcontents

\section{Introduction}
\label{sec:introduction}

We begin with the observation that there are really two sorts of
bicategories (or 2-categories).  This fact is well appreciated in
2-categorical circles, but not as widely known as it ought to be.  (In
fact, there are other sorts of bicategory, but we will only be
concerned with two.)

The first sort is exemplified by the 2-category \Cat\ of categories,
functors, and natural transformations.  Here, the 0-cells are
`objects', the 1-cells are maps between them, and the 2-cells are
`maps between maps.'  This sort of bicategory is well-described by the
slogan ``a bicategory is a category enriched over categories.''

The second sort is exemplified by the bicategory \calMod\ of rings,
bimodules, and bimodule homomorphisms.  Here, the 1-cells are
themselves `objects', the 2-cells are maps between them, and the
0-cells are a different sort of `object' which play a `bookkeeping'
role in organizing the relationships between the 1-cells.  This sort
of bicategory is well-described by the slogan ``a bicategory is a
monoidal category with many objects.''

Many notions in bicategory theory work as well for one sort as for the
other.  For example, the notion of 2-functor (including lax 2-functors
as well as pseudo ones) is well-suited to describe morphisms of either
sort of bicategory.  Other notions, such as that of internal
adjunction (or `dual pair'), are useful in both situations, but their
meaning in the two cases is very different.

However, some bicategorical ideas make more sense for one sort of
bicategory than for the other, and frequently it is the second sort
that gets slighted.  A prime example is the notion of equivalence of
0-cells in a bicategory.  This specializes in \Cat\ to equivalence of
categories, which is unquestionably the fundamental notion of
`sameness' for categories.  But in \calMod\ it specializes to Morita
equivalence of rings, which, while very interesting, is not the most
fundamental sort of `sameness' for rings; isomorphism is.

This may not seem like such a big deal, since if we want to talk about
when two rings are isomorphic, we can use the category of rings
instead of the bicategory \calMod.  However, it becomes more acute
when we consider the notion of biequivalence of bicategories, which
involves pseudo 2-functors $F$ and $G$, and equivalences $X\eqv GFX$
and $Y\eqv FGY$.  This is fine for \Cat-like bicategories, but for
\calMod-like bicategories, the right notion of equivalence ought to
include something corresponding to ring isomorphisms instead.  This
problem arose in~\cite[19.3.5]{pht}, where two \calMod-like
bicategories were clearly `equivalent', yet the language did not exist
to describe what sort of equivalence was meant.

Similar problems arise in many other situations, such as the
following.
\begin{enumerate}
\item \Cat\ is a monoidal bicategory in the usual sense, which entails
  (among other things) natural equivalences $(C\times D)\times E \eqv
  C\times (D\times E)$.  But although \calMod\ is `morally monoidal'
  under tensor product of rings, the associativity constraint is
  really a ring isomorphism $(R\ten S)\ten T \iso R\ten (S\ten T)$,
  not an invertible bimodule (although it can be made into one).
\item For \Cat-like bicategories, the notions of pseudonatural
  transformation and modification, making bicategories into a
  tricategory, are natural and useful.  But for \calMod-like
  bicategories, it is significantly less clear what the right sort of
  higher morphisms are.
\item The notion of `biadjunction' is well-suited to adjunctions
  between \Cat-like bicategories, but fails badly for \calMod-like
  bicategories.  Attempts to solve this problem have resulted in some
  work, such as~\cite{verity:base-change,basechange-i,basechange-ii},
  which is closely related to ours.
\end{enumerate}

These problems all stem from essentially the same source: the
bicategory structure does not include the correct `maps between
0-cells', since the 1-cells of the bicategory are being used for
something else.  In this paper, we show how to use an abstract
structure to deal with this sort of situation by incorporating the
maps of 0-cells separately from the 1-cells.  This structure forms a
pseudo double category with extra properties, which we call a
\emph{framed bicategory}.

The first part of this paper is devoted to framed bicategories.  In
\S\S\ref{sec:double-categories}--\ref{sec:duality-theory} we review
basic notions about double categories and fibrations, define framed
bicategories, and prove some basic facts about them.  Then in
\S\S\ref{sec:2cat-frbi}--\ref{sec:involutions} we apply framed
bicategories to resolve the problems mentioned above.  We define lax,
oplax, and strong framed functors and framed transformations, and
thereby obtain three 2-categories of framed bicategories.  We then
apply general 2-category theory to obtain useful notions of framed
equivalence, framed adjunction, and monoidal framed bicategory.

The second part of the paper, consisting of
\S\S\ref{sec:modules}--\ref{sec:monfib-to-framed-ii}, deals with two
important ways of constructing framed bicategories.  The first, which
we describe in \S\ref{sec:modules}, starts with a framed bicategory
\bbD\ and constructs a new framed bicategory $\Mod(\bbD)$ of monoids
and modules in \bbD.  The second starts with a different `parametrized
monoidal structure' called a \emph{monoidal fibration}, and is
essentially the same as the construction of the bicategory of
parametrized spectra in~\cite{pht}.  In
\S\S\ref{sec:monoidal-fibrations}--\ref{sec:closed-monfib} we
introduce monoidal fibrations, and in \S\ref{sec:monfib-to-framed} we
explain the connection to framed bicategories.  Then in
\S\ref{sec:monoids-monfib}, we combine these two constructions and
thereby obtain a natural theory of `categories which are both internal
and enriched'.
\S\S\ref{sec:beck-chev}--\S\ref{sec:monfib-to-framed-ii} are devoted
to the proofs of the main theorems in \S\ref{sec:monfib-to-framed}.

Finally, in the appendices we treat the relationship of framed
bicategories to other work.  This includes the theory of
\emph{connection pairs} and \emph{foldings} in double categories,
various parts of pure bicategory theory, and the bicategorical theory
of \emph{equipments}.  Our conclusion is that they are all, in
suitable senses, equivalent, but each has advantages and
disadvantages, and we believe that framed bicategories are a better
choice for many purposes.

There are two important themes running throughout this paper.  One is
a preoccupation with defining 2-categories and making constructions
2-functorial.  Assembling objects into 2-categories allows us to apply
the theory of adjunctions, equivalences, monads, and so on, internal
to these 2-categories.  Thus, without any extra work, we obtain
notions such as framed adjunctions and framed monads, which behave
much like ordinary adjunctions and monads.  Making various
constructions 2-functorial makes it easy to obtain framed adjunctions
and monads from more ordinary ones.

We do not use very much 2-category theory in this paper, so a passing
acquaintance with it should suffice.  Since we are not writing
primarily for category theorists, we have attempted to avoid or
explain the more esoteric categorical concepts which arise.  A classic
reference for 2-category theory is~\cite{r2cats}; a more modern and
comprehensive one (going far beyond what we will need)
is~\cite{steve:companion}.

The second important theme of this paper is the mixture of `algebraic'
and `nonalgebraic' structures.  A monoidal category is an algebraic
structure: the product is a specified operation on objects.  On the
other hand, a category with cartesian products is a nonalgebraic
structure: the products are characterized by a universal property, and
merely assumed to exist.  We can always make a choice of products to
make a category with products into a monoidal category, but there are
many possible choices, all isomorphic.

There are many technical advantages to working with nonalgebraic
structures.  For example, no coherence axioms are required of a
category with products, whereas a monoidal category requires several.
This advantage becomes more significant as the coherence axioms
multiply.  On the other hand, when doing concrete work, one often
wants to make a specific choice of the structure and work with it
algebraically.  Moreover, not all algebraic structure satisfies an
obvious universal property, and while it can usually be tortured into
doing so, frequently it is easier in these cases to stick with the
algebraic version.

Framed bicategories are a mixture of algebraic and nonalgebraic
notions; the composition of 1-cells is algebraic, while the base
change operations are given nonalgebraically, using a `categorical
fibration'.  Our experience shows that this mixture is very
technically convenient, and we hope to convince the reader of this
too.  In particular, the proof of \xref{monfib->framed} is much
simpler than it would be if we used fully algebraic definitions.  This
is to be contrasted with the similar structures we will consider in
appendices~\ref{sec:folding} and~\ref{sec:equipments}, which are
purely algebraic.

Our intent in this paper is not to present any one particular result,
but rather to argue for the general proposition that framed
bicategories, and related structures, provide a useful framework for
many different kinds of mathematics.  Despite the length of this
paper, we have only had space in it to lay down the most basic
definitions and ideas, and much remains to be said.

The theory of framed bicategories was largely motivated by the desire
to find a good categorical structure for the theory of parametrized
spectra in~\cite{pht}.  The reader familiar with~\cite{pht} should
find the idea of a framed bicategory natural; it was realized clearly
in~\cite{pht} that existing categorical structures were inadequate to
describe the combination of a bicategory with base change operations
which arose naturally in that context.  Another motivation for this
work came from the bicategorical `shadows' of~\cite{kate:traces}, and
a desire to explain in what way they are actually the same as the
horizontal composition in the bicategory; we will do this in the
forthcoming~\cite{me-kate:traces}.

I would like to thank my advisor, Peter May, as well as Kate Ponto,
for many useful discussions about these structures; Tom Fiore, for the
idea of using double categories; and Joachim Kock and Stephan Stolz
for pointing out problems with the original version of \xref{eg:cob}.
The term `framed bicategory' was suggested by Peter May.

\section{Double categories}
\label{sec:double-categories}

As mentioned in the introduction, most of the problems with
\calMod-like bicategories can be traced to the fact that the
`morphisms' of the 0-cells are missing.  Thus, a natural replacement
which suggests itself is a \emph{double category}, a structure which
is like a 2-category, except that it has two types of 1-cells, called
`vertical' and `horizontal', and its 2-cells are shaped like squares.
Double categories go back originally to Ehresmann
in~\cite{ehresmann-cat-str}; a brief introduction can be found
in~\cite{r2cats}.  Other references
include~\cite{multi_funct_i,double-limits,double-adjoints,garner:double-clubs}.

In this section, we introduce basic notions of double categories.  Our
terminology and notation will sometimes be different from that
commonly used.  For example, usually the term `double category' refers
to a strict object, and the weak version is called a `pseudo double
category'.  Since we are primarily interested in the weak version, we
will use the term \emph{double category} for these, and add the word
`strict' if necessary.

\begin{defn}
  A \textbf{double category} \bbD\ consists of a `category of objects'
  $\bbDz$ and a `category of arrows' $\bbDo$, with structure functors
  \begin{align*}
    U&\maps \bbDz\to \bbDo\\
    L,R&\maps \bbDo\rightrightarrows \bbDz\\
    \odot&\maps \bbDo\times_{\bbDz}\bbDo\to \bbDo
  \end{align*}
  (where the pullback is over
  $\bbDo\too[R]\bbDz\overset{L}{\longleftarrow} \bbDo$) such that
  \begin{align*}
    L(U_A) &= A\\
    R(U_A) &= A\\
    L(M\odot N) &= LM\\
    R(M\odot N) &= RM
  \end{align*}
  equipped with natural isomorphisms
  \begin{align*}
    \fa &: (M\odot N) \odot P \too[\iso] M \odot (N \odot P)\\
    \fl &: U_A \odot M \too[\iso] M\\
    \fr &: M \odot U_B \too[\iso] M
  \end{align*}
  such that $L(\fa)$, $R(\fa)$, $L(\fl)$, $R(\fl)$, $L(\fr)$, and
  $R(\fr)$ are all identities, and such that the standard coherence
  axioms for a monoidal category or bicategory (such as Mac Lane's
  pentagon; see~\cite{maclane}) are satisfied.
\end{defn}

We can think of a double category as an internal category in \Cat\
which is suitably weakened, although this is not strictly true because
\Cat\ contains only small categories while we allow \bbDz\ and \bbDo\
to be large categories (but still locally small, that is, having only
a set of morphisms between any two objects).

We call the objects of $\bbDz$ \textbf{objects} or \textbf{0-cells}, and
we call the morphisms of $\bbDz$ \textbf{vertical arrows} and write them
as $f\maps A\to B$.  We call the objects of $\bbDo$ \textbf{horizontal
  1-cells} or just \textbf{1-cells}.  If $M$ is a horizontal 1-cell with
$L(M)=A$ and $R(M)=B$, we write $M\maps A\hto B$, and say that $A$ is
the \textbf{left frame} of $M$ and $B$ is the \textbf{right frame}.  We
use this terminology in preference to the more usual `source' and
`target' because of our philosophy that the horizontal 1-cells are not
`morphisms', but rather objects in their own right which just happen
to be `labeled' by a pair of objects of another type.

A morphism $\alpha\maps M\to N$ of $\bbDo$ with $L(\alpha)=f$ and
$R(\alpha)=g$ is called a \textbf{2-cell}, written $\alpha\maps M
\sto{f}{g} N$, or just $M\too[\alpha] N$, and drawn as follows:
\begin{equation}\label{eq:square}
  \xymatrix{
    A \ar[r]|{|}^{M}  \ar[d]_f \ar@{}[dr]|{\Downarrow\alpha}&
    B\ar[d]^g\\
    C \ar[r]|{|}_N & D
  }.
\end{equation}
We say that $M$ and $N$ are the \textbf{source} and \textbf{target} of
$\alpha$, while $f$ and $g$ are its \textbf{left frame} and \textbf{right
  frame}.  We write the composition of vertical arrows $A\too[f]
B\too[g] C$ and the vertical composition of 2-cells $M\too[\alpha]
N\too[\beta] P$ with juxtaposition, $gf$ or $\beta\alpha$, but we
write the horizontal composition of horizontal 1-cells as $M\odot N$
and that of 2-cells as $\alpha\odot\beta$.

We write horizontal composition `forwards' rather than backwards: for
$M\maps A\hto B$ and $N\maps B\hto C$, we have $M\odot N\maps A\hto
C$.  This is also called `diagrammatic order' and has several
advantages.  First, in examples such as that of rings and bimodules
(\xref{eg:bimodules}), we can define a horizontal 1-cell $M\maps A\hto
B$ to be an $(A,B)$-bimodule, rather than a $(B,A)$-bimodule, and
still preserve the order in the definition $M\odot N = M \ten_B N$ of
horizontal composition.  It also makes it easier to avoid mistakes in
working with 2-cell diagrams; it is easier to compose
\[\xymatrix{A \ar[r]|{|}^M & B \ar[r]|{|}^N & C}\]
and get
\[\xymatrix{A \ar[rr]|{|}^{M\odot N} && C}\]
than to remember to switch the order in which $M$ and $N$ appear every
time horizontal 1-cells are composed.  Finally, it allows us to say
that an adjunction $M\adj N$ in the horizontal bicategory is the same
as a `dual pair' $(M,N)$ (see \S\ref{sec:duality-theory}), with the
left adjoint also being the left dual.

Every object $A$ of a double category has a vertical identity $1_A$
and a horizontal unit $U_A$, every horizontal 1-cell $M$ has an
identity 2-cell $1_M$, every vertical arrow $f$ has a horizontal unit
$U_f$, and we have $1_{U_A} = U_{1_A}$ (by the functoriality of $U$).
We will frequently abuse notation by writing $A$ or $f$ instead of
$U_A$ or $U_f$ when the context is clear.  The important point to
remember is that vertical composition is strictly associative and
unital, while horizontal composition is associative and unital only up
to specified coherent isomorphisms.

Note that if \bbDz\ is the terminal category, then the definition of
double category just says that \bbDo\ is a monoidal category.  We call
such double categories \textbf{vertically trivial}.

We call $\bbDz$ the \textbf{vertical category} of \bbD.  We say that two
objects are isomorphic if they are isomorphic in $\bbDz$, and that two
horizontal 1-cells are isomorphic if they are isomorphic in $\bbDo$.
We will never refer to a horizontal 1-cell as an isomorphism.  A
2-cell whose left and right frames are identities is called
\textbf{globular}.  Note that the constraints $\fa,\fl,\fr$ are globular
isomorphisms, but they are natural with respect to all 2-cells, not
just globular ones.

Every double category \bbD\ has a \textbf{horizontal bicategory} \calD\
consisting of the objects, horizontal 1-cells, and globular 2-cells.
If $A$ and $B$ are objects of \bbD, we write $\bbD(A,B)$ for the
\emph{set} of vertical arrows from $A$ to $B$ and $\calD(A,B)$ for the
\emph{category} of horizontal 1-cells and globular 2-cells from $A$ to
$B$.  It is standard in bicategory theory to say that something holds
\textbf{locally} when it is true of all hom-categories $\calD(A,B)$, and
we will extend this usage to double categories.

We also write $_f\bbD_g(M,N)$ for the set of 2-cells $\alpha$ of the
shape~(\ref{eq:square}).  If $f$ and $g$ are identities, we write
instead $\calD(M,N)$ for the \emph{set} of globular 2-cells from $M$
to $N$.  This may be regarded as shorthand for $\calD(A,B)(M,N)$ and
is standard in bicategory theory.

We now consider some examples.  Note that unlike 1-categories, which
we generally name by their objects, we generally name double
categories by their horizontal 1-cells.

\begin{eg}\xlabel{eg:bimodules}
  Let \Mod\ be the double category defined as follows.  Its objects
  are (not necessarily commutative) rings and its vertical morphisms
  are ring homomorphisms.  A 1-cell $M:A\hto B$ is an
  $(A,B)$-bimodule, and a 2-cell $\alpha:M\sto{f}{g} N$ is an
  \emph{$(f,g)$-bilinear map} $M\to N$, i.e.\ an abelian group
  homomorphism $\alpha:M\to N$ such that $\alpha(amb) =
  f(a)\alpha(m)g(b)$.  This is equivalent to saying $\alpha$ is a map
  of $(A,B)$-bimodules $M\to \,_fN_g$, where $_fN_g$ is $N$ regarded
  as an $(A,B)$-bimodule by means of $f$ and $g$.  The horizontal
  composition of bimodules $M\maps A\hto B$ and $N\maps B\hto C$ is
  given by their tensor product, $M\odot N = M\ten_B N$.  For 2-cells
  \[\xymatrix{A \ar[r]^M\ar[d]_f \ar@{}[dr]|{\alpha} &
    C \ar[d]|g \ar[r]^P \ar@{}[dr]|{\beta} &
    E \ar[d]^h\\
    B \ar[r]_N & D \ar[r]_Q & F}\]
  we define $\alpha\odot\beta$ to be the composite
  \begin{equation*}
    \xymatrix{
      M\ten_C P \ar[rr]^<>(.5){\alpha\ten\beta} &&
      {_fN_g}\ten_C \,_gQ_h \ar@{->>}[r] &
      {_fN} \ten_D Q_h \iso \,_f(N\ten_D Q)_h
    }.
  \end{equation*}

  This example may be generalized by replacing \Ab\ with any monoidal
  category \sC\ that has coequalizers preserved by its tensor product,
  giving the double category $\Mod(\sC)$ of monoids, monoid
  homomorphisms, and bimodules in \sC.  If $\sC=\mathbf{Mod}_R$ is the
  category of modules over a commutative ring $R$, then the resulting
  double category $\Mod(\sC) = \Mod(R)$ is made of $R$-algebras,
  $R$-algebra homomorphisms, and bimodules over $R$-algebras.

  Similarly, we define the double category \CMod\ whose objects are
  commutative rings, and if \sC\ is a symmetric monoidal category, we
  have $\CMod(\sC)$.
\end{eg}

\begin{eg}\xlabel{eg:spans}
  Let \sC\ be a category with pullbacks, and define a double category
  $\Span(\sC)$ whose vertical category is \sC, whose 1-cells $A\hto B$
  are spans $A\leftarrow C \rightarrow B$ in \sC, and whose 2-cells
  are commuting diagrams:
  \[\xymatrix{
    A\ar[d] & C\ar[l]\ar[d]\ar[r] & B\ar[d] \\
    D       & F\ar[l]\ar[r]       & E \\
  }\]
  in \sC.  Horizontal composition is by pullback.
\end{eg}

\begin{eg}\xlabel{eg:ex}
  There is a double category of parametrized spectra called \Ex, whose
  construction is essentially contained in~\cite{pht}.  The vertical
  category is a category of (nice) topological spaces, and a 1-cell
  $A\hto B$ is a spectrum parametrized over $A\times B$ (or $B\times
  A$; see the note above about the order of composition).

  In~\cite{pht} this structure is described only as a bicategory with
  `base change operations', but it is pointed out there that existing
  categorical structures do not suffice to describe it.  We will see
  in \S\ref{sec:monfib-to-framed} how this sort of structure gives
  rise, quite generally, not only to a double category, but to a
  framed bicategory, which supplies the missing categorical structure.
\end{eg}

\begin{eg}\xlabel{eg:distributors}
  Let \V\ be a complete and cocomplete closed symmetric monoidal
  category, such as \Set, \Ab, $\mathbf{Cat}$, or a convenient
  cartesian closed subcategory of topological spaces, and define a
  double category $\Dist(\V)$ as follows.  Its objects are (small)
  categories enriched over \V, or \V-categories.  Its vertical arrows
  are \V-functors, its 1-cells are \emph{\V-distributors}, and its
  2-cells are \V-natural transformations.  (Good references for
  enriched category theory include~\cite{kelly} and~\cite{dubuc}.)  A
  \V-distributor $H\maps \sB\hto\sA$ is simply a \V-functor $H\maps
  \sA\op\otimes\sB\to\V$.  When \sA\ and \sB\ have one object, they
  are just monoids in \V, and a distributor between them is a bimodule
  in \V; thus we have an inclusion $\Mod(\V)\into\Dist(\V)$.
  Horizontal composition of distributors is given by the coend
  construction, also known as `tensor product of functors'.

  In the bicategorical literature, distributors are often called
  `bimodules' or just `modules', but we prefer to reserve that term
  for the classical one-object version.  The term
  `distribut\emph{or}', due to Benabou, is intended to suggest a
  generalization of `funct\emph{or}', just as in analysis a
  `distribut\emph{ion}' is a generalized `funct\emph{ion}'.  The term
  `profunctor' is also used for these objects, but we prefer to avoid
  it because a distributor is nothing like a pro-object in a functor
  category.
\end{eg}

\begin{eg}\xlabel{eg:cob}
  We define a double category \nCob\ as follows.  Its vertical
  category consists of oriented $(n-1)$-manifolds without boundary and
  diffeomorphisms.  A 1-cell $M\hto N$ is a (possibly thin)
  $n$-dimensional cobordism from $M$ to $N$, and a 2-cell is a
  compatible diffeomorphism.  Horizontal composition is given by
  gluing of cobordisms.

  More formally, if $A$ and $B$ are oriented $(n-1)$-manifolds, a
  horizontal 1-cell $M\maps B\hto A$ is \emph{either} a diffeomorphism
  $A\iso B$ (regarded as a `thin' cobordism from $A$ to $B$), or an
  $n$-manifold with boundary $M$ equipped with a `collar' map
  \[(A\op\sqcup B)\times [0,1) \to M\]
  which is a diffeomorphism onto its image and restricts to a
  diffeomorphism
  \[A\op \sqcup B \iso \del M.
  \]
  (Here $A\op$ means $A$ with the opposite orientation.)  The unit is
  the identity $1_A$, regarded as a thin cobordism.
\end{eg}

\begin{eg}\xlabel{eg:adj}
  The following double category is known as \Adj.  Its objects are
  categories, and its \emph{horizontal} 1-cells are functors.  Its
  vertical arrows $C\to D$ are adjoint pairs of functors $f_!\maps
  C\toot D\spam f^*$.  We then seem to have two choices for the
  2-cells; a 2-cell with boundary
  \[\xymatrix{A \ar[r]^h\ar@/^1mm/[d]^{f_!} & C \ar@/^1mm/[d]^{g_!}\\
    B\ar[r]_k \ar@/^1mm/[u]^{f^*} & D \ar@/^1mm/[u]^{g^*}}\]
  could be chosen to be either a natural transformation $g_!h\to kf_!$
  or a natural transformation $hf^*\to g^* k$.  However, it turns out
  that there is a natural bijection between natural transformations
  $g_!h\to kf_!$ and $hf^*\to g^* k$ which respects composition, so it
  doesn't matter which we pick.  Pairs of natural transformations
  corresponding to each other under this bijection are called
  \textbf{mates}; the mate of a transformation $\alpha\maps hf^*\to g^*
  k$ is given explicitly as the composite
  \[\xymatrix{g_!h \ar[r]^-{g_!h \eta} &
    g_!hf^*f_! \ar[r]^-{g_!\alpha f_!} &
    g_!g^*kf_! \ar[r]^-{\ep k f_!} & k f_!}\]
  where $\eta$ is the unit of the adjunction $f_!\adj f^*$ and $\ep$
  is the counit of the adjunction $g_!\adj g^*$.  The inverse
  construction is dual.

  More generally, if \calK\ is any (strict) 2-category, we can define
  the notion of an \emph{adjunction} internal to \calK: it consists of
  morphisms $f\maps A\to B$ and $g\maps B\to A$ together with 2-cells
  $\eta\maps 1_A \Rightarrow gf$ and $\ep\maps fg \Rightarrow 1_B$
  satisfying the usual triangle identities.  We can then define a
  double category $\Adj(\calK)$ formed by objects, morphisms,
  adjunctions, and mate-pairs internal to \calK.

  These double categories have a different flavor than the others
  introduced above.  We mention them partly to point out that double
  categories have uses other than those we are interested in, and
  partly because we will need the notion of mates later on.  More
  about mate-pairs in 2-categories and their relationship to \Adj\ can
  be found in~\cite{r2cats}; one fact we will need is that if $h$ and
  $k$ are identities, then $\alpha$ is an isomorphism if and only if
  its mate is an isomorphism.
\end{eg}

\section{Review of the theory of fibrations}
\label{sec:fibrations}

Double categories incorporate both the 1-cells of a \calMod-like
bicategory and the `morphisms of 0-cells', but there is something
missing.  An important feature of all our examples is that the 1-cells
can be `base changed' along the vertical arrows.  For example, in
\Mod, we can extend and restrict scalars along a ring homomorphism.

An appropriate abstract structure to describe these base change
functors is the well-known categorical notion of a `fibration'.  In
this section we will review some of the theory of fibrations, and then
in \S\ref{sec:framed} we will apply it to base change functors in
double categories.  All the material in this section is standard.  The
theory of fibrations is originally due to Grothendieck and his school;
see, for example~\cite[Expos\'e~VI]{sga1}.  Modern references
include~\cite[B1.3]{elephant1} and~\cite[Ch.~8]{borceaux:handbook-2}.
More abstract versions can be found in the 2-categorical literature,
such as~\cite{street:fibi}.

\begin{defn}
  Let $\Phi:\sA\to\sB$ be a functor, let $f:A\to C$ be an arrow in
  \sB, and let $M$ be an object of \sA\ with $\Phi(M) = C$.  An arrow
  $\phi:f^*M\to M$ in \sA\ is \textbf{cartesian over $f$} if, firstly,
  $\Phi(\phi) = f$:
  \begin{equation*}
    \xymatrix{
      f^*M \ar[r]^\phi \ar@{ |-> }[d] & M \ar@{ |-> }[d]\\
      A \ar[r]_f & C
    }
  \end{equation*}
  and secondly, whenever $\psi:N\to M$ is an arrow in \sA\ and $g:B\to
  A$ is an arrow in \sB\ such that $\Phi(\psi) = fg$, there is a
  unique $\chi$ such that $\psi = \phi\chi$ and $\Phi(\chi)=g$:
  \begin{equation*}
    \xymatrix{
      N \ar[drrr]^\psi \ar@{ |-> }[dd] \ar@{-->}[dr]_\chi\\
      &f^*M \ar[rr]_\phi \ar@{ |-> }[dd] && M \ar@{ |-> }[dd]\\
      B \ar[dr]_g \ar[drrr]\\
      &A \ar[rr]_f && C
    }
  \end{equation*}
  We say that $\Phi$ is a \textbf{fibration} if for every $f\maps A\to
  C$ and $M$ with $\Phi(M)=C$, there exists a cartesian arrow
  $\phi_{f,M}:f^*M\to M$ over $f$.  If $\Phi$ is a fibration, a
  \textbf{cleavage} for $\Phi$ is a choice, for every $f$ and $M$, of
  such a $\phi_{f,M}$.  The cleavage is \textbf{normal} if $\phi_{1_A,M}
  = 1_M$; it is \textbf{split} if $\phi_{g,M}\phi_{f,g^*M} =
  \phi_{gf,M}$ for all composable $f,g$.
\end{defn}

For an arrow $f\maps A\to B$, we think of $f^*$ as a `base change'
operation that maps the fiber $\sA_B$ (consisting of all objects over
$B$ and morphisms over $1_B$) to the fiber $\sA_A$.  We think of $\sA$
as `glued together' from the fiber categories $\sA_B$ as $B$ varies,
using the base change operations $f^*$.  We think of the whole
fibration as `a category parametrized by \sB'.

\begin{eg}\xlabel{eg:fib-rings}
  Let $\mathbf{Ring}$ be the category of rings, and $\mathbf{Mod}$ be
  the category of pairs $(R,M)$ where $R$ is a ring and $M$ is an
  $R$-module, with morphisms consisting of a ring homomorphism $f$ and
  an $f$-equivariant module map.  Then the forgetful functor
  $\ttMod\maps \mathbf{Mod}\to\mathbf{Ring}$, which sends $(R,M)$ to
  $R$, is a fibration.  If $M$ is an $R$-module and $f\maps S\to R$ is
  a ring homomorphism, then if we denote by $f^*M$ the abelian group
  $M$ regarded as an $S$-module via $f$, the identity map of $M$
  defines an $f$-equivariant map $f^*M\to M$, which is a cartesian
  arrow over $f$.

  Note that the fiber $\sA_R$ is the ordinary category of $R$-modules.
  Thus we may say that modules form a category parametrized by rings.
\end{eg}

\begin{eg}\xlabel{eg:fib-arr}
  Let \sC\ be a category with pullbacks, let $\sC^\dn$ denote the
  category of arrows in \sC\ (whose morphisms are commutative
  squares), and let $\ttArr_\sC\maps \sC^\dn\to\sC$ take each arrow to
  its codomain.  Then $\ttArr_\sC$ is a fibration; a commutative
  square is a cartesian arrow in $\sC^\dn$ precisely when it is a
  pullback square.  This fibration is sometimes referred to as the
  \emph{self-indexing} of \sC.
\end{eg}

We record some useful facts about fibrations.

\begin{prop}\xlabel{fib-facts}
  Let $\Phi\maps \sA\to\sB$ be a fibration.
  \begin{enumerate}
  \item The composite of cartesian arrows is cartesian.\label{item:comp-cart}
  \item If $\phi\maps (fg)^*M\to M$ and $\psi\maps g^*M\to M$ are
    cartesian over $fg$ and $g$ respectively, and $\chi\maps
    (fg)^*M\to g^*M$ is the unique factorization of $\psi$ through
    $\phi$ lying over $f$, then $\chi$ is cartesian.\label{item:cart-fact-cart}
  \item If if $\phi\maps f^*M\to M$ and $\phi'\maps (f^*M)'\to M$ are
    two cartesian lifts of $f$, then there is a unique isomorphism
    $f^*M\iso (f^*M)'$ commuting with $\phi$ and $\phi'$.\label{item:cart-uniq}
  \item Any isomorphism in \sA\ is cartesian.\label{item:iso-cart}
  \item If $f$ is an isomorphism in \sB, then any cartesian lift of
    $f$ is an isomorphism.\label{item:iso-lift-iso}
  \end{enumerate}
\end{prop}

In \xref{eg:fib-rings}, there is a `canonical' choice of a cleavage,
but this is not true in \xref{eg:fib-arr}, since pullbacks are only
defined up to isomorphism.  \xref{fib-facts}\ref{item:cart-uniq} tells
us that more generally, cleavages in a fibration are unique up to
canonical isomorphism.  Thus, a fibration is a `nonalgebraic' approach
to defining base change functors: the operation $f^*$ is characterized
by a universal property, and the definition merely stipulates that an
object satisfying that property exists, rather than choosing a
particular such object as part of the structure.  In the terminology
of~\cite{makkai:comparing}, they are \emph{virtual} operations.

The `algebraic' notion corresponding to a fibration $\Phi\maps
\sA\to\sB$ is a pseudofunctor $P\maps \sB\op \to \Cat$.  Given a
fibration $\Phi$, if we choose a cleavage, then we obtain, for each
$f\maps A\to B$ in \sB, a functor $f^*\maps \sA_B\to \sA_A$.  If we
define $P(A) = \sA_A$ and $P(f) = f^*$, the uniqueness-up-to-iso of
cartesian lifts makes $P$ into a pseudofunctor.  Conversely, given a
pseudofunctor $P\maps \sB\op \to \Cat$, we can build a fibration over
\sB\ whose fiber over $A$ is $P(A)$.  (This is sometimes called the
`Grothendieck construction'.)

In order to state the full 2-categorical sense in which these
constructions are inverse equivalences, we need to introduce the
morphisms and transformations between fibrations.  Consider a
commuting square of functors
\begin{equation}
  \xymatrix{\sA \ar[r]^{F_1}\ar[d]_{\Phi} & \sA' \ar[d]^{\Phi'}\\
    \sB \ar[r]_{F_0} & \sB'}\label{eq:mor-fib}
\end{equation}
where $\Phi$ and $\Phi'$ are fibrations, and let $\phi\maps g^*M\to M$
be cartesian in $\sA$ over $g$.  Then we have $F_1(\phi)\maps
F_1(g^*M)\to F_1M$ in $\sA'$ over $F_0(g)$.  But since $\Phi'$ is a
fibration, there is a cartesian arrow $\psi\maps (F_0g)^*(F_1M) \to
F_1M$ over $F_0(g)$, so $F_1(\phi)$ factors uniquely through it,
giving a canonical map
\begin{equation}
  F_1(g^*M) \too (F_0g)^*(F_1M)\label{eq:mor-fib-oplax}
\end{equation}
which is an isomorphism if and only if $F_1(\phi)$ is cartesian.

It should thus be unsurprising that any commuting
square~(\ref{eq:mor-fib}) gives rise to an oplax natural
transformation between the corresponding pseudofunctors.  Recall that
an \emph{oplax natural transformation} between pseudofunctors
$P,Q\maps \sB\op\to\Cat$ consists of functors $\phi_x\maps Px\to Qx$
and natural transformations
\[\xymatrix{Px \ar[r]^{Pg}\ar[d]_{\phi_x} \drtwocell\omit{\phi_g}
  & Py \ar[d]^{\phi_y}\\
  Qx \ar[r]_{Qg} & Qy}\]
satisfying appropriate coherence conditions.  In a \emph{lax natural
  transformation}, the 2-cells go the other direction, and in a
\emph{pseudo natural transformation} the 2-cells are invertible.

\begin{defn}
  Any commuting square of functors~(\ref{eq:mor-fib}) is called an
  \textbf{oplax morphism of fibrations}.  It is a \textbf{strong morphism
    of fibrations} if whenever $\phi$ is a cartesian arrow in $\sA$
  over $g$, then $F_1(\phi)$ is cartesian in $\sA'$ over $F_0(g)$.  If
  $F_0$ is an identity $\sB=\sB'$, then we say $F_1$ is a morphism
  \textbf{over \sB}.

  A \textbf{transformation of fibrations} between two (oplax) morphisms
  of fibrations is just a pair of natural transformations, one lying
  above the other.  If the two morphisms are over \sB, the
  transformation is \textbf{over \sB} if its downstairs component is the
  identity.
\end{defn}

\begin{prop}\xlabel{fib-psfr-eqv}
  Let $\FibcB$ denote the 2-category of fibrations over \sB, oplax
  morphisms of fibrations over \sB, and transformations over \sB, and let
  $[\sB\op,\Cat]_{\mathit{op}\ell}$ denote the 2-category of pseudofunctors
  $\sB\op\to\Cat$, oplax natural transformations, and modifications.
  Then the above constructions define an equivalence of 2-categories
  \[\FibcB\eqv [\sB\op,\Cat]_{\mathit{op}\ell}.\]
  If we restrict to the strong morphisms of fibrations over \sB\ on
  the left and the pseudo natural transformations on the right, we
  again have an equivalence
  \[\FibB\eqv [\sB\op,\Cat].\]
\end{prop}

Compared to pseudofunctors, fibrations have the advantage that they
incorporate all the base change functors $f^*$ and all their coherence
data automatically.  We must remember, however, that the functors
$f^*$ are not determined uniquely by the fibration, only up to natural
isomorphism.

If $\Phi$ is a functor such that $\Phi\op\maps \sA\op\to\sB\op$ is a
fibration, we say that $\Phi$ is an \textbf{opfibration}.  (The term
`cofibration' used to be common, but this carries the wrong intuition
for homotopy theorists, since an opfibration is still characterized by
a \emph{lifting} property.)  The cartesian arrows in $\sA\op$ are
called \textbf{opcartesian} arrows in \sA.  A cleavage for an
opfibration consists of opcartesian arrows $M\to f_!M$, giving rise to
a functor $f_!\maps \sA_A\to \sA_B$ for each morphism $f\maps A\to B$
in \sB.

For any opfibration, the collection of functors $f_!$ forms a
\emph{covariant} pseudofunctor $\sB\to\Cat$, and conversely, any
covariant pseudofunctor gives rise to an opfibration.  A commutative
square~(\ref{eq:mor-fib}) in which $\Phi$ and $\Phi'$ are opfibrations
is called a \textbf{lax morphism of opfibrations}, and it is
\textbf{strong} if $F_1$ preserves opcartesian arrows; these correspond
to lax and pseudo natural transformations, respectively.

\begin{prop}\xlabel{fib-opfib-adjt}
  A fibration $\Phi\maps \sA\to\sB$ is also an opfibration precisely
  when all the functors $f^*$ have left adjoints $f_!$.
\end{prop}
\begin{proof}
  By definition of $f^*$, there is a natural bijection between
  morphisms $M\to N$ in \sA\ lying over $f\maps A\to B$ and morphisms
  $M\to f^*N$ in the fiber $\sA_A$.  But if $\Phi$ is also an
  opfibration, these morphisms are also bijective to morphisms
  $f_!M\to N$ in $\sA_B$, so we have an adjunction $\sA_A(M,f^*N)\iso
  \sA_B(f_!M,N)$ as desired.  The converse is straightforward.
\end{proof}

We will refer to a functor which is both a fibration and an
opfibration as a \textbf{bifibration}.  A square~(\ref{eq:mor-fib}) in
which $\Phi$ and $\Phi'$ are bifibrations is called a \textbf{lax
  morphism of bifibrations} if $F_1$ preserves cartesian arrows, an
\textbf{oplax morphism of bifibrations} if it preserves opcartesian
arrows, and a \textbf{strong morphism of bifibrations} if it preserves
both.

\begin{egs}
  The fibration $\ttMod\maps \mathbf{Mod}\to\mathbf{Ring}$ is in fact
  a bifibration; the left adjoint $f_!$ is given by extension of
  scalars.  For any category \sC\ with pullbacks, the fibration
  $\ttArr_\sC\maps \sC^\dn\to\sC$ is also a bifibration; the left
  adjoint $f_!$ is given by composing with $f$.
\end{egs}

In many cases, the functors $f^*$ also have \emph{right} adjoints,
usually written $f_*$.  These functors are not as conveniently
described by a fibrational condition, but we will see in
\S\ref{sec:duality-theory} that in a framed bicategory, they can be
described in terms of base change objects and a closed structure.  We
say that a fibration is a \textbf{$*$-fibration} if all the functors
$f^*$ have right adjoints $f_*$.  Similarly we have a
\textbf{$*$-bifibration}, in which every morphism $f$ gives rise to an
adjoint string $f_!\adj f^*\adj f_*$.

\begin{egs}
  \ttMod\ is a $*$-bifibration; the right adjoints are given by
  coextension of scalars.  $\ttArr_\sC$ is a $*$-bifibration precisely
  when \sC\ is \emph{locally cartesian closed} (that is, each slice
  category $\sC/X$ is cartesian closed).
\end{egs}

Often the mere existence of left or right adjoints is insufficient,
and we need to require a commutativity condition.  We will explore
this further in \S\ref{sec:beck-chev}.

\section{Framed bicategories}
\label{sec:framed}

Morally speaking, a framed bicategory is a double category in which
the 1-cells can be restricted and extended along the vertical arrows.
We will formalize this by saying that $L$ and $R$ are bifibrations.
Thus, for any $f\maps A\to B$ in \bbDz, there will be two different
functors which should be called $f^*$, one arising from $L$ and one
from $R$.  We distinguish by writing the first on the left and the
second on the right.  In other words, $f^*M$ is a horizontal 1-cell
equipped with a cartesian 2-cell
\[\xymatrix{A \ar[r]|{|}^{f^*M}\ar[d]_f \ar@{}[dr]|{\mathrm{cart}} & D \ar@{=}[d]\\
  B \ar[r]|{|}_M & D}\]
while $Mg^*$ is equipped with a cartesian 2-cell
\[
\begin{array}{c}
  \xymatrix{B \ar[r]|{|}^{Mg^*}\ar@{=}[d] \ar@{}[dr]|{\mathrm{cart}} & C \ar[d]^g\\
    B \ar[r]|{|}_M & D}
\end{array}.
\]
A general cartesian arrow in \bbDo\ lying over $(f,g)$ in
$\bbDz\times\bbDz$ can then be written as $f^*Mg^*\sto{f}{g} M$.  We do
similarly for opcartesian arrows and the corresponding functors $f_!$.
We refer to $f^*$ as \textbf{restriction} and to $f_!$ as
\textbf{extension}.  If $f^*$ also has a right adjoint $f_*$, we refer
to it as \textbf{coextension}.

It is worth commenting explicitly on what it means for a 2-cell in a
double category to be cartesian or opcartesian.  Suppose given a
`niche' of the form
\[\xymatrix{A \ar[d]_f & C \ar[d]^g\\
  B \ar[r]|{|}_M & D}\]
in a double category \bbD.  This corresponds to an object $M\in\bbDo$
and a morphism $(f,g)\maps (A,C)\to (B,D) = (L,R)(M)$ in
$\bbDz\times\bbDz$.  A cartesian lifting of this morphism is a 2-cell
\[\xymatrix{A \ar[r]|{|}^{f^*Mg^*}\ar[d]_f \ar@{}[dr]|{\mathrm{cart}} &
  C \ar[d]^g\\
  B \ar[r]|{|}_M & D}\]
such that any 2-cell of the form
\[\xymatrix{E \ar[r]|{|}^N\ar[d]_{fh} \ar@{}[dr]|{\Downarrow} & F \ar[d]^{gk}\\
  B\ar[r]|{|}_M & D}\]
factors uniquely as follows:
\[\xymatrix@C=3pc{E \ar[r]|{|}^N \ar[d]_h \ar@{}[dr]|{\Downarrow} & F \ar[d]^{k}\\
  A \ar[r]|{f^*Mg^*}\ar[d]_f \ar@{}[dr]|{\mathrm{cart}} &
  C \ar[d]^g\\
  B \ar[r]|{|}_M & D.}\]
In particular, if $h=1_A$ and $k=1_C$, this says that any 2-cell
\[\xymatrix{A \ar[r]|{|}^N\ar[d]_{f} \ar@{}[dr]|{\Downarrow} & C \ar[d]^{g}\\
  B\ar[r]|{|}_M & D}\]
can be represented by a globular 2-cell
\[\xymatrix{A \ar[r]|{|}^N\ar@{=}[d] \ar@{}[dr]|{\Downarrow} & B \ar@{=}[d]\\
  A \ar[r]|{|}_{f^*Mg^*} & B.}\]
Therefore, `all the information' about the 2-cells in a framed
bicategory will be carried by the globular 2-cells and the base change
functors.  In particular, we can think of \bbD\ as `the bicategory
\calD\ equipped with base change functors'.  This can be made precise;
see appendix~\ref{sec:equipments}.

The interaction of fibrational conditions with the double category
structure has further implications.  It is reasonable to expect that
restriction and extension will commute with horizontal composition;
thus we will have $f^*(M\odot N)g^* \iso f^*M\odot Ng^*$.  This
implies, however, that for any 1-cell $M\maps B\hto C$ and arrow
$f\maps A\to B$, we have
\[f^*M \iso f^*(U_B\odot M) \iso f^*U_B \odot M,\]
and hence the base change \emph{functor} $f^*$ can be represented by
horizontal composition with the special object $f^*U_B$, which we call
a \emph{base change object}.

In the case of \Mod, this is the standard fact that restricting along
a ring homomorphism $f\maps A\to B$ is the same as tensoring with the
$(A,B)$-bimodule ${}_fB$, by which we mean $B$ regarded as an
$(A,B)$-bimodule via $f$ on the left.  For this reason, we write
${}_fB$ for the base change object $f^*U_B$ in any double category.
Similarly, we write $B_f$ for $U_B f^*$.

The existence of such base change objects, suitably formalized, turns
out to be sufficient to ensure that \emph{all} restrictions exist.
This formalization of base change objects can be given in an
essentially diagrammatic way, which moreover is self-dual.  Thus, it
it is also equivalent to the existence of \emph{extensions}.  This is
the content of the following result.

\begin{thm}\xlabel{thm:framed}
  The following conditions on a double category \bbD\ are equivalent.
  \begin{enumerate}
  \item $(L,R):\bbDo\to \bbDz\times\bbDz$ is a fibration.\label{item:frbi-fib}
  \item $(L,R):\bbDo\to \bbDz\times\bbDz$ is an opfibration.\label{item:frbi-opfib}
  \item For every vertical arrow $f\maps A\to B$, there exist 1-cells\label{item:frbi-compconj}
    ${}_fB\maps A\hto B$ and $B_f\maps B\hto A$ together with 2-cells
    \begin{equation}
      \begin{array}{c}
        \xymatrix{
          \ar[r]|<>(.5){|}^<>(.5){{}_fB} \ar[d]_f \ar@{}[dr]|\Downarrow
          & \ar@{=}[d]\\
          \ar[r]|<>(.5){|}_<>(.5){U_B} & }
      \end{array}\quad,\quad
      \begin{array}{c}
        \xymatrix{
          \ar[r]|<>(.5){|}^<>(.5){B_f} \ar@{=}[d] \ar@{}[dr]|\Downarrow
          & \ar[d]^f\\
          \ar[r]|<>(.5){|}_<>(.5){U_B} & }
      \end{array}\quad,\quad
      \begin{array}{c}
        \xymatrix{
          \ar[r]|<>(.5){|}^<>(.5){U_A} \ar[d]_f \ar@{}[dr]|\Downarrow
          & \ar@{=}[d]\\
          \ar[r]|<>(.5){|}_<>(.5){B_f} & }
      \end{array}, \quad\text{and}\quad
      \begin{array}{c}
        \xymatrix{
          \ar[r]|<>(.5){|}^<>(.5){U_A} \ar@{=}[d] \ar@{}[dr]|\Downarrow
          & \ar[d]^f\\
          \ar[r]|<>(.5){|}_<>(.5){{}_fB} & }
      \end{array}\label{eq:frbi-compconj}
    \end{equation}
    such that the following equations hold.
    \begin{align}
      \begin{array}{c}
        \xymatrix{
          \ar[r]|<>(.5){|}^<>(.5){U_A} \ar@{=}[d] \ar@{}[dr]|\Downarrow
          & \ar[d]^f\\
          \ar[r]|<>(.5){{}_fB} \ar[d]_f \ar@{}[dr]|\Downarrow
          & \ar@{=}[d]\\
          \ar[r]|<>(.5){|}_<>(.5){U_B} & }
      \end{array} &= 
      \begin{array}{c}
        \xymatrix{ \ar[r]|<>(.5){|}^<>(.5){U_A} \ar[d]_f
          \ar@{}[dr]|{\Downarrow U_f} &  \ar[d]^f\\
          \ar[r]|<>(.5){|}_<>(.5){U_B} & }
      \end{array}
      &
      \begin{array}{c}
        \xymatrix{
          \ar[r]|<>(.5){|}^<>(.5){U_A} \ar[d]_f \ar@{}[dr]|\Downarrow
          & \ar@{=}[d]\\
          \ar[r]|<>(.5){B_f} \ar@{=}[d] \ar@{}[dr]|\Downarrow
          & \ar[d]^f\\
          \ar[r]|<>(.5){|}_<>(.5){U_B} & }
      \end{array} &= 
      \begin{array}{c}
        \xymatrix{ \ar[r]|<>(.5){|}^<>(.5){U_A} \ar[d]_f
          \ar@{}[dr]|{\Downarrow U_f} &  \ar[d]^f\\
          \ar[r]|<>(.5){|}_<>(.5){U_B} & }
      \end{array}\label{eq:frbi-compconj-1}
      \\
      \begin{array}{c}
        \xymatrix{
          \ar[rr]|<>(.5){|}^<>(.5){{}_fB} \ar@{}[drr]|\iso \ar@{=}[d] &&
          \ar@{=}[d] \\
          \ar[r]|<>(.5){|}^<>(.5){U_A} \ar@{=}[d] \ar@{}[dr]|\Downarrow &
          \ar[r]|<>(.5){|}^<>(.5){{}_fB} \ar[d]_f \ar@{}[dr]|\Downarrow
          & \ar@{=}[d]\\
          \ar[r]|<>(.5){|}_<>(.5){{}_fB} &
          \ar[r]|<>(.5){|}_<>(.5){U_B} &\\
          \ar[rr]|<>(.5){|}_<>(.5){{}_fB} \ar@{}[urr]|\iso \ar@{=}[u] &&
          \ar@{=}[u]}
      \end{array} &=
      \begin{array}{c}
        \xymatrix{
          \ar[r]|<>(.5){|}^<>(.5){{}_fB} \ar@{=}[d]
          & \ar@{=}[d]\\
          \ar[r]|<>(.5){|}_<>(.5){{}_fB} & }
      \end{array}
      &
      \begin{array}{c}
        \xymatrix{
          \ar[rr]|<>(.5){|}^<>(.5){B_f} \ar@{}[drr]|\iso \ar@{=}[d] &&
          \ar@{=}[d] \\
          \ar[r]|<>(.5){|}^<>(.5){{}_fB} \ar@{=}[d] \ar@{}[dr]|\Downarrow &
          \ar[r]|<>(.5){|}^<>(.5){U_A} \ar[d]_f \ar@{}[dr]|\Downarrow
          & \ar@{=}[d]\\
          \ar[r]|<>(.5){|}_<>(.5){U_B} &
          \ar[r]|<>(.5){|}_<>(.5){B_f} &\\
          \ar[rr]|<>(.5){|}_<>(.5){B_f} \ar@{}[urr]|\iso \ar@{=}[u] &&
          \ar@{=}[u]}
      \end{array} &=
      \begin{array}{c}
        \xymatrix{
          \ar[r]|<>(.5){|}^<>(.5){B_f} \ar@{=}[d]
          & \ar@{=}[d]\\
          \ar[r]|<>(.5){|}_<>(.5){B_f} & }
      \end{array}\label{eq:frbi-compconj-2}
    \end{align}
  \end{enumerate}
\end{thm}
\begin{proof}
  We first show
  that~\ref{item:frbi-fib}$\Rightarrow$\ref{item:frbi-compconj}.  As
  indicated above, if $(L,R)$ is a fibration we define ${}_fB =
  f^*U_B$ and $B_f = U_Bf^*$, and we let the first two 2-cells
  in~(\ref{eq:frbi-compconj}) be the cartesian 2-cells characterizing
  these two restrictions.  The unique factorizations of $U_f$ through
  these two 2-cells then gives us the second two 2-cells
  in~(\ref{eq:frbi-compconj}) such that the
  equations~(\ref{eq:frbi-compconj-1}) are satisfied by definition.

  We show that the first equation in~(\ref{eq:frbi-compconj-2}) is
  satisfied.  If we compose the left side of this equation with the
  cartesian 2-cell defining ${}_fB$, we obtain
  \begin{align*}
    \begin{array}{c}
      \xymatrix{ \ar[rr]|{|}^{{}_fB} \ar@{=}[d] \ar@{}[drr]|{\iso} &&  \ar@{=}[d]\\
         \ar[r]|{|}^{U_A} \ar@{=}[d] \ar@{}[dr]|{\Downarrow} &
         \ar[r]|{|}^{{}_fB} \ar[d]^f \ar@{}[dr]|{\Downarrow} &
         \ar@{=}[d] \\
         \ar[r]|{|}_{{}_fB} \ar@{=}[d] \ar@{}[drr]|{\iso} &  \ar[r]|{|}_{U_B} &  \ar@{=}[d]\\
         \ar[rr]|{{}_fB} \ar[d]_f \ar@{}[drr]|{\Downarrow} &&  \ar@{=}[d]\\
         \ar[rr]|{|}_{U_B} && 
      }\end{array} \quad =\quad
    \begin{array}{c}
      \xymatrix{ \ar[rr]|{|}^{{}_fB} \ar@{=}[d] \ar@{}[drr]|{\iso} &&  \ar@{=}[d]\\
         \ar[r]|{|}^{U_A} \ar@{=}[d] \ar@{}[dr]|{\Downarrow} &
         \ar[r]|{|}^{{}_fB} \ar[d]^f \ar@{}[dr]|{\Downarrow} &
         \ar@{=}[d] \\
         \ar[r]|{{}_fB} \ar[d]_f \ar@{}[dr]|{\Downarrow} &
         \ar[r]|{|}_{U_B} \ar@{=}[d] & \ar@{=}[d]\\
         \ar[r]|{|}_{U_B} \ar@{=}[d] \ar@{}[drr]|{\iso} &  \ar[r]|{|}_{U_B} &  \ar@{=}[d]\\
         \ar[rr]|{|}_{U_B} && 
      }
    \end{array}
    \quad = \quad
    \begin{array}{c}
      \xymatrix{ \ar[rr]|{|}^{{}_fB} \ar@{=}[d] \ar@{}[drr]|{\iso} &&  \ar@{=}[d]\\
         \ar[r]|{|}^{U_A} \ar[d]_f \ar@{}[dr]|{\Downarrow U_f} &
         \ar[r]|{|}^{{}_fB} \ar[d]^f \ar@{}[dr]|{\Downarrow} &
         \ar@{=}[d] \\
         \ar[r]|{|}_{U_B} \ar@{=}[d] \ar@{}[drr]|{\iso} &  \ar[r]|{|}_{U_B} &  \ar@{=}[d]\\
         \ar[rr]|{|}_{U_B} && 
      }
    \end{array}
    \quad = \quad
    \begin{array}{c}
      \xymatrix{
        \ar[r]|<>(.5){|}^<>(.5){{}_fB} \ar[d]_f \ar@{}[dr]|\Downarrow
        & \ar@{=}[d]\\
        \ar[r]|<>(.5){|}_<>(.5){U_B} & }
    \end{array},
  \end{align*}
  which is once again the cartesian 2-cell defining ${}_fB$.  However,
  we also have
  \[\begin{array}{c}
    \xymatrix{
      \ar[r]|<>(.5){|}^<>(.5){{}_fB} \ar@{=}[d]
      & \ar@{=}[d]\\
      \ar[r]|<>(.5){|}^<>(.5){{}_fB} \ar[d]_f \ar@{}[dr]|\Downarrow
      & \ar@{=}[d]\\
      \ar[r]|<>(.5){|}_<>(.5){U_B} & }
  \end{array},
  \quad = \quad
  \begin{array}{c}
    \xymatrix{
      \ar[r]|<>(.5){|}^<>(.5){{}_fB} \ar[d]_f \ar@{}[dr]|\Downarrow
      & \ar@{=}[d]\\
      \ar[r]|<>(.5){|}_<>(.5){U_B} & }
  \end{array},
  \]
  Thus, the uniqueness of factorizations through cartesian arrows
  implies that the given 2-cell is equal to the identity, as desired.
  This shows the first equation in~(\ref{eq:frbi-compconj-2}); the
  second is analogous.  Thus
  \ref{item:frbi-fib}$\Rightarrow$\ref{item:frbi-compconj}.

  Now assume~\ref{item:frbi-compconj}, and let $M\maps B\hto D$ be a
  1-cell and $f\maps A\to B$ and $g\maps C\to D$ be vertical arrows;
  we claim that the following composite is cartesian:
  \begin{equation}
    \xymatrix{
      \ar[r]|{|}^{{}_fB} \ar[d]_f \ar@{}[dr]|{\Downarrow} &
      \ar[r]|{|}^M \ar@{=}[d] \ar@{}[dr]|{1_M} &
      \ar@{=}[d] \ar[r]|{|}^{D_g} \ar@{}[dr]|{\Downarrow} &
      \ar[d]^g\\
      \ar[r]|{|}_{U_B} & \ar[r]|M & \ar[r]|{|}_{U_D} & \\
      \ar@{=}[u] \ar@{}[urrr]|\iso \ar[rrr]|{|}_M &&&  \ar@{=}[u] .}\label{eq:frbi-compconj-cart}
  \end{equation}
  To show this, suppose that
  \[\xymatrix{\ar[r]|{|}^N\ar[d]_{fh} \ar@{}[dr]|{\Downarrow\alpha} &  \ar[d]^{gk}\\
    \ar[r]|{|}_M & }\]
  is a 2-cell; we must show that it factors uniquely
  through~(\ref{eq:frbi-compconj-cart}).  Consider the composite
  \begin{equation}
    \xymatrix{
      \ar[rrr]|{|}^{N} \ar@{}[drrr]|\iso \ar@{=}[d] &&&
      \ar@{=}[d]\\
      \ar[r]|{|}^{U_A} \ar[d]_h \ar@{}[dr]|{\Downarrow U_h} & \ar[r]|{|}^N \ar[d]^h
      \ar@{}[ddr]|{\Downarrow\alpha} &
      \ar[d]_k\ar[r]|{|}^{U_C} \ar@{}[dr]|{\Downarrow U_k} & \ar[d]^k\\
      \ar[r]|{U_A} \ar@{=}[d] \ar@{}[dr]|{\Downarrow} &
      \ar[d]^f  &
      \ar[d]_g \ar[r]|{U_C} \ar@{}[dr]|{\Downarrow} &
      \ar@{=}[d]\\
      \ar[r]|{|}_{{}_fB} &
      \ar[r]|{|}_M &
      \ar[r]|{|}_{D_g} &.}\label{eq:frbi-compconj-result}
  \end{equation}
  Composing this with~(\ref{eq:frbi-compconj-cart}) and using the
  equations~(\ref{eq:frbi-compconj-1}) on each side, we get $\alpha$
  back again.  Thus,~(\ref{eq:frbi-compconj-result}) gives a
  factorization of $\alpha$ through~(\ref{eq:frbi-compconj-cart}).  To
  prove uniqueness, suppose that we had another factorization
  \begin{equation}
    \begin{array}{c}
      \xymatrix{
        \ar[rrr]|{|}^{N} \ar@{}[drrr]|{\Downarrow\beta} \ar[d]_h &&&
        \ar[d]^k\\
        \ar[r]|{|}^{{}_fB} \ar[d]_f \ar@{}[dr]|{\Downarrow} &
        \ar[r]|M \ar@{=}[d] \ar@{}[dr]|{1_M} &
        \ar@{=}[d] \ar[r]|{|}^{D_g} \ar@{}[dr]|{\Downarrow} &
        \ar[d]^g\\
        \ar[r]|{|}_{U_B} & \ar[r]|M & \ar[r]|{|}_{U_D} & \\
        \ar@{=}[u] \ar@{}[urrr]|\iso \ar[rrr]|{|}_M &&&  \ar@{=}[u].}
    \end{array} =
    \begin{array}{c}
      \xymatrix{\ar[r]|{|}^N\ar[d]_{fh} \ar@{}[dr]|{\Downarrow\alpha} &  \ar[d]^{gk}\\
        \ar[r]|{|}_M & }
    \end{array}\label{eq:frbi-compconj-alt}
  \end{equation}
  Then if we substitute the left-hand side
  of~(\ref{eq:frbi-compconj-alt}) for $\alpha$
  in~(\ref{eq:frbi-compconj-result}) and use the
  equations~(\ref{eq:frbi-compconj-2}) on the left and right, we see
  that everything cancels and we just get $\beta$.  Hence, $\beta$ is
  equal to~(\ref{eq:frbi-compconj-result}), so the factorization is
  unique.  This proves that~(\ref{eq:frbi-compconj-cart}) is
  cartesian, so
  \ref{item:frbi-compconj}$\Rightarrow$\ref{item:frbi-fib}.  The proof
  that \ref{item:frbi-opfib}$\Leftrightarrow$\ref{item:frbi-compconj}
  is exactly dual.
\end{proof}

\begin{defn}\xlabel{def:framed}
  When the equivalent conditions of \xref{thm:framed} are satisfied,
  we say that \bbD\ is a \textbf{framed bicategory}.
\end{defn}

Thus, a framed bicategory has both restrictions and extensions.  By
the construction for
\ref{item:frbi-compconj}$\Rightarrow$\ref{item:frbi-fib}, we see that
in a framed bicategory we have
\begin{equation}
  f^*Mg^* \iso {}_fB \odot M \odot D_g.\label{eq:restriction-by-odot}
\end{equation}
The dual construction for
\ref{item:frbi-compconj}$\Rightarrow$\ref{item:frbi-opfib} shows that
\begin{equation}
  f_!Ng_! \iso {B}_f \odot N \odot {}_gD.\label{eq:extn-by-odot}
\end{equation}
In particular, taking $N=U_B$, we see that addition to
\begin{equation}
  {}_fB \iso f^*U_B \quad\text{and}\quad B_f \iso U_Bf^*,\label{eq:bco-via-restr}
\end{equation}
we have
\begin{equation}
  {}_fB \iso U_Af_! \quad\text{and}\quad B_f \iso f_!U_A.\label{eq:bco-via-extn}
\end{equation}
More specifically, the first two 2-cells in~(\ref{eq:frbi-compconj})
are always cartesian and the second two are opcartesian.  It thus
follows that from the uniqueness of cartesian and opcartesian arrows
that 1-cells ${}_fB$ and $B_f$ equipped with the data of
\xref{thm:framed}\ref{item:frbi-compconj} are unique up to
isomorphism.  In fact, if ${}_fB$ and $\widetilde{{}_fB}$ are two such 1-cells,
the canonical isomorphism ${}_fB\too[\iso] \widetilde{{}_fB}$ is given explicitly
by the composite
\[\begin{array}{c}
  \xymatrix{
    \ar[rr]|<>(.5){|}^<>(.5){{}_fB} \ar@{}[drr]|\iso \ar@{=}[d] &&
    \ar@{=}[d] \\
    \ar[r]|<>(.5){|}^<>(.5){U_A} \ar@{=}[d] \ar@{}[dr]|\Downarrow &
    \ar[r]|<>(.5){|}^<>(.5){{}_fB} \ar[d]_f \ar@{}[dr]|\Downarrow
    & \ar@{=}[d]\\
    \ar[r]|<>(.5){|}_<>(.5){\widetilde{{}_fB}} &
    \ar[r]|<>(.5){|}_<>(.5){U_B} &\\
    \ar[rr]|<>(.5){|}_<>(.5){\widetilde{{}_fB}} \ar@{}[urr]|\iso \ar@{=}[u] &&
    \ar@{=}[u]}
\end{array}\]
The case of $B_f$ is similar.

We can now prove the expected compatibility between base change and
horizontal composition.

\begin{cor}
  In a framed bicategory, we have
  \begin{align*}
    f^*(M\odot N)g^* &\iso f^*M\odot Ng^*\qquad\text{and}\\
    f_!(M\odot N)g_! &\iso f_!M\odot Ng_!.
  \end{align*}
\end{cor}
\begin{proof}
  Use~(\ref{eq:restriction-by-odot}) and~(\ref{eq:extn-by-odot}),
  together with the associativity of $\odot$.
\end{proof}

On the other hand, if coextensions exist, we have a canonical morphism
\begin{equation}
  \label{eq:odot-notpres-coextn}
  f_*M \odot Ng_* \too f_*(M\odot N)g_*
\end{equation}
given by the adjunct of the composite
\[f^*\big(f_*M \odot Ng_*\big)g^* \iso
f^*f_*M \odot Ng_*g^* \too
M\odot N,
\]
but it is rarely an isomorphism.  In general, coextension is often
less well behaved than restriction and extension, which partly
justifies our choice to use a formalism in which it is less natural.

\begin{egs}\xlabel{eg:framed}
  All of the double categories we introduced in
  \S\ref{sec:double-categories} are actually framed bicategories, and
  many of them have coextensions as well.
  \begin{itemize}
  \item In \Mod, if $M$ is an $(A,B)$-bimodule and $f\maps C\to A$,
    $g\maps D\to B$ are ring homomorphisms, then the restriction
    $f^*Mg^*$ is $M$ regarded as a $(C,D)$-bimodule via $f$ and $g$.
    Similarly, $f_!$ is given by extension of scalars and $f_*$ by
    coextension of scalars.  The base change objects ${}_fB$ and $B_f$
    are $B$ regarded as an $(A,B)$-bimodule and $(B,A)$-bimodule,
    respectively, via the map $f$.
  \item In $\Span(\sC)$, restrictions are given by pullback and
    extensions are given by composition.  The base change objects
    ${}_fB$ and $B_f$ for a map $f\maps A\to B$ in \sC\ are the spans
    $A\oot[1_A] A\too[f] B$ and $B\oot[f] A \too[1_A] B$,
    respectively.  These are often known as the \emph{graph} of $f$.
    Coextensions exist when \sC\ is locally cartesian closed.
  \item In \Ex, the base change functors are defined in~\cite[\S{}11.4
    and \S{}12.6]{pht}, and the base change objects are a version of
    the sphere spectrum described in~\cite[\S17.2]{pht}.
  \item In $\Dist(\sV)$, restrictions are given by precomposition, and
    extensions and coextensions are given by left and right Kan
    extension, respectively.  For a \V-functor $f\maps A\to B$, the
    base change objects ${}_fB$ and $B_f$ are the distributors
    $B(-,f-)$ and $B(f-,-)$, respectively.
  \item In \nCob, restriction, extension, and coextension are all
    given by composing a diffeomorphism of $(n-1)$-manifolds with the
    given diffeomorphism onto a collar of the boundary.  The base
    change objects of a diffeomorphism $f\maps A\iso B$ are $f$ and
    its inverse, regarded as thin cobordisms.
  \item In \Adj, restriction and extension are given by composing with
    suitable adjoints.  For example, given $h\maps B\to D$ and
    adjunctions $f_!\maps A \toot B\spam f^*$ and $g_!\maps C \toot
    D\spam g^*$, then a cartesian 2-cell is given by the square
    \[\xymatrix{A \ar[r]|{|}^{g^*hf_!}\ar[d]_{f_!} \drtwocell\omit{\ep} &
      C \ar[d]^{g_!}\\
      B \ar[r]|{|}_h & D}\]
    where $\ep$ is the counit of the adjunction $g_!\adj g^*$.  The
    base change objects for an adjunction $f_!\adj f^*$ are $f_!$ and
    $f^*$, respectively.  Coextensions do not generally exist.
  \end{itemize}
\end{egs}

We also observe that the base change objects are pseudofunctorial.
This is related to, but distinct from, the pseudofunctoriality of the
base change \emph{functors}.  Pseudofunctoriality of base change
functors means that for $A\too[f] B \too[g] C$, we have
\[f^*(g^*(M)) \iso (gf)^*M
\]
coherently, while pseudofunctoriality of base change objects means
that we have
\[{}_fB \odot {}_gC \iso {}_{gf}C
\]
coherently.  However, since base change objects represent all base
change functors, either implies the other.

\begin{prop}\xlabel{bco-psfr}
  If \bbD\ is a framed bicategory with a chosen cleavage, then the
  operation $f\mapsto {}_fB$ defines a pseudofunctor $\bbDz\to\calD$
  which is the identity on objects.  Similarly, the operation
  $f\mapsto B_f$ defines a contravariant pseudofunctor
  $\bbDz\op\to\calD$.
\end{prop}

\section{Duality theory}
\label{sec:duality-theory}

We mentioned in \xref{eg:adj} that the notion of an \emph{adjunction}
can be defined internal to any 2-category.  In fact, the definition
can easily be extended to any bicategory: an adjunction in a
bicategory \calB\ is a pair of 1-cells $F\maps A\hto B$ and $G\maps
B\hto A$ together with 2-cells $\eta\maps U_B\to G\odot F$ and
$\ep\maps F\odot G \to U_A$, satisfying the usual triangle identities
with appropriate associativity and unit isomorphisms inserted.

An internal adjunction is an example of a formal concept which is
useful in both types of bicategories discussed in the introduction,
but its \emph{meaning} is very different in the two cases.  In
\Cat-like bicategories, adjunctions behave much like ordinary adjoint
pairs of functors; in fact, we will use them in this way in
\S\ref{sec:framed-adjunctions}.  In \calMod-like bicategories, on the
other hand, adjunctions encode a notion of \emph{duality}.

In particular, if \sC\ is a monoidal category, considered as a
one-object bicategory, an adjunction in \sC\ is better known as a
\emph{dual pair} in \sC, and one speaks of an object $Y$ as being
\emph{left} or \emph{right dual} to an object $X$; see, for
example,~\cite{may:picard}.  When \sC\ is symmetric monoidal, left
duals and right duals coincide.

\begin{egs}
  When $\sC=\mathbf{Mod}_R$ for a commutative ring $R$, the dualizable
  objects are the finitely generated projectives.  When \sC\ is the
  stable homotopy category, the dualizable objects are the finite CW spectra.
\end{egs}

The terminology of dual pairs was extended in~\cite{pht} to
adjunctions in \calMod-like bicategories, which behave more like dual
pairs in monoidal categories than they do like adjoint pairs of
functors.  Of course, now the distinction between left and right
matters.  Explicitly, we have the following.

\begin{defn}
  A \textbf{dual pair} in a bicategory \calD\ is a pair $(M,N)$, with
  $M\maps A\hto B$, $N\maps B\hto A$, together with `evaluation' and
  `coevaluation' maps 
  \[N\odot M \to U_B
  \qquad\text{and}\qquad
  U_A \to M\odot N
  \]
  satisfying the triangle identities.  We say that $N$ is the
  \textbf{right dual} of $M$ and that $M$ is \textbf{right dualizable},
  and dually.
\end{defn}

The definition of dual pair given in~\cite[16.4.1]{pht} is actually
reversed from ours, although it doesn't look it, because of our
different conventions about which way to write horizontal composition.
But because we also turn around the horizontal 1-cells in all the
examples, the terms `right dualizable' and `left dualizable' refer to
the same actual objects as before.  Our convention has the advantage
that the right dual is also the right adjoint.

Because a dual pair is formally the same as an adjunction, all formal
properties of the latter apply as well to the former.  One example is
the calculus of mates, as defined in \xref{eg:adj}: if $(M,N)$ and
$(P,Q)$ are dual pairs, then there is a natural bijection between
morphisms $M\to P$ and $Q\to N$.

We define a dual pair in a framed bicategory \bbD\ to be just a dual
pair in its underlying horizontal bicategory \calD.  In this case, we
have natural examples coming from the base change objects.

\begin{prop}\xlabel{bco-dual-pairs}
  If $f\maps A\to B$ is a vertical arrow in a framed bicategory \bbD,
  then $({}_fB,B_f)$ is naturally a dual pair.
\end{prop}
\begin{proof}
  Since the base change \emph{functor} $f_!$ is left adjoint to $f^*$,
  we have equivalences
  \[\calD(M \odot {}_fB,N) \eqv \calD(Mf_!,N) \eqv \calD(M,Nf^*) \eqv \calD(M,N \odot B_f).
  \]
  By the bicategorical Yoneda lemma (see, for
  example,~\cite{street:fibi}), which applies to dual pairs just as it
  applies to adjunctions, this implies the desired result.

  Alternately, the unit and counit can be constructed directly from
  the data in \xref{thm:framed}\ref{item:frbi-compconj}; the unit is
  \[\xymatrix{
    \ar[rr]|{|}^{U_A} \ar@{}[drr]|\iso \ar@{=}[d] &&\ar@{=}[d]\\
    \ar[r]|{|}^{U_A}\ar@{=}[d]  \ar@{}[dr]|\Downarrow &
    \ar[r]|{|}^{U_A}\ar[d]^f \ar@{}[dr]|\Downarrow & \ar@{=}[d]\\
    \ar[r]|{|}_{{}_fB} & \ar[r]|{|}_{B_f} & }\]
  and the counit is
  \[\xymatrix{
    \ar[r]|{|}^{B_f} & \ar[r]|{|}^{{}_fB} \ar[d]^f & \\
    \ar[r]|{|}_{U_B}\ar@{=}[u]  \ar@{}[ur]|\Downarrow &
    \ar[r]|{|}_{U_B} \ar@{}[ur]|\Downarrow & \ar@{=}[u]\\
    \ar[rr]|{|}_{U_B} \ar@{}[urr]|\iso \ar@{=}[u] &&\ar@{=}[u].}\]
  Equations~(\ref{eq:frbi-compconj-1}) and~(\ref{eq:frbi-compconj-2})
  are then exactly what is needed to prove the triangle identities.
\end{proof}

In particular, each of the base change objects ${}_fB$ and $B_f$
determines the other up to isomorphism.  Combining these dual pairs
with another general fact about adjunctions in a bicategory, we have
the following generalization of~\cite[17.3.3--17.3.4]{pht}.

\begin{prop}\xlabel{base-change-of-duals}
  Let $(M,N)$ be a dual pair in a framed bicategory with $M\maps A\hto
  B$, $N\maps B\hto A$, and let $f\maps B\to C$ be a vertical arrow.
  Then $(Mf_!,f_!N)$ is also a dual pair.  Similarly, for any $g\maps
  D\to A$, $(g^*M,Ng^*)$ is a dual pair.
\end{prop}
\begin{proof}
  We have $Mf_! \iso M\odot {}_fB$ and $f_!N \iso B_f\odot N$, so the
  result follows from the fact that the composite of adjunctions in a
  bicategory is an adjunction.  The other case is analogous.
\end{proof}

This implies the following generalization of the calculus of mates.

\begin{prop}
  Let $(M,N)$ and $(P,Q)$ be dual pairs in a framed bicategory.  Then
  there is a natural bijection between 2-cells of the following forms:
  \[
  \begin{array}{c}
    \xymatrix{A \ar[r]^M|{|}\ar[d]_f \ar@{}[dr]|{\Downarrow} & B \ar[d]^g\\
      C\ar[r]_P|{|} & D}
  \end{array}
  \qquad\text{and}\qquad
  \begin{array}{c}
    \xymatrix{B \ar[r]^N|{|}\ar[d]_g \ar@{}[dr]|{\Downarrow} & A \ar[d]^f\\
      D\ar[r]_Q|{|} & C}
  \end{array}.
  \]
\end{prop}
\begin{proof}
  A 2-cell of the former shape is equivalent to a globular 2-cell
  $Mg_! \to f^*P$, and a 2-cell of the latter shape is equivalent to a
  globular 2-cell $g_!N\to Qf^*$.  By \xref{base-change-of-duals}, we
  have dual pairs $(Mg_!,g_!N)$ and $(f^*P,Qf^*)$, so the ordinary
  calculus of mates applies.
\end{proof}

\begin{egs}\xlabel{egs:dual-pairs}
  Dual pairs behave significantly differently in many of our examples.
  \begin{enumerate}
  \item If $R$ is a not-necessarily commutative ring, then a right
    $R$-module $M\maps \bbZ\hto R$ is right dualizable in \Mod\ when it
    is finitely generated projective.
  \item The only dual pairs in $\Span(\sC)$ are the base change dual
    pairs.  (This is easy in \Set, and we can then apply the Yoneda
    lemma for arbitrary \sC.)
  \item If $M\maps A\hto B$ is a right dualizable distributor in
    $\Dist(\sV)$, and $B$ satisfies a mild cocompleteness condition
    depending on \sV\ (called `Cauchy completeness'), then $M$ is
    necessarily also of the form ${}_fB$ for some \sV-functor $f\maps
    A\to B$.  When $\sV=\Set$, Cauchy completeness just means that
    every idempotent splits.  When $\sV=\Ab$, it means that
    idempotents split and finite coproducts exist.
    See~\cite[\S5.5]{kelly} for more about Cauchy completion of
    enriched categories.\label{item:dual-pairs-dist}
  \item Dualizable objects in \Ex\ are studied extensively
    in~\cite[Ch.~18]{pht}.
  \end{enumerate}
\end{egs}

\begin{rmk}
  There is also a general notion of \emph{trace} for endomorphisms of
  a dualizable object in a symmetric monoidal category: if $(X,Y)$ is
  a dual pair and $f\maps X\to X$, then the trace of $f$ is the
  composite
  \[I\too[\eta] X\ten Y \too[f\ten 1] X\ten Y \too[\iso] Y\ten X
  \too[\ep] I.
  \]
  Traces were extended to dual pairs in a bicategory
  in~\cite{kate:traces}, by equipping the bicategory with a suitable
  structure, called a \emph{shadow}, to take the place of the symmetry
  isomorphism.  In~\cite{me-kate:traces} we will consider shadows in
  framed bicategories.
\end{rmk}

Duality in symmetric monoidal categories is most interesting when the
monoidal category is closed.  There is also a classical notion of
\emph{closed bicategory}, which means that the composition of 1-cells
has adjoints on both sides:
\[\calB(M\odot N,P) \iso \calB(M,N\rhd P) \iso \calB(N,P\lhd M).\]
Recall that $\calB(M\odot N,P)$ denotes the \emph{set} of globular
2-cells from $M\odot N$ to $P$.  In 2-categorical language, this says
that right Kan extensions and right Kan liftings exist in the
bicategory \calB.

It is proven in~\cite[\S{}16.4]{pht}, extending classical results for
symmetric monoidal categories, that when $M\maps A\hto B$ is right
dualizable, its right dual is always (isomorphic to) the `canonical
dual' $D_rM = M\rhd U_B$; and conversely, whenever the canonical map
$M \odot D_rM \to M\rhd M$ is an isomorphism, then $M$ is right
dualizable.  This can also be stated as the generalization to
bicategories of the fact (see~\cite[X.7]{maclane}) that a functor $G$
has an adjoint when the Kan extension of the identity along $G$ exists
and is preserved by $G$, and in that case the Kan extension gives the
adjoint.

\begin{defn}
  A framed bicategory \bbD\ is \textbf{closed} just when its underlying
  horizontal bicategory \calD\ is closed.
\end{defn}

\begin{egs}
  Many of our examples of framed bicategories are closed.
  \begin{itemize}
  \item $\Mod$ is closed; its hom-objects are given by
    \begin{align*}
      P\lhd M &= \Hom_C(M,P)\\
      N\rhd P &= \Hom_A(N,P).
    \end{align*}
  \item As long as \V\ is closed and complete, then $\Dist(\V)$ is
    closed; its hom-objects are given by the cotensor product of
    distributors (the end construction).
  \item $\Span(\sC)$ is closed precisely when $\sC$ is locally
    cartesian closed.
  \item \Ex\ is also closed.  This is proven in~\cite[\S17.1]{pht}; we
    will describe the general method of proof in
    \S\ref{sec:monfib-to-framed}.
  \end{itemize}
\end{egs}

\begin{rmk}
  A monoidal category is closed (on both sides) just when its
  corresponding vertically trivial framed bicategory is closed.  On
  the other hand, if a monoidal category is symmetric, then the left
  and right internal-homs are isomorphic.  In \S\ref{sec:involutions}
  we will prove an analogue of this fact for framed bicategories
  equipped with an `involution', which includes all of our examples.
\end{rmk}

It is not surprising that there is some relationship between
closedness and base change.

\begin{prop}\xlabel{closed-bco}
  Let \bbD\ be a closed framed bicategory.  Then for any $f:A\vto C$,
  $g:B\vto D$, and $M:C\hto D$, we have
  \begin{equation*}
    \begin{split}
      f^*Mg^* &\iso ({}_gD\rhd M)\lhd C_f \\
      &\iso {}_gD\rhd (M\lhd C_f)
    \end{split}
  \end{equation*}
  and in particular
  \begin{equation*}
    \begin{split}
      {}_fC &\iso C\lhd C_f \\
      D_g   &\iso {}_gD\rhd D
    \end{split}
  \end{equation*}
\end{prop}
\begin{proof}
  Straightforward adjunction arguments.
\end{proof}

Note that this implies, by uniqueness of adjoints, that the
restriction functor $f^*$ can also be described as $f^*N = N \lhd
C_f$.  Of course there are corresponding versions for composing on the
other side.

Moreover, if coextensions exist, then uniqueness of adjoints also
implies that we have
\begin{equation}
  f_*M \iso M \lhd {}_fC.\label{eq:coext-via-closed}
\end{equation}
Conversely, if \bbD\ is closed, then~(\ref{eq:coext-via-closed})
defines a right adjoint to $f^*$; thus coextensions exist in any
closed framed bicategory, and also have a natural description in terms
of the base change objects.

\section{The 2-category of framed bicategories}
\label{sec:2cat-frbi}

We now introduce the morphisms between framed bicategories.  To begin
with, it is easy to define morphisms of double categories by analogy
with monoidal categories.

\begin{defn}
  Let \bbD\ and \bbE\ be double categories.  A \textbf{lax double
    functor} $F\maps \bbD\to \bbE$ consists of the following.
  \begin{itemize}
  \item Functors $F_0\maps \bbD_0 \to \bbE_0$ and $F_1\maps \bbD_1 \to
    \bbE_1$ such that $L\circ F_1 = F_0\circ L$ and $R\circ F_1 =
    F_0\circ R$.
  \item Natural transformations $F_\odot\maps F_1M \odot F_1N \to
    F_1(M\odot N)$ and $F_U\maps U_{F_0 A} \to F_1(U_A)$, whose
    components are globular, and which satisfy the usual coherence
    axioms for a lax monoidal functor or 2-functor (see, for
    example,~\cite[\S{}XI.2]{maclane}).
  \end{itemize}
  Dually, we have the definition of an \textbf{oplax double functor},
  for which $F_\odot$ and $F_U$ go in the opposite direction.  A
  \textbf{strong double functor} is a lax double functor for which
  $F_\odot$ and $F_U$ are (globular) isomorphisms.  If just $F_U$ is
  an isomorphism, we say that $F$ is \textbf{normal}.
\end{defn}

We occasionally abuse notation by writing just $F$ for either $F_0$ or
$F_1$.  Observe that a lax double functor preserves vertical
composition and identities strictly, but preserves horizontal
composition and identities only up to constraints.  Like the
constraints $\fa,\fl,\fr$ for a double category, the maps $F_\odot$
and $F_U$ are globular, but must be natural with respect to all
2-cells, not only globular ones.

If \bbD\ and \bbE\ are just monoidal categories, then a double functor
$F\maps \bbD\to\bbE$ is the same as a monoidal functor (of whichever
sort).  The terms `lax', `oplax', and `strong' are chosen to
generalize this situation; some authors refer to strong double
functors as \emph{pseudo double functors}.  Since the monoidal
functors which arise in practice are most frequently lax, many authors
refer to these simply as `monoidal functors'.  It is also true for
framed bicategories that the lax morphisms are often those of most
interest, but we will always keep the adjectives for clarity.

\begin{eg}\xlabel{eg:monfr->dblfr}
  Let $F\maps \sC\to\sD$ be a lax monoidal functor, where the monoidal
  categories \sC\ and \sD\ both have coequalizers preserved by $\ten$.
  Then it is well known that $F$ preserves monoids, monoid
  homomorphisms, bimodules, and equivariant maps.  Moreover, if
  $M\maps A\hto B$ and $N\maps B\hto C$ are bimodules in \sC, so that
  their tensor product is the coequalizer
  \[\xymatrix{
    M\ten B\ten N \ar@<1mm>[r]\ar@<-1mm>[r] &
    M\ten N \ar[r] &
    M\odot N
  }\]
  then we have the commutative diagram
  \begin{equation}
    \xymatrix{
      FM\ten FB\ten FN \ar@<1mm>[r]\ar@<-1mm>[r] \ar[d] &
      FM\ten FN \ar[r] \ar[d] &
      FM\odot FN \ar@{.>}[d] \\
      F(M\ten B\ten N) \ar@<1mm>[r]\ar@<-1mm>[r] &
      F(M\ten N) \ar[r] &
      F(M\odot N)
    }\label{eq:mod-lax-cmp}
  \end{equation}
  in which the top diagram is the coequalizer defining the tensor
  product of bimodules in \sD, and hence the dotted map is induced.
  Moreover, since $U_A$ in $\Mod(\sC)$ is just $A$ regarded as an
  $(A,A)$-bimodule, we have $F(U_A)\iso U_{FA}$.  It is
  straightforward to check that this isomorphism and the dotted map
  in~\eqref{eq:mod-lax-cmp} are the data for a normal lax double
  functor $\Mod(F)\maps \Mod(\sC)\to\Mod(\sD)$.

  Note that $F$ does not need to preserve coequalizers, so the bottom
  row of~\eqref{eq:mod-lax-cmp} need not be a coequalizer diagram.
  However, if $F$ does preserve coequalizers, and is moreover a strong
  monoidal functor, so that the left and middle vertical maps are
  isomorphisms, then so is the right vertical map; hence $\Mod(F)$ is
  a strong double functor in this case.

  In particular, if $\sC=\mathbf{Mod}_R$ and $\sD=\mathbf{Mod}_S$ for
  commutative rings $R$ and $S$ and $f\maps R\to S$ is a homomorphism
  of commutative rings, then the extension-of-scalars functor
  $f_!\maps \mathbf{Mod}_R \to \mathbf{Mod}_S$ is strong monoidal and
  preserves coequalizers, hence induces a strong double functor.  The
  restriction-of-scalars functor $f^*\maps \mathbf{Mod}_S\to
  \mathbf{Mod}_R$, on the other hand, is only lax monoidal, and hence
  induces a normal lax double functor.
\end{eg}

\begin{eg}\xlabel{eg:lex->dblfr}
  Let $F\maps \sC\to\sD$ be any functor between two categories with
  pullbacks.  Then we have an induced normal \emph{oplax} double
  functor $\Span(F)\maps \Span(\sC)\to\Span(\sD)$.  If $F$ preserves
  pullbacks, then $\Span(F)$ is strong.
\end{eg}

Now suppose that \bbD\ and \bbE\ are framed bicategories.  Since the
characterization of base change objects in
\xref{thm:framed}\ref{item:frbi-compconj} only involves horizontal
composition with units, any \emph{normal} lax (or oplax) framed
functor will preserve base change objects up to isomorphism; that is,
$F({}_fB) \iso {}_{Ff}(FB)$.  If it is \emph{strong}, then it will
also preserve restrictions and extensions, since we have $f^*Mg^*\iso
{}_fB\odot M \odot D_g$ and similarly.

More generally, any lax or oplax double functor $F\maps \bbD\to\bbE$
between framed bicategories automatically induces comparison 2-cells
such as
\begin{align}
  (Ff)_!(FM) &\too F(f_!M) \qquad\text{and}\label{eq:lax-extn-cmp}\\
  F(f^*N) &\too (Ff)^*(FN),\label{eq:lax-restr-cmp}
\end{align}
by unique factorization through cartesian and opcartesian arrows.  As
remarked in \S\ref{sec:fibrations}, the first of these goes in the
`lax direction' while the second goes in the `oplax direction'.  Thus,
for the whole functor to deserve the adjective `lax', the second of
these must be an isomorphism, so that it has an inverse which goes in
the lax direction.  This happens just when $F$ preserves cartesian
2-cells.

However, it turns out that this is automatic: \emph{any} lax double
functor between framed bicategories preserves cartesian 2-cells, so
that~(\ref{eq:lax-restr-cmp}) is always an isomorphism when $F$ is
lax.  Dually, any oplax double functor preserves opcartesian 2-cells,
so that~(\ref{eq:lax-extn-cmp}) is an isomorphism when $F$ is oplax.

To prove this, we first observe that for any lax double functor
$F\maps \bbD\to\bbE$ and any arrow $f\maps A\to B$ in \bbD, we have
the following diagram of 2-cells in \bbE:
\begin{equation}
  \xymatrix{U_{FA} \ar[r]^{F_U}\ar[d]_{\mathrm{opcart}} &
    F(U_A) \ar[d]^{F(\mathrm{opcart})}\\
    {}_{Ff}(FB)\ar@{.>}[r]^{{}_fF} \ar[d]_{\mathrm{cart}} &
    F({}_fB) \ar[d]\ar[d]^{F(\mathrm{cart})}\\
    U_{FB}\ar[r]_{F_U} & F(U_B).}\label{eq:lax-cart-defn}
\end{equation}
The dotted arrow, given by unique factorization through the
opcartesian one, is the special case of~(\ref{eq:lax-extn-cmp}) when
$M=U_A$; we denote it by ${}_fF$.  The upper square
in~(\ref{eq:lax-cart-defn}) commutes by definition, and the lower
square also commutes by uniqueness of the factorization.  Similarly,
we have a 2-cell $(FB)_{Ff} \too[F_f] F(B_f)$.

If $F$ is oplax instead, we have 2-cells in the other direction.  If
$F$ is strong (or even just normal), the transformations exist in both
directions, are inverse isomorphisms, and each is the mate of the
inverse of the other.

\begin{prop}\xlabel{dbl-functor->framed}
  Any lax double functor between framed bicategories preserves
  cartesian 2-cells, and any oplax double functor between framed
  bicategories preserves opcartesian 2-cells.
\end{prop}
\begin{proof}
  Let $M\maps B\hto D$ in \bbD\ and $f\maps A\to B$, $g\maps C\to D$.
  Then the following composite is cartesian in \bbD:
  \begin{equation}
    \xymatrix{A \ar[r]|{|}^{{}_{f}B}\ar[d]_{f} \ar@{}[dr]|{\mathrm{cart}} &
      B \ar@{=}[d] \ar[r]|{|}^{M} &
      D \ar@{=}[d] \ar[r]|{|}^{D_{g}} \ar@{}[dr]|{\mathrm{cart}} &
      C \ar[d]^{g}\\
      B \ar[r]|{|}_{U_B} &
      B \ar[r]|{|}_{M} &
      D \ar[r]|{|}_{U_D}  &
      D}\label{eq:lax-cart-source}
  \end{equation}
  and the following composite is cartesian in \bbE:
  \begin{equation}
    \xymatrix{FA \ar[r]|{|}^{{}_{Ff}(FB)}\ar[d]_{Ff} \ar@{}[dr]|{\mathrm{cart}} &
      FB \ar@{=}[d] \ar[r]|{|}^{FM} &
      FD \ar@{=}[d] \ar[r]|{|}^{(FD)_{Fg}} \ar@{}[dr]|{\mathrm{cart}} &
      FC \ar[d]^{Fg}\\
      FB \ar[r]|{|}_{U_{FB}} &
      FB \ar[r]|{|}_{FM} &
      FD \ar[r]|{|}_{U_{FD}}  &
      FD.}\label{eq:lax-cart-dest}
  \end{equation}
  Applying $F$ to~(\ref{eq:lax-cart-source}) and factoring the result
  through~(\ref{eq:lax-cart-dest}), we obtain a comparison map
  \begin{equation}
    F({}_fB \odot M \odot D_g) \too {}_{Ff}(FB)\odot FM \odot (FD)_{Fg},\label{eq:lax-cart-cmp}
  \end{equation}
  which we want to show to be an isomorphism.  We have an obvious
  candidate for its inverse, namely the following composite.
  \begin{equation}
    \label{eq:lax-cart-inverse}
    \xymatrix{\ar[r]|{|}^{{}_{Ff}(FB)}\ar@{=}[d] \ar@{}[dr]|{{}_fF\Downarrow} &
      \ar[r]|{|}^{FM} \ar@{=}[d] &
      \ar[r]|{|}^{(FD)_{Fg}} \ar@{=}[d]  \ar@{}[dr]|{\Downarrow F_f} &
      \ar@{=}[d] \\
      \ar[r]|{|}_{F({}_fB)} \ar@{=}[d] \ar@{}[drrr]|{\Downarrow F_\odot} &
      \ar[r]|{|}^{FM} &
      \ar[r]|{|}_{F(D_g)} &
      \ar@{=}[d]\\
      \ar[rrr]|{|}_{F({}_fB \odot M \odot D_g)} &&&}
  \end{equation}
  Consider first the composite of~(\ref{eq:lax-cart-inverse}) followed
  by~(\ref{eq:lax-cart-cmp}):
  \[\xymatrix{\ar[r]|{|}^{{}_{Ff}(FB)}\ar@{=}[d] \ar@{}[dr]|{{}_fF\Downarrow} &
    \ar[r]|{|}^{FM} \ar@{=}[d] &
    \ar[r]|{|}^{(FD)_{Fg}} \ar@{=}[d]  \ar@{}[dr]|{\Downarrow F_f} &
    \ar@{=}[d] \\
    \ar[r]|{|}_{F({}_fB)} \ar@{=}[d] \ar@{}[drrr]|{\Downarrow F_\odot} &
    \ar[r]|{|}^{FM} &
    \ar[r]|{|}_{F(D_g)} &
    \ar@{=}[d]\\
    \ar[rrr]|{F({}_fB \odot M \odot D_g)} \ar@{=}[d] \ar@{}[drrr]|\Downarrow &&& \ar@{=}[d]\\
    \ar[r]|{|}_{{}_{Ff}(FB)} &
    \ar[r]|{|}_{FM} &
    \ar[r]|{|}_{(FD)_{Fg}} &
  }\]
  If we postcompose this with~(\ref{eq:lax-cart-dest}), then by
  definition of~(\ref{eq:lax-cart-cmp}), we obtain
  \[\xymatrix{\ar[r]|{|}^{{}_{Ff}(FB)}\ar@{=}[d] \ar@{}[dr]|{{}_fF\Downarrow} &
    \ar[r]|{|}^{FM} \ar@{=}[d] &
    \ar[r]|{|}^{(FD)_{Fg}} \ar@{=}[d]  \ar@{}[dr]|{\Downarrow F_f} &
    \ar@{=}[d] \\
    \ar[r]|{|}_{F({}_fB)} \ar@{=}[d] \ar@{}[drrr]|{\Downarrow F_\odot} &
    \ar[r]|{|}^{FM} &
    \ar[r]|{|}_{F(D_g)} &
    \ar@{=}[d]\\
    \ar[rrr]|{F({}_fB \odot M \odot D_g)} \ar[d]_{Ff} \ar@{}[drrr]|{F(\mathrm{cart})}
    &&& \ar[d]^{Fg}\\
    \ar[rrr]|{|}_{FM} &&&.}\]
  By naturality of the lax constraint for $F$, this is equal to
  \[\xymatrix@C=3pc{\ar[r]|{|}^{{}_{Ff}(FB)}\ar@{=}[d] \ar@{}[dr]|{{}_fF\Downarrow} &
    \ar[r]|{|}^{FM} \ar@{=}[d] &
    \ar[r]|{|}^{(FD)_{Fg}} \ar@{=}[d]  \ar@{}[dr]|{\Downarrow F_f} &
    \ar@{=}[d] \\
    \ar[r]|{F({}_fB)} \ar[d]_{Ff} \ar@{}[dr]|{F(\mathrm{cart})} &
    \ar[r]|{|}^{FM}\ar@{=}[d] &
    \ar[r]|{F(D_g)} \ar@{=}[d] \ar@{}[dr]|{F(\mathrm{cart})} &
    \ar[d]^{Fg}\\
    \ar[r]|{|}_{F(U_B)} \ar@{=}[d] \ar@{}[drrr]|{\Downarrow F_\odot} &
    \ar[r]|{|}^{FM} &
    \ar[r]|{|}_{F(U_D)}  &
    \ar@{=}[d]\\
    \ar[rrr]|{|}_{FM} &&& 
  }\]
  Because the lower square in~(\ref{eq:lax-cart-defn}) commutes, this
  is equal to
  \[\xymatrix@C=3pc{\ar[r]|{|}^{{}_{Ff}(FB)}\ar[d]_{Ff} \ar@{}[dr]|{\mathrm{cart}} &
    \ar[r]|{|}^{FM} \ar@{=}[d] &
    \ar[r]|{|}^{(FD)_{Fg}} \ar@{=}[d] \ar@{}[dr]|{\mathrm{cart}} &
    \ar[d]^{Fg} \\
    \ar[r]|{U_{FB}} \ar@{=}[d] \ar@{}[dr]|{\Downarrow F_U} &
    \ar[r]|{|}^{FM}\ar@{=}[d] &
    \ar[r]|{U_{FD}} \ar@{=}[d] \ar@{}[dr]|{\Downarrow F_U} &
    \ar@{=}[d]\\
    \ar[r]|{|}_{F(U_B)} \ar@{=}[d] \ar@{}[drrr]|{\Downarrow F_\odot} &
    \ar[r]|{|}^{FM} &
    \ar[r]|{|}_{F(U_D)}  &
    \ar@{=}[d]\\
    \ar[rrr]|{|}_{FM} &&&
  }\]
  which is equal to~(\ref{eq:lax-cart-dest}), by the coherence axioms
  for $F$.  Thus, by unique factorization
  through~(\ref{eq:lax-cart-dest}), we conclude
  that~(\ref{eq:lax-cart-inverse}) followed by~(\ref{eq:lax-cart-cmp})
  is the identity.

  Now consider the composite of~(\ref{eq:lax-cart-cmp}) followed
  by~(\ref{eq:lax-cart-inverse}).  By the construction of
  factorizations in \xref{thm:framed},~(\ref{eq:lax-cart-cmp}) can be
  computed by composing horizontally with opcartesian 2-cells; thus
  our desired composite is
  \[\xymatrix@C=4pc{
    \ar[r]|{|}^{U_{FA}} \ar@{=}[d] \ar@{}[dr]|{\mathrm{opcart}} &
    \ar[rr]|{|}^{F({}_fB\odot M\odot D_g)} \ar[d]^{Ff} \ar@{}[drr]|{F(\mathrm{cart})} &&
    \ar[r]|{|}^{U_{FC}} \ar[d]_{Fg} \ar@{}[dr]|{\mathrm{opcart}} &
    \ar@{=}[d] \\
    \ar[r]|{{}_{Ff}(FB)}\ar@{=}[d] \ar@{}[dr]|{{}_fF\Downarrow} &
    \ar[rr]|{FM} \ar@{=}[d] &&
    \ar[r]|{(FD)_{Fg}} \ar@{=}[d] \ar@{}[dr]|{\Downarrow F_f} &
    \ar@{=}[d] \\
    \ar[r]|{|}_{F({}_fB)} \ar@{=}[d] \ar@{}[drrrr]|{\Downarrow F_\odot} &
    \ar[rr]|{|}^{FM} &&
    \ar[r]|{|}_{F(D_g)} &
    \ar@{=}[d]\\
    \ar[rrrr]|{|}_{F({}_fB\odot M\odot D_g)} & && &.
  }\]
  By definition of ${}_fF$ and $F_f$, this is equal to
  \[\xymatrix@C=3pc{
    \ar[r]|{|}^{U_{FA}} \ar@{=}[d] \ar@{}[dr]|{\Downarrow F_U} &
    \ar[rr]|{|}^{F({}_fB\odot M\odot D_g)} \ar@{=}[d]  &&
    \ar[r]|{|}^{U_{FC}} \ar@{=}[d] \ar@{}[dr]|{\Downarrow F_U} &
    \ar@{=}[d] \\
    \ar[r]|{F(U_A)}\ar@{=}[d] \ar@{}[dr]|{F(\mathrm{opcart})} &
    \ar[rr]|{|}^{F({}_fB\odot M\odot D_g)} \ar[d] \ar@{}[drr]|{F(\mathrm{cart})} &&
    \ar[r]|{F(U_C)} \ar[d] \ar@{}[dr]|{F(\mathrm{opcart})} &
    \ar@{=}[d] \\
    \ar[r]|{|}_{F({}_fB)} \ar@{=}[d] \ar@{}[drrrr]|{\Downarrow F_\odot} &
    \ar[rr]|{FM} &&
    \ar[r]|{|}_{F(D_g)} &
    \ar@{=}[d]\\
    \ar[rrrr]|{|}_{F({}_fB\odot M\odot D_g)} & && &.
  }\]
  By naturality of $F_\odot$, this is equal to
  \begin{equation}
    \xymatrix@C=3pc{
      \ar[r]|{|}^{U_{FA}} \ar@{=}[d] \ar@{}[dr]|{\Downarrow F_U} &
      \ar[rr]|{|}^{F({}_fB\odot M\odot D_g)} \ar@{=}[d]  &&
      \ar[r]|{|}^{U_{FC}} \ar@{=}[d] \ar@{}[dr]|{\Downarrow F_U} &
      \ar@{=}[d] \\
      \ar[r]|{|}_{F(U_A)}\ar@{=}[d] \ar@{}[drrrr]|{\Downarrow F_\odot} &
      \ar[rr]|{|}^{F({}_fB\odot M\odot D_g)} &&
      \ar[r]|{|}_{F(U_C)} &
      \ar@{=}[d] \\
      \ar[rrrr]|{F({}_fB\odot M\odot D_g)} \ar@{=}[d]
      \ar@{}[drrrr]|{F(\mathrm{stuff})} & && & \ar@{=}[d]\\
      \ar[rrrr]|{|}_{F({}_fB\odot M\odot D_g)} & && &
    }\label{eq:dblfr-fr-result}
  \end{equation}
  where `$\mathrm{stuff}$' is the composite
  \[\xymatrix{\ar[r]|{|}^{U_A}\ar@{=}[d] \ar@{}[dr]|{\mathrm{opcart}} &
    \ar[r]|{|}^{{}_fB} \ar[d]|f \ar@{}[dr]|{\mathrm{cart}} &
    \ar[r]|{|}^M \ar@{=}[d] &
    \ar[r]|{|}^{D_g} \ar@{=}[d] \ar@{}[dr]|{\mathrm{cart}} &
    \ar[r]|{|}^{U_C} \ar[d]|g \ar@{}[dr]|{\mathrm{opcart}} &
    \ar@{=}[d]\\
    \ar[r]|{|}_{{}_fB} &
    \ar[r]|{|}_{U_B} &
    \ar[r]|{|}_M &
    \ar[r]|{|}_{U_D} &
    \ar[r]|{|}_{D_g} & }\]
  which is equal (modulo constraints) to the identity on ${}_fB \odot
  M \odot D_g$.  Thus, applying the coherence axioms for $F$
  again,~\eqref{eq:dblfr-fr-result} reduces to the identity of
  $F({}_fB \odot M \odot D_g)$.
  Therefore,~(\ref{eq:lax-cart-inverse}) is a two-sided inverse
  for~(\ref{eq:lax-cart-cmp}), so the latter is an isomorphism; hence
  $F$ preserves cartesian 2-cells.  The oplax case is dual.
\end{proof}

Here we see again the advantage of using fibrations rather than
introducing base change functors explicitly: since fibrations are
`non-algebraic', all their constraints and coherence come for free.
This leads us to the following definition.

\begin{defn}\xlabel{def:framed-frs}
  A \textbf{lax framed functor} is a lax double functor between framed
  bicategories.  Similarly, an \textbf{oplax} or \textbf{strong framed
    functor} is a double functor of the appropriate type between
  framed bicategories.
\end{defn}

We observed in \S\ref{sec:introduction} that while 2-functors give a
good notion of morphism between both sorts of bicategories, the right
notion of transformation for \calMod-like bicategories is rather
murkier.  Once we include the vertical arrows to get a framed
bicategory, however, it becomes much clearer what the transformations
should be.

\begin{defn}
  A \textbf{double transformation} between two lax double functors
  $\alpha: F\to G:\bbD\to\bbE$ consists of natural transformations
  $\alpha_0\maps F_0\to G_0$ and $\alpha_1\maps F_1\to G_1$ (both
  usually written as $\alpha$), such that $L(\alpha_{M}) =
  \alpha_{LM}$ and $R(\alpha_{M}) = \alpha_{RM}$, and such that
  \[\begin{array}{c}
    \xymatrix{
      FA \ar@{=}[d] \ar[r]|{|}^{FM}
      \ar@{}[drr]|{\Downarrow F_\odot} &
      FB \ar[r]|{|}^{FN} &
      FC \ar@{=}[d]\\
      FA \ar[rr]|{F(M\odot N)} \ar[d]_{\alpha_A}
      \ar@{}[drr]|{\Downarrow \alpha_{M\odot N}} &&
      FC \ar[d]^{\alpha_C}\\
      GA \ar[rr]|{|}_{G(M\odot N)} && GC
    }\end{array} =
  \begin{array}{c}
    \xymatrix{
      FA \ar[d]_{\alpha_A} \ar@{}[dr]|{\Downarrow \alpha_M} \ar[r]|{|}^{FM} &
      FB \ar[d]|{\alpha_B} \ar@{}[dr]|{\Downarrow \alpha_N} \ar[r]|{|}^{FN} &
      FC \ar[d]^{\alpha_C}\\
      GA \ar@{=}[d] \ar[r]|{|}_{GM} \ar@{}[drr]|{\Downarrow G_\odot} &
      GB \ar[r]|{|}_{GN} &
      GC \ar@{=}[d]\\
      GA \ar[rr]|{|}_{G(M\odot N)} && GC
    }\end{array}\]
  and
  \[\begin{array}{c}
    \xymatrix{
      FA \ar[rr]|{|}^{U_{FA}} \ar@{=}[d]
      \ar@{}[drr]|{\Downarrow F_0} &&
      FA \ar@{=}[d]\\
      FA \ar[rr]|{F(U_A)} \ar[d]_{\alpha_A}
      \ar@{}[drr]|{\Downarrow \alpha_{U_A}} &&
      FA \ar[d]^{\alpha_A}\\
      GA \ar[rr]|{|}_{G(U_A)} && GA
    }
  \end{array} =
  \begin{array}{c}
    \xymatrix{
      FA \ar[rr]|{|}^{U_{FA}} \ar[d]_{\alpha_A}
      \ar@{}[drr]|{\Downarrow U_{\alpha_A}} &&
      FA \ar[d]^{\alpha_A}\\
      GA \ar[rr]|{U_{GA}} \ar@{=}[d]
      \ar@{}[drr]|{\Downarrow F_0} &&
      GA \ar@{=}[d]\\
      GA \ar[rr]|{|}_{G(U_A)} && GA
    }.
  \end{array}\]
\end{defn}

The framed version of this definition requires no modification at all.

\begin{defn}\xlabel{def:framed-transf}
  A \textbf{framed transformation} between two lax framed functors is
  simply a double transformation between their underlying lax double
  functors.
\end{defn}

We leave it to the reader to define transformations between oplax
functors.  In the case of ordinary bicategories, there is also a
notion of `modification', or morphism between transformations, but
with framed bicategories we usually have no need for these.  Thus, the
framed bicategories, framed functors, and framed transformations form
a \Cat-like bicategory, which is in fact a strict 2-category.

\begin{prop}\xlabel{dbl-2cat}
  Small framed bicategories, lax framed functors, and framed
  transformations form a strict 2-category \FrBil.  If we restrict to
  strong framed functors, we obtain a 2-category \FrBi, and if we use
  oplax framed functors instead, we obtain a 2-category \FrBic.
\end{prop}

Of course, double categories, double functors, and double
transformations also form larger 2-categories \Dbll, \Dbl, and \Dblc.

\begin{rmk}
  Recall that we can regard a monoidal category as a framed bicategory
  whose vertical category is trivial, and that the framed functors
  between vertically trivial framed bicategories are precisely the
  monoidal functors (whether lax, oplax, or strong).  It is easy to
  check that framed transformations are also the same as monoidal
  transformations; thus \MonCat\ is equivalent to a full
  sub-2-category of \FrBi.

  This is to be contrasted with the situation for ordinary `unframed'
  bicategories.  We can also consider monoidal categories to be
  bicategories with just one 0-cell, and 2-functors between such
  bicategories do also correspond to monoidal functors, but most
  transformations between such 2-functors do not give rise to anything
  resembling a monoidal transformation; see~\cite{cg:degeneracy-i}.
  Thus, framed bicategories are a better generalization of monoidal
  categories than ordinary bicategories are.
\end{rmk}

\begin{eg}\xlabel{eg:montr->dbltr}
  Let \sC\ and \sD\ be monoidal categories with coequalizers preserved
  by $\ten$, and let $\alpha\maps F\Rightarrow G \maps \sC\to\sD$ be a
  monoidal natural transformation between lax monoidal functors.  We
  have already seen that $F$ and $G$ give rise to lax framed functors.
  Moreover, the fact that $\alpha$ is a monoidal transformation
  implies that if $A$ is a monoid in \sC, $\alpha_A\maps FA\to GA$ is
  a monoid homomorphism in \sD, and similarly for bimodules.
  Therefore, we have an induced framed transformation
  \[\Mod(\alpha)\maps \Mod(F)\to\Mod(G).
  \]
  This makes $\Mod(-)$ into a strict 2-functor.  Its domain is the
  2-category of monoidal categories with coequalizers preserved by
  $\ten$, lax monoidal functors, and monoidal transformations, and its
  codomain is \FrBil.  If we restrict the domain to strong monoidal
  functors which preserve coequalizers, the image lies in \FrBi.
\end{eg}

\begin{eg}\xlabel{eg:cart->dbltr}
  Let $\sC,\sD$ be categories with pullbacks and $\alpha\maps
  F\Rightarrow G\maps \sC\to\sD$ a natural transformation.  Then
  $\alpha$ induces a framed transformation
  \[\Span(\alpha)\maps \Span(F)\to\Span(G)
  \]
  in an obvious way.  This makes $\Span$ into a strict 2-functor from the
  2-category of categories with pullbacks, all functors, and all
  natural transformations, to $\FrBic$.  If we restrict the domain to
  functors which preserve pullbacks, the image lies in \FrBi.
\end{eg}

\begin{rmk}
  It is easy to see that any framed functor induces a 2-functor of the
  appropriate type between horizontal bicategories, but the situation
  for framed transformations is less clear.  We will consider this
  further in appendix~\ref{sec:frbi-vs-bicat}.
\end{rmk}

\section{Framed equivalences}
\label{sec:framed-equivalences}

All the usual notions of 2-category theory apply to the study of
framed bicategories via the 2-categories \FrBil, \FrBi, and \FrBic,
and generally reduce to elementary notions when expressed explicitly.
Since, as remarked above, the lax framed functors are often those of
most interest, we work most frequently in \FrBil, but analogous
results are always true for the other two cases.

One important 2-categorical notion is that of \emph{internal
  equivalence}.  This is defined to be a pair of morphisms $F\maps
D\to E$ and $G\maps E\to D$ with 2-cell isomorphisms $FG\iso \Id$ and
$GF\iso \Id$.  The notion of equivalence for framed bicategories we
obtain in this way solves another of the problems raised in
\S\ref{sec:introduction}.

\begin{defn}
  A \textbf{framed equivalence} is an internal equivalence in \FrBil.
\end{defn}

Thus, a framed equivalence consists of lax framed functors
$F:\bbD\rightleftarrows\bbE:G$ with framed natural isomorphisms
$\eta:\Id_\bbD\iso GF$ and $\varepsilon:FG\iso\Id_\bbE$.  It might
seem strange not to require $F$ and $G$ to be \emph{strong} framed
functors in this definition, but in fact this is automatic.

\begin{prop}\xlabel{eqv-strong}
  In a framed equivalence as above, $F$ and $G$ are automatically
  strong framed functors (hence give an equivalence in \FrBi).
\end{prop}

We will prove this in the next section as \xref{eqv-strong-2}.

Since strict 2-functors preserve internal equivalences, our
2-functorial ways of constructing framed bicategories give us a ready
supply of framed equivalences.  For example, any monoidal equivalence
$\sC\eqv \sD$ of monoidal categories with coequalizers preserved by
$\ten$ induces a framed equivalence $\Mod(\sC)\eqv \Mod(\sD)$.
Similarly, any equivalence of categories with pullbacks induces a
framed equivalence between framed bicategories of spans.

As for ordinary categories, we can characterize the framed functors
which are equivalences as those which are `full, faithful, and
essentially surjective'.  First we introduce the terminology,
beginning with double categories.  Recall that we write
$_f\bbD_g(M,N)$ for the set of 2-cells of the form
\[\xymatrix{ A \ar[r]|{|}^{M} \ar[d]_f \ar@{}[dr]|{\Downarrow\alpha}&
  B\ar[d]^g\\
  C \ar[r]|{|}_N & D.}\]

\begin{defn}\xlabel{def:ff}
  A lax or oplax double functor is \textbf{full} (resp.\
  \textbf{faithful}) if it is full (resp.\ faithful) on vertical
  categories and each map
  \begin{equation}
    F\maps {}_g\bbD_f(M,N)\too {}_{Fg}\bbE_{Ff}(FM,FN)\label{eq:ff-2cell-map}
  \end{equation}
  is surjective (resp.\ injective).
\end{defn}

In the case of a \emph{framed} functor, however, the notions simplify
somewhat.

\begin{prop}\xlabel{framed-ff}
  A lax or oplax \emph{framed} functor $F\maps \bbD\to\bbE$ is full
  (resp.\ faithful) in the sense of \xref{def:ff} if and only if it is
  full (resp.\ faithful) on vertical categories and each functor
  $\calD(A,B)\to\calE(FA,FB)$ is full (resp.\ faithful).
\end{prop}
\begin{proof}
  \xref{def:ff} clearly implies the given condition.  Conversely,
  suppose that $F\maps \bbD\to\bbE$ is a lax framed functor.  We have
  a natural bijection
  \[_g\bbD_f(M,N) \iso \calD(M,f^*Ng^*)
  \]
  which is preserved by $F$, since it preserves restriction.  In other
  words, the diagram
  \[\xymatrix{{}_g\bbD_f(M,N) \ar[r]^\iso \ar[d]_F &
    \calD(M,f^*Ng^*) \ar[d]^F\\
    {}_{Fg}\bbE_{Ff}(FM,FN) \ar[dr]^\iso &
    \calE(FM,F(f^*Ng^*)) \ar[d]^\iso\\
    & \calE(FM,(Ff)^*(FN)(Fg)^*)}\]
  commutes.  Thus, if the right-hand map is surjective (resp.\
  injective), so is the left-hand map.  An analogous argument works
  for an oplax framed functor, using extension instead of restriction.
\end{proof}

This is yet another expression of the fact that in a framed
bicategory, the globular 2-cells carry the information about all the
2-cells.  A similar thing happens for essential surjectivity.

\begin{defn}\xlabel{def:eso}
  A lax or oplax double functor $F:\bbD\to\bbE$ is \textbf{essentially
    surjective} if we can simultaneously make the following choices:
  \begin{itemize}
  \item For each object $C$ of \bbE, an object $A_C$ of \bbD\ and a
    vertical isomorphism $\alpha_C:F(A_C)\iso C$, and
  \item For each horizontal 1-cell $N:C\hto D$ in \bbE, a horizontal
    arrow $M_N\maps A_C\hto A_D$ in \bbE\ and a 2-cell isomorphism
    \[\xymatrix{F(A_C) \ar[r]^{F(M_N)}\ar[d]_{\alpha_C}|{|} \ar@{}[dr]|{\alpha_M\,\iso}&
      F(A_D) \ar[d]^{\alpha_D}\\
      C \ar[r]_N|{|} & D}.\]
  \end{itemize}
\end{defn}

\begin{prop}\xlabel{vert-invar-eso}
  A lax or oplax \emph{framed} functor is essentially surjective, in
  the sense of \xref{def:eso}, if and only if it is essentially
  surjective on vertical categories and each functor
  $\calD(A,B)\to\calE(FA,FB)$ is essentially surjective.
\end{prop}
\begin{proof}
  Clearly \xref{def:eso} implies the given condition.  Conversely,
  suppose that $F$ satisfies the given condition.  Choose isomorphisms
  $\alpha_C:F(A_C)\iso C$ for each object $C$ of \bbE, which exist because
  $F$ is essentially surjective on vertical categories.  Then given
  $N:C\hto D$, we have $\alpha_C^*N\alpha_D^*\maps F(A_C)\hto
  F(A_D)$, so since $F\maps \calD(A_C,A_D)\to \calE(F(A_C),F(A_D))$ is
  essentially surjective, we have an $M_N\maps A_C\to A_D$ and a
  globular isomorphism $F(M_N)\iso \alpha_C^*N\alpha_D^*$.  Composing
  this with the cartesian 2-cell defining $\alpha_C^*N\alpha_D^*$, we
  obtain the desired $\alpha_M$.
\end{proof}

The following theorem and its corollary are the main points of this
section.  Of course, we define a \textbf{double equivalence} to be an
internal equivalence in \Dbll.

\begin{thm}\xlabel{dbl-eqv}
  A strong double functor $F:\bbD\to\bbE$ is part of a double
  equivalence if and only if it is full, faithful, and essentially
  surjective.
\end{thm}

\begin{proof}
  We sketch a construction of an inverse equivalence $G:\bbE\to \bbD$
  for $F$.  Make choices as in \xref{def:eso}, and define $GC=A_C$ and
  $GN=M_N$.  Define $G$ on vertical arrows and 2-cells by composing
  with the chosen isomorphisms; vertical functoriality follows from
  $F$ being full and faithful.  We produce the constraint cells for
  $G$ by composing these isomorphisms with the \emph{inverses} of the
  constraint cells for $F$ and using that $F$ is full and faithful;
  this is why we need $F$ to be strong.
  
  The choices from the definition of essentially surjective then give
  directly a double natural isomorphism $FG\iso\Id_\bbE$, and we can
  produce a double natural isomorphism $GF\iso\Id_\bbD$ by reflecting
  identity maps in \bbE.  Thus $G$ and $F$ form a double equivalence.
\end{proof}

\begin{cor}
  A strong framed functor $F\maps \bbD\to\bbE$ is part of a framed
  equivalence precisely when
  \begin{itemize}
  \item It induces an equivalence $F_0\maps \bbDz\to\bbEz$ on vertical
    categories, and
  \item Each functor $F\maps \calD(A,B)\to \calE(FA,FB)$ is an
    equivalence of categories.
  \end{itemize}
\end{cor}
\begin{proof}
  Combine \xref{framed-ff} and \xref{vert-invar-eso} with
  \xref{dbl-eqv} to see that $F$ has an inverse which is a strong
  double functor, hence also a strong framed functor by
  \xref{dbl-functor->framed}.
\end{proof}

A framed equivalence $F\maps \bbD\toot \bbE\spam G$ clearly includes
an equivalence $F_0\maps \bbDz\toot \bbEz \spam G_0$ of vertical
categories.  It is less clear that it induces a biequivalence
$\calD\eqv \calE$ of horizontal bicategories.  We will see in
appendix~\ref{sec:frbi-vs-bicat}, however, that this is true, though
not trivial.  This lack of triviality, in the following example, was
one of the original motivations for this work.

\begin{eg}\xlabel{eg:ex-not-equiv}
  There are a number of framed bicategories related to \Ex, such as a
  fiberwise version $\Ex_B$ where the objects are already parametrized
  over some space $B$, and an equivariant version $G\Ex$ in which
  everything carries an action by some fixed group $G$.
  In~\cite[19.3.5]{pht} it was observed (essentially) that
  $G\Ex_{G/H}$, the framed bicategory of $G$-equivariant parametrized
  spectra all over the coset space $G/H$, and $H\Ex$, the framed
  bicategory of $H$-equivariant parametrized spectra, are equivalent.

  However, as observed in~\cite{pht}, the language of bicategories
  does not really suffice to describe this fact.  On objects, the
  equivalence goes as follows: if $X$ is a $G$-space over $G/H$, the
  fiber $X_e$ is an $H$-space; while if $Y$ is an $H$-space,
  $G\times_H Y$ is a $G$-space over $G/H$.  But the composites in
  either direction are only homeomorphic, not equal, whereas the
  bicategory \calEx\ described in~\cite{pht} does not include any
  information about homeomorphisms of base spaces.
\end{eg}

\section{Framed adjunctions}
\label{sec:framed-adjunctions}

Adjunctions are one of the most important tools of category theory.
Thus, from a categorical point of view, one of the most serious
problems with \calMod-like bicategories is the lack of a good notion
of adjunction between them.  For example, Ross Street wrote the
following in a review of~\cite{basechange-i}:

\begin{quote}
  Nearly two decades after J. W. Gray's work~\cite{gray:formal-ct},
  the most useful general notion of adjointness for morphisms between
  2-categories has still not emerged.  Perhaps the good notion should
  depend on the kind of 2-categories in mind; 2-categories whose
  arrows are functions or functors are of a different nature from
  those whose arrows are relations or profunctors.
\end{quote}

In fact, motivated by the desire for a good notion of
adjunction,~\cite{basechange-i} and related papers such
as~\cite{verity:base-change,basechange-ii} come very close to our
definition of framed bicategory.  In appendix~\ref{sec:equipments} we
will make a formal comparison; for now we simply develop the theory of
framed adjunctions.

\begin{defn}
  A \textbf{framed adjunction} $F\adj G$ is an internal adjunction in
  the 2-category \FrBil.  Explicitly, it consists of lax framed
  functors $F\maps \bbD\to\bbE$ and $G\maps \bbE\to\bbD$, together
  with framed transformations $\eta\maps \Id_\bbD\to GF$ and $\ep\maps
  FG\to\Id_\bbE$ satisfying the usual triangle identities.  Similarly,
  an \textbf{op-framed adjunction} is an internal adjunction in \FrBic.
\end{defn}

Experience shows that adjunctions in \FrBil\ arise more frequently
than the other two types, hence deserve the unadorned name.  However,
we have the following fundamental result.

\begin{prop}\xlabel{ladj-strong}
  In any framed adjunction $F\adj G$, the left adjoint $F$ is always a
  \emph{strong} framed functor.
\end{prop}
\begin{proofsk}
  This actually follows formally from a general 2-categorical result
  known as `doctrinal adjunction'; see~\cite{kelly:doc-adjn}.  For the
  non-2-categorically inclined reader we sketch a more concrete
  version of the proof.  We first show that the following composite is
  an inverse to $F_\odot\maps FM\odot FN \to F(M\odot N)$:
  \begin{equation}
    F(M\odot N) \xrightarrow{F(\eta\odot \eta)}
    F(GFM \odot GFN) \too[G_\odot]
    FG(FM\odot FN) \too[\ep]
    FM\odot FN
    \label{eq:doc-adjn-odot}
  \end{equation}
  For example, the following diagram shows that the composite in one
  direction is the identity.
  \[\xymatrix{F(GFM\odot GFN) \ar[r]^<>(.5){G_\odot} &
    FG(FM\odot FN) \ar[r]^<>(.5){\ep} \ar[d]^{FG(F_\odot)} &
    FM\odot FN \ar[d]^{F_\odot}\\
    F(M\odot N) \ar[r]^<>(.5){F(\eta)} \ar@(dr,dl)[rr]_{\mathrm{id}}
    \ar[u]^{F(\eta\odot\eta)} &
    FGF(M\odot N) \ar[r]^<>(.5){\ep} &
    F(M\odot N).}\]
  The right-hand square commutes by naturality of $\ep$, the left-hand
  square commutes because $\eta$ is a framed transformation, and the
  lower triangle is one of the triangle identities.  The other
  direction is analogous.

  Similarly, we show that the following composite is an inverse to
  $F_U\maps U_{FA}\to F(U_A)$:
  \begin{equation}
    \xymatrix{F(U_A) \ar[r]^{F(U_\eta)} &
      FU_{GFA} \ar[r]^{G_U} &
      FGU_{FA} \ar[r]^<>(.5){\ep} &
      U_{FA},}\label{eq:doc-adjn-unit}
  \end{equation}
  so that $F$ is strong.
\end{proofsk}

The similarity between~(\ref{eq:doc-adjn-odot})
and~(\ref{eq:doc-adjn-unit}) is obvious.  In fact, these composites
are the mates of the constraint cells for $G$ under an adjunction in a
suitable 2-category; the reader may consult~\cite{kelly:doc-adjn} for
details.  Of course, in an op-framed adjunction, the right adjoint is
strong.

We can now prove \xref{eqv-strong}.

\begin{cor}\xlabel{eqv-strong-2}
  Both functors in a framed equivalence are strong framed functors.
\end{cor}
\begin{proof}
  It is well-known that any classical equivalence of categories can be
  improved to an `adjoint equivalence', meaning an equivalence in
  which the isomorphisms $FG\iso \Id$ and $\Id\iso GF$ are also the
  unit and counit of an adjunction $F\adj G$, and hence their inverses
  are the unit and counit of an adjunction $G\adj F$.  This fact can
  easily be `internalized' to any 2-category, such as \FrBil.  Thus,
  any lax framed functor which is part of a framed equivalence is a
  framed left adjoint, and hence by \xref{ladj-strong} is strong.
\end{proof}

As is the case for categories, we can also characterize framed
adjunctions using universal arrows.  A similar result for double
categories was given in~\cite{garner:double-clubs}.

Recall that given a functor $G\maps \sE\to\sD$, a \emph{universal
  arrow} to $G$ is an arrow $\eta\maps A\to GFA$ in \sD, for some
object $FA\in\sE$, such that any other arrow $A\to GY$ factors through
$\eta$ via a unique map $FA\to Y$ in \sE.  Similarly, if $G\maps
\bbE\to\bbD$ is a framed functor, we define a \emph{universal 2-cell}
to be a 2-cell $\eta\maps M\to GFM$ in \bbD, not in general globular,
whose left and right frames are universal arrows in \bbDz, and such
that any 2-cell $M\to GN$ factors through $\eta$ via a unique 2-cell
$FM\to N$ in \bbE.

\begin{prop}\xlabel{dbladj-univ}
  Let $G\maps \bbE\to\bbD$ be a lax framed functor.  Then $G$ has a
  framed left adjoint if and only if the following are true.
  \begin{enumerate}
  \item For every object $A$ in \bbD, there is a universal arrow $A\to
    GFA$.\label{item:dbladj-univ-arrow}
  \item For every horizontal 1-cell $M\maps A\hto B$ in \bbD, there is
    a universal 2-cell $M\to GFM$, as described
    above.\label{item:dbladj-univ-twocell}
  \item If $M\to GFM$ and $N\to GFN$ are universal 2-cells, then so is
    the composite\label{item:dbladj-comp}
    \[\xymatrix{
      \ar[d]\ar[rr]|{|}^M \ar@{}[drr]|{\mathrm{univ}} &&
      \ar[d]\ar[rr]|{|}^N \ar@{}[drr]|{\mathrm{univ}} &&
      \ar[d]\\
      \ar[rr]|{GFM} \ar@{=}[d] \ar@{}[drrrr]|{\Downarrow G_\odot} &&
      \ar[rr]|{GFN} && \ar@{=}[d]\\
      \ar[rrrr]|{|}_{G(FM\odot FN)} &&&&}
    \]
  \item If $A\to GFA$ is universal, then so is the
    composite\label{item:dbladj-unit}
    \[\xymatrix{
      \ar[rr]|{|}^{U_A} \ar[d]_{\mathrm{univ}}
      \ar@{}[drr]|{\Downarrow U_\mathrm{univ}} &&
      \ar[d]^{\mathrm{univ}}\\
      \ar[rr]|{U_{GFA}} \ar@{=}[d]
      \ar@{}[drr]|{\Downarrow G_0} &&
      \ar@{=}[d]\\
      \ar[rr]|{|}_{G(U_{FA})} &&
    }\]
  \end{enumerate}
  If $G$ is strong, then~\ref{item:dbladj-comp} simplifies to `the
  horizontal composite of universal 2-cells is universal'
  and~\ref{item:dbladj-unit} simplifies to `the horizontal unit of a
  universal arrow is a universal 2-cell'.
\end{prop}
\begin{proofsk}
  It is straightforward to show that if $G$ has a left adjoint, then
  the conditions are satisfied.  Conversely,
  conditions~\ref{item:dbladj-univ-arrow}
  and~\ref{item:dbladj-univ-twocell} clearly guarantee that $G_0$ and
  $G_1$ both have left adjoints $F_0$ and $F_1$, and that $LF_1 \iso
  F_0L$ and $RF_1 \iso F_0R$.  Since \bbE\ is a framed bicategory, we
  can redefine $F_1$ by restricting along these isomorphisms to ensure
  that $LF_1=F_0L$ and $RF_1 = F_0R$.

  Conditions~\ref{item:dbladj-comp} and~\ref{item:dbladj-unit} then
  supply the constraints to make $F$ into a strong framed functor.
  The universal cells give a double transformation $\eta\maps \Id\to
  GF$ and the counit $\ep\maps FG\to \Id$ is constructed as usual.
  The last statement follows because anything isomorphic to a
  universal arrow is universal.
\end{proofsk}

Since strict 2-functors preserve internal adjunctions, our
2-functorial ways of constructing framed bicategories also give us a
ready supply of framed adjunctions.

\begin{eg}\xlabel{eg:monadj->dbladj}
  Since \Mod\ is a 2-functor, any monoidal adjunction between monoidal
  categories with coequalizers preserved by $\ten$ gives rise to a
  framed adjunction.  Here by a \emph{monoidal adjunction} we mean an
  adjunction in the 2-category \MonCatl\ of monoidal categories and
  lax monoidal functors.

  For example, if $f\maps R\to S$ is a homomorphism of
  commutative rings, we have an induced monoidal adjunction
  \[f_!\maps \mathbf{Mod}_R\rightleftarrows \mathbf{Mod}_S\spam f^*\]
  and therefore a framed adjunction
  \[\Mod(f_!)\maps \Mod(R)\rightleftarrows \Mod(S)\spam \Mod(f^*).\]
\end{eg}

\begin{eg}\xlabel{eg:lexadj->dbladj}
  Since \Span\ is a strict 2-functor, any adjunction $f^*\maps
  \sE\toot\sF\spam f_*$ between categories with pullbacks gives rise
  to an op-framed adjunction $\Span(\sE)\toot\Span(\sF)$.  If $f^*$
  also preserves pullbacks, then this adjunction lies in \FrBi, hence
  is also a framed adjunction.
\end{eg}

\section{Monoidal framed bicategories}
\label{sec:monoidal-structures}

Most of our examples also have an `external' monoidal structure.  For
example, if $M$ is an $(A,B)$-bimodule and $N$ is a $(C,D)$-bimodule,
we can form the $(A\ten C, B\ten D)$-bimodule $M\ten N$.  The
definition of a `monoidal bicategory' involves many coherence axioms
(see~\cite{tricats,nick:tricats}), but for \emph{framed} bicategories
we can simply invoke general 2-category theory once again.
\nocite{garner:double-clubs}

In any 2-category with finite products, we have the notion of a
\emph{pseudo-monoid}: this is an object $A$ equipped with
multiplication $A\times A\to A$ and unit $1\to A$ satisfying the usual
monoid axioms up to coherent isomorphism.  A pseudo-monoid in \Cat\ is
precisely an ordinary monoidal category.  Thus, it makes sense to
define a monoidal framed bicategory to be a pseudo-monoid in \FrBi.
What this means is essentially the following.

\begin{defn}
  A \textbf{monoidal framed bicategory} is a framed bicategory equipped
  with a strong framed functor $\ten\maps \bbD\times\bbD\to\bbD$, a
  unit $I\in\bbDz$, and framed natural constraint isomorphisms
  satisfying the usual axioms.
\end{defn}

If we unravel this definition more explicitly, it says the following.
\begin{enumerate}
\item \bbDz\ and \bbDo\ are both monoidal categories.
\item $I$ is the monoidal unit of \bbDz\ and $U_I$ is the monoidal
  unit of \bbDo.
\item The functors $L$ and $R$ are strict monoidal.
\item We have an `interchange' isomorphism
  \[\fx\maps (M\ten P)\odot (N\ten Q)\iso (M\odot N)\ten (P\odot Q)\]
  and a unit isomorphism
  \[\fu\maps U_{A\ten B} \iso (U_A \ten U_B)\]
  satisfying appropriate axioms (these arise from the constraint data
  for the strong framed functor $\ten$).
\item The associativity and unit isomorphisms for $\ten$ are framed
  transformations.
\end{enumerate}

As we saw in \S\ref{sec:2cat-frbi}, a strong framed functor such as
$\ten$ preserves cartesian and opcartesian arrows.  Thus we
automatically have isomorphisms such as $f^*M\ten g^*N \iso (f\ten
g)^*(M\ten N)$.

\begin{egs}\xlabel{egs:monoidal}
  Many of our examples of framed bicategories are in fact monoidal.
  \begin{itemize}
  \item The framed bicategory \Mod, and more generally $\Mod(\sC)$ for
    a symmetric monoidal \sC, is monoidal under the tensor product of
    rings and bimodules.  Note that the tensor product of bimodules
    referred to here is `external': if $M$ is an $(R,S)$-bimodule and
    $N$ is a $(T,V)$-bimodule, then $M\ten N$ is an $(R\ten T, S\ten
    V)$-bimodule.
  \item If \sC\ has finite limits, then $\Span(\sC)$ is a monoidal
    framed bicategory under the cartesian product of objects and
    spans.
  \item \Ex\ is monoidal under the cartesian product of spaces and the
    `external smash product' $\overline{\sm}$ of parametrized spectra.
  \item \nCob\ is monoidal under disjoint union of manifolds and
    cobordisms.
  \item $\Dist(\V)$ is monoidal under the tensor product of
    \V-categories (see~\cite[\S1.4]{kelly}).
  \end{itemize}
\end{egs}

\begin{eg}
  Recalling that monoidal categories can be identified with vertically
  trivial framed bicategories, it is easy to check that a vertically
  trivial monoidal framed bicategory is the same as a category with
  two interchanging monoidal structures.  More generally, if \bbD\ is
  any monoidal framed bicategory, then the category $\calD(I,I)$
  inherits two interchanging monoidal structures $\odot$ and $\ten$.
  By the Eckmann-Hilton argument, any two such interchanging monoidal
  structures agree up to isomorphism and are braided.
\end{eg}

We emphasize that the associativity and unit constraints are
\emph{vertical} isomorphisms.  For example, in the monoidal framed
bicategory \Mod, the associativity constraint on objects is the ring
isomorphism $(A\ten B)\ten C \iso A\ten (B\ten C)$.  This is to be
contrasted with the classical notion of `monoidal bicategory' in which
the constraints are 1-cells, which would correspond to bimodules in
this case.  So while a framed bicategory obviously has an underlying
horizontal bicategory, it requires proof that a monoidal framed
bicategory has an underlying monoidal bicategory; see
appendix~\ref{sec:frbi-vs-bicat}.

We observe, in passing, that an external monoidal structure
automatically preserves dual pairs.

\begin{prop}
  If $(M,N)$ and $(P,Q)$ are dual pairs in a monoidal framed
  bicategory, then so is $(M\ten P,N\ten Q)$.
\end{prop}
\begin{proof}
  It is easy to see that any strong framed functor preserves dual
  pairs, and $\ten$ is a strong framed functor.
\end{proof}

Now, just as an ordinary monoidal category can be braided or
symmetric, so can a pseudo-monoid in an arbitrary 2-category with
products.  We define a \textbf{braided} or \textbf{symmetric} monoidal
framed bicategory to be essentially a braided or symmetric
pseudo-monoid in \FrBi.

More explicitly, a braided monoidal framed bicategory is a monoidal
framed bicategory such that \bbDz\ and \bbDo\ are braided monoidal
with braidings \fs, the functors $L$ and $R$ are braided monoidal, and
the following diagrams commute:
\[\xymatrix{(M\odot N)\ten (P\odot Q) \ar[r]^\fs\ar[d]_\fx &
  (P\odot Q)\ten (M \odot N)\ar[d]^\fx\\
  (M\ten P)\odot (N\ten Q) \ar[r]_\fs &
  (P\ten M) \odot (Q \ten N)}
\]
\[\xymatrix{(U_A \ten U_B) \ar[r]^(0.55)\fu \ar[d]_\fs &
  U_{A\ten B} \ar[d]^{U_\fs}\\
  U_B\ten U_A \ar[r]_(0.55)\fu &
  U_{B\ten A}}.
\]
A symmetric monoidal framed bicategory is a braided monoidal framed
bicategory such that \bbDz\ and \bbDo\ are symmetric.

\begin{egs}
  All the examples of monoidal framed bicategories given in
  \xref{egs:monoidal} are in fact symmetric monoidal.
\end{egs}

\begin{eg}
  If \bbD\ is a braided or symmetric monoidal framed bicategory, then
  $\calD(I,I)$ inherits two interchanging monoidal structures, one of
  which is braided, and therefore it is essentially a symmetric
  monoidal category.  Conversely, the vertically trivial monoidal
  framed bicategory corresponding to any symmetric monoidal category
  is a symmetric monoidal framed bicategory.
\end{eg}

We now define the morphisms between monoidal framed bicategories.  As
usual, these come in three flavors.

\begin{defn}
  A \textbf{lax monoidal framed functor} between monoidal framed
  bicategories $\bbD,\bbE$ consists of the following structure and
  properties.
  \begin{itemize}
  \item A lax framed functor $F\maps \bbD\to\bbE$.
  \item The structure of a lax monoidal functor on $F_0$ and $F_1$.
  \item Equalities $LF_1=F_0L$ and $RF_1=F_0R$ of lax monoidal functors.
  \item The composition constraints for the lax framed functor $F$ are
    monoidal natural transformations.
  \end{itemize}
  It is \textbf{strong} if $F$ is a strong framed functor and $F_0$ and
  $F_1$ are strong monoidal functors.  If \bbD\ and \bbE\ are braided
  (resp.\ symmetric), then $F$ is \textbf{braided} (resp.\
  \textbf{symmetric}) if $F_0$ and $F_1$ are.  We have a dual definition
  of \textbf{oplax monoidal framed functor}.  A \textbf{monoidal framed
    transformation} is a framed transformation such that $\alpha_0$
  and $\alpha_1$ are monoidal transformations.
\end{defn}

These definitions give various 2-categories, each of which has its own
attendant notion of equivalence and adjunction.  We will not spell
these out explicitly.

\begin{egs}
  The 2-functor \Mod\ lifts to a 2-functor from symmetric monoidal
  categories with coequalizers preserved by $\ten$ to symmetric
  monoidal framed bicategories.  Similarly, \Span\ lifts to a
  2-functor landing in symmetric monoidal framed bicategories.
\end{egs}

Finally, we consider what it means for a framed bicategory to be
`closed monoidal'.

\begin{defn}
  A monoidal framed bicategory \bbD\ is \textbf{externally closed} if
  for any objects $A,B,C,D$, the functor
  \[\ten\maps \calD(A,C)\times \calD(B,D) \too \calD(A\ten B, C\ten D)\]
  has right adjoints in each variable, which we write $\xlhd$ and $\xrhd$.
\end{defn}

Explicitly, this means that for horizontal 1-cells $M\maps A\hto C$,
$N\maps B\hto D$, and $P\maps A\ten B\hto C\ten D$, there are 1-cells
$N\xrhd P$ and $P\xlhd M$ and bijections
\[\calD(M\ten N, P) \iso \calD(M, N\xrhd P) \iso \calD(N, P\xlhd M).\]
Of course, if \bbD\ is symmetric, then $\xlhd$ and $\xrhd$ agree,
modulo suitable isomorphisms.

\begin{egs}
  The monoidal framed bicategory \Mod\ is externally closed, as are
  \Ex\ and $\Dist(\sV)$.  If \sC\ is locally cartesian closed, then
  $\Span(\sC)$ is also externally closed.
\end{egs}

\section{Involutions}
\label{sec:involutions}

In most of our examples, the `directionality' of the horizontal 1-cells
is to some extent arbitrary.  For example, an $(A,B)$-bimodule could
just as well be regarded as a $(B\op,A\op)$-bimodule.  We now define a
structure which encodes this fact formally.

If \bbD\ is a framed bicategory, we write $\bbD\hop$ for its
`horizontal dual': $\bbD\hop$ has the same vertical category as \bbD,
but a horizontal 1-cell from $A$ to $B$ in $\bbD\hop$ is a horizontal
1-cell from $B$ to $A$ in \bbD, and the 2-cells are similarly  flipped
horizontally.

\begin{defn}
  An \textbf{involution} on a framed bicategory \bbD\ consists of the
  following.
  \begin{enumerate}
  \item A strong framed functor $(-)\op\maps \bbD\hop \to \bbD$.
  \item A framed natural isomorphism $\xi\maps ((-)\op)\op \iso
    \Id_\bbD$ such that $(\xi_A)\op = \xi_{A\op}$; thus $\xi$ and
    $\xi^{-1}$ make $(-)\op$ into an adjoint equivalence.
  \end{enumerate}
  We say an involution is \textbf{vertically strict} if the vertical
  arrow components of $\xi$ are identities.  If \bbD, $(-)\op$, and
  $\xi$ are all monoidal (resp.\ symmetric monoidal), we say that the
  involution is \textbf{monoidal} (resp.\ \textbf{symmetric monoidal}).
\end{defn}

The strong functoriality of $(-)\op$ implies that we have
\begin{align*}
  (U_A)\op &\iso U_{A\op}\\
  (M\odot N)\op &\iso N\op \odot M\op.\\
  (f^*Mg^*)\op &\iso (g\op)^*(M\op)(f\op)^*\\
  (f_!Mg_!)\op &\iso (g\op)_!(M\op)(f\op)_!.
\end{align*}
In particular, we have $(A_f)\op \iso {}_{f\op}(A\op)$ and dually.
If the involution is monoidal, we also have
\begin{align*}
  (A\ten B)\op &\iso A\op \ten B\op\\
  I\op &\iso I.
\end{align*}

\begin{egs}\xlabel{egs:involution}
  Most of our examples are equipped with vertically strict symmetric
  monoidal involutions.
  \begin{itemize}
  \item The involution on \Mod\ takes a ring $A$ to the opposite ring
    $A\op$, and an $(A,B)$-bimodule to the same abelian group regarded
    as a $(B\op,A\op)$-bimodule.
  \item The involution on $\Dist(\sV)$ takes a \V-category to its
    opposite and reverses distributors in an obvious way.
  \item The involution on $\Span(\sC)$ takes each object to itself,
    and a span $A\oot[f] X \too[g] B$ to the span $B\oot[g] X \too[f]
    A$.
  \item The involution on \nCob\ takes a manifold $M$ to the manifold
    $M\op$ with the opposite orientation, and reverses the direction
    of cobordisms in an obvious way.
  \end{itemize}
  In all these cases, the 2-cell components of $\xi$ can also be
  chosen to be identities, but this is not true for all involutions,
  even vertically strict ones.
  \begin{itemize}
  \item The involution on \Ex\ takes each space to itself, but takes a
    spectrum $E$ parametrized over $B\times A$ to the pullback
    $\fs^*E$ over $A\times B$, where $\fs$ is the symmetry isomorphism
    $A\times B\iso B\times A$.  Here $\fs^*\fs^*E$ is only canonically
    isomorphic to $E$, by pseudofunctoriality.
  \end{itemize}
\end{egs}

In~\cite[16.2.1]{pht} an involution on a \emph{bicategory} was defined
to be essentially a pseudofunctor $(-)\op\maps \calB\op\to\calB$
equipped with a pseudonatural transformation $\xi\maps ((-)\op)\op
\iso \Id_\calB$ whose 1-cell components are identities (although the
unit axiom for $\xi$ was omitted).  It is easy to see that any
vertically strict involution on \bbD\ gives rise to an involution on
\calD.  All the above examples are vertically strict, but in
\S\ref{sec:modules} and \S\ref{sec:monoids-monfib} we will see
examples which are not.

Any symmetric monoidal category, considered as a vertically trivial
framed bicategory, has a canonical involution.  The functor $(-)\op$
is the identity on 1-cells (the objects of the monoidal category), and
its composition constraint is the symmetry isomorphism:
\[(A\odot B)\op = A\odot B \too[\iso] B\odot A = B\op \odot A\op.
\]
All the components of $\xi$ are identities.  In fact, to give an
involution on a vertically trivial framed bicategory which is the
identity on 1-cells and for which $\xi$ is an identity is essentially
to give a symmetry for the corresponding monoidal category.  Thus, we
may view an involution on a framed bicategory as a generalization of a
symmetry on a monoidal category.

One consequence of a monoidal category's being symmetric is that if it
is closed, then the left and right internal-homs are isomorphic.  The
original motivation in~\cite{pht} for introducing involutions was to
obtain a similar result for closed bicategories;
see~\cite[16.3.5]{pht}.  Of course, this is also true for framed
bicategories.

\begin{prop}
  If \bbD\ is a closed framed bicategory equipped with an involution,
  then we have
  \[M\rhd N \iso (N\op \lhd M\op)\op.\]
\end{prop}
\begin{proof}
  Since $(-)\op$ is a framed equivalence, it is locally full and
  faithful.  Thus, if $M\maps A\hto B$, $N\maps C\hto B$, and $P\maps
  C\hto A$, we have
  \begin{align*}
    \calD(C,A)(P,M\rhd N)
    &\iso \calD(C,B)(P\odot M, N)\\
    &\iso \calD(B\op,C\op)((P\odot M)\op,N\op)\\
    &\iso \calD(B\op,C\op)(M\op\odot P\op,N\op)\\
    &\iso \calD(A\op,C\op)(P\op, N\op \lhd M\op)\\
    &\iso \calD\Big((C\op)\op,(A\op)\op\Big)\Big((P\op)\op,(N\op \lhd M\op)\op\Big)\\
    &\iso \calD(C,A)(P,(N\op \lhd M\op)\op)
  \end{align*}
  so the result follows by the Yoneda lemma.
\end{proof}

\section{Monoids and modules}
\label{sec:modules}

In most of our examples of monoidal framed bicategories, the external
monoidal structure and the horizontal composition are more closely
related than is captured by the interchange isomorphism: namely, the
horizontal composition $M\odot N$ is a subobject or quotient of the
external product $M\ten N$.  For example, in \Mod\ the tensor product
$M\ten_R N$ is a quotient of the external product $M\ten N$, while in
\Span\ the pullback $M\times_B N$ is a subobject of the external
product $M\times N$.  An analogous relationship holds between the
bicategorical homs $\lhd,\rhd$ and the external homs $\xlhd,\xrhd$.

In this section we will generalize the construction of the framed
bicategory $\Mod(\sC)$ of monoids and modules from
\xref{eg:bimodules}, replacing the monoidal category \sC\ with a
framed bicategory \bbD.  This describes one general class of examples
in which the horizontal composition of `bimodules' is defined as a
coequalizer.  In
\S\S\ref{sec:monoidal-fibrations}--\ref{sec:monfib-to-framed}, we
will investigate framed bicategories constructed in a way analogous to
\Span.  We will then combine these two constructions in
\S\ref{sec:monoids-monfib} to define framed bicategories of internal and
enriched categories.

\begin{defn}\xlabel{def:mod}
  Let \bbD\ be a framed bicategory.
  \begin{itemize}
  \item A \textbf{monoid} in \bbD\ consists of an object $R$, a
    horizontal 1-cell $A\maps R\hto R$, and globular 2-cells $e\maps
    R\to A$ and $m\maps A\odot A\to A$ called `unit' and
    `multiplication' such that the standard diagrams commute.  Thus it
    is just a monoid in the ordinary monoidal category $\calD(R,R)$.
  \item A \textbf{monoid homomorphism} $(R,A)\to (S,B)$ consists of a
    vertical arrow $f\maps R\to S$ and a 2-cell $\phi\maps A
    \sto{f}{f} B$ such that $\phi\circ e = e$ and $\phi\circ m = m
    \circ (\phi\odot\phi)$.
  \item A \textbf{bimodule} from a monoid $(R,A)$ to a monoid $(S,B)$ is
    a horizontal 1-cell $M\maps R\hto S$ together with action maps
    $a_\ell\maps A\odot M\to M$ and $a_r\maps M\odot B\to M$ obeying
    the obvious compatibility axioms.
  \item Let $(f,\phi)\maps (R,A)\to (S,B)$ and $(g,\psi)\maps (T,C)\to
    (U,D)$ be monoid homomorphisms and $M\maps (R,A)\hto (T,C)$,
    $N\maps (S,B)\hto (U,D)$ be bimodules.  A
    \textbf{$(\phi,\psi)$-equivariant map} is a 2-cell $\alpha\maps
    M\sto{f}{g} N$ such that $a_\ell(\phi\odot\alpha) = \alpha a_\ell$
    and $a_r(\alpha\odot\psi) = \alpha a_r$.
  \item Let $M\maps R\hto S$ be an $(A,B)$-bimodule and $N\maps S\hto
    T$ be a $(B,C)$-bimodule.  Their \textbf{tensor product} is the
    following coequalizer in $\calD(A,C)$, if it exists:
    \[M \odot B \odot N \toto M \odot N \to M\odot_B N.\]
  \end{itemize}
\end{defn}

Of course, if \bbD\ is a monoidal category, these notions reduce to
the usual ones.

\begin{eg}
  If \sC\ has pullbacks, then a monoid in $\Span(\sC)$ is an internal
  category in \sC, and a monoid homomorphism is an internal functor.
  A bimodule in $\Span(\sC)$ is an `internal distributor'.
\end{eg}

\begin{eg}
  A monoid in $\Mod$ consists of a ring $R$ together with an
  $R$-algebra $A$, and a monoid homomorphism $(R,A)\to (S,B)$ consists
  of a ring homomorphism $f\maps R\to S$ and an $f$-equivariant
  algebra map $A\to B$.  A bimodule in \Mod\ is just a bimodule for
  the algebras.
\end{eg}

In order to define a framed bicategory of monoids and bimodules in
\bbD, we need to know that coequalizers exist and are well-behaved.

\begin{defn}
  A framed bicategory \bbD\ has \textbf{local coequalizers} if each
  category $\calD(A,B)$ has coequalizers and $\odot$ preserves
  coequalizers in each variable.  We introduce the following notations.
  \begin{itemize}
  \item $\FrBilq$ denotes the \emph{full} sub-2-category of \FrBil\
    determined by the framed bicategories with local coequalizers.
  \item $\FrBilnq$ denotes the locally full sub-2-category of \FrBil\
    determined by the framed bicategories with local coequalizers and
    the \emph{normal} lax framed functors.
  \item $\FrBiq$ denotes the locally full sub-2-category of \FrBi\
    determined by the framed bicategories with local coequalizers and
    the strong framed functors which preserve local coequalizers.
  \end{itemize}
\end{defn}

Note that if \bbD\ is closed, as defined in
\S\ref{sec:duality-theory}, then $\odot$ preserves all colimits since
it is a left adjoint.  The following omnibus theorem combines all our
results about monoids and modules in framed bicategories.

\begin{thm}\xlabel{frbi->mod}
  Let \bbD\ be a framed bicategory with local coequalizers.  Then
  there is a framed bicategory $\Mod(\bbD)$ of monoids, monoid
  homomorphisms, bimodules, and equivariant maps in \bbD.  Moreover:
  \begin{itemize}
  \item $\Mod(\bbD)$ also has local coequalizers.
  \item If \bbD\ is closed and each category $\calD(A,B)$ has
    equalizers, then $\Mod(\bbD)$ is closed.
  \item If \bbD\ is monoidal and its external product $\ten$ preserves
    local coequalizers, then $\Mod(\bbD)$ has both of these
    properties.  If \bbD\ is symmetric, so is $\Mod(\bbD)$.  If \bbD\
    is externally closed and each category $\calD(A,B)$ has
    equalizers, then $\Mod(\bbD)$ is externally closed.
  \item If \bbD\ is equipped with an involution, so is $\Mod(\bbD)$.
    If the involution of \bbD\ is monoidal or symmetric monoidal, so
    is that of $\Mod(\bbD)$.
  \item $\Mod$ defines 2-functors $\FrBilq\to\FrBilnq$ and
    $\FrBiq\to\FrBiq$, and similarly for the monoidal versions.
  \end{itemize}
\end{thm}

Even if $F$ is a strong framed functor, $\Mod(F)$ is only lax unless
$F$ preserves local coequalizers.  If $F$ is oplax, we cannot even
define $\Mod(F)$.  Of course, there is a dual construction \Comod, but
it arises much less frequently in practice.

\begin{eg}
  If \sC\ is a monoidal category with coequalizers preserved by
  $\ten$, then $\Mod(\sC)$ has local coequalizers, so we have a framed
  bicategory $\Mod(\Mod(\sC))$ of algebras and bimodules in \sC.
\end{eg}

\begin{eg}\xlabel{eg:internal-dist}
  If \sC\ is a category with pullbacks and coequalizers preserved by
  pullback, then $\Span(\sC)$ has local coequalizers, so we have a
  framed bicategory $\Mod(\Span(\sC))$ of internal categories and
  distributors in \sC.
\end{eg}

\begin{eg}\xlabel{eg:mat->dist}
  When \sV\ is a cocomplete closed monoidal category, we can also
  construct the framed bicategory $\Dist(\sV)$ of enriched categories
  and distributors in this way.  We first define the framed bicategory
  $\Mat(\sV)$ as follows: its vertical category is \Set, and the
  category $\calMat(\sV)(A,B)$ is the category of $A\times B$ matrices
  $(M_{ab})_{a\in A,b\in B}$ of objects of \sV.  Composition is by
  `matrix multiplication'.  It is then easy to check that $\Mat(\sV)$
  has local coequalizers and that $\Mod(\Mat(\sV))\iso \Dist(\sV)$.
  The monoidal category $\calMat(\sV)(A,A)$ is also called the
  category of \emph{\sV-graphs} with object set $A$.
\end{eg}

\begin{eg}
  Unlike these examples, \Ex\ does not have local coequalizers.  We
  will see a replacement for `$\Mod(\Ex)$' in
  \S\ref{sec:monoids-monfib}.
\end{eg}

The rest of this section is devoted to the proof of \xref{frbi->mod},
breaking it up into a series of propositions for clarity.  Although
long, the proof is routine and follow-your-nose, so it can easily be
skipped.

\begin{prop}
  If \bbD\ is a framed bicategory with local coequalizers, then there
  is a framed bicategory $\Mod(\bbD)$ of monoids, monoid
  homomorphisms, bimodules, and equivariant maps in \bbD, and it also
  has local coequalizers.
\end{prop}
\begin{proof}
  The proof that $\Mod(\bbD)$ is a double category is similar to the
  case of a monoidal category.  For example, we need the fact that
  $\odot$ preserves coequalizers to show that $M\odot N$ is a bimodule
  and that the tensor product is associative.  To define the
  horizontal composite of bimodule maps $\alpha\maps M
  \sto{\phi}{\psi} N$ and $\beta\maps P \sto{\psi}{\chi} Q$ (where
  $\phi\maps A\sto{f}{f} D$, $\psi\maps B\sto{g}{g}E$, and $\chi\maps
  C\sto{h}{h} F$ are monoid homomorphisms), we start with the
  composite
  \begin{equation}
    \xymatrix{R \ar[r]^M\ar[d]_f \ar@{}[dr]|{\alpha} &
    S \ar[d]|g \ar[r]^P \ar@{}[dr]|{\beta} & T\ar[d]^h\\
    U \ar[r]|N\ar@{=}[d] \ar@{}[drr]|{\mathrm{coeq}} &
    V \ar[r]|Q & W \ar@{=}[d]\\
    U \ar[rr]_{N\odot_E Q} && W}\label{eq:bimod-comp}
  \end{equation}
  which we would like to factor through the coequalizer defining
  $M\odot_B P$.  However, that coequalizer lives in $\calD(R,T)$,
  whereas~(\ref{eq:bimod-comp}) is not globular.  But since \bbD\ is a
  framed bicategory, we can factor~(\ref{eq:bimod-comp}) through a
  cartesian arrow to get a map
  \[M\odot P \to f^*(N\odot_E Q)g^*
  \]
  in $\calD(R,T)$, and then apply the universal property of the
  coequalizer to get a map $M\odot_B P\to f^*(N\odot_E Q)g^*$, and
  hence $M\odot_B P\sto{f}{g}(N\odot_E Q)$.  This defines a
  $(\phi,\chi)$-equivariant map which we call the horizontal composite
  $\alpha\odot_\psi \beta$.  The axioms for a double category follow
  directly.

  We now show that $\Mod(\bbD)$ is a framed bicategory.  By
  \xref{thm:framed}, it suffices to show that it has restrictions.
  Thus, suppose that $A\maps R\hto R$, $B\maps S\hto S$, $C\maps T\hto
  T$, and $D\maps U\hto U$ are monoids in \bbD, $M\maps S\hto U$ is a
  $(B,D)$-bimodule, and $\phi\maps A\sto{f}{f} B$ and $\psi\maps
  C\sto{g}{g} D$ are monoid homomorphisms.  We then have the
  restriction $f^*Mg^*\maps R\hto T$ in \bbD.  By composing the
  cartesian arrow in \bbD\ with $\phi$ or $\psi$ and using the actions
  of $B$ and $D$ on $M$, then factoring through the cartesian arrow,
  we obtain actions of $A$ and $C$ on $f^*Mg^*$.  For example, the
  action of $A$ on $f^*Mg^*$ is determined by the equality
  \[\begin{array}{c}
    \xymatrix{
      R \ar[r]^A \ar[d]_f \ar@{}[dr]|\phi &
      R \ar[r]^{f^*Mg^*} \ar[d]|f \ar@{}[dr]|{\mathrm{cart}} &
      T \ar[d]^g\\
      S \ar[r]|B \ar@{=}[d] \ar@{}[drr]|{\mathrm{act}} &
      S \ar[r]_M & U \ar@{=}[d]\\
      S \ar[rr]_M && U
    }\end{array} = \begin{array}{c}\xymatrix{
      R \ar[r]^A \ar@{=}[d] \ar@{}[drr]|{\mathrm{act}} &
      R \ar[r]^{f^*Mg^*} & T \ar@{=}[d]\\
      R \ar[rr]|{f^*Mg^*} \ar[d]_f \ar@{}[drr]|{\mathrm{cart}} && T \ar[d]^g\\
      S \ar[rr]_M && U
    }\end{array}
  \]
  It is straightforward to check that with this structure, the
  cartesian arrow $f^*Mg^*\sto{f}{g} M$ in \bbD\ defines a cartesian
  arrow $\phi^*M\psi^*\sto{\phi}{\psi} M$ in $\Mod(\bbD)$.
\end{proof}

\begin{prop}
  \Mod\ defines a 2-functor $\FrBilq\to\FrBilnq$, which restricts to a
  2-functor $\FrBiq\to\FrBiq$.
\end{prop}
\begin{proof}
  Let $\bbD,\bbE\in\FrBilq$ and let $F\maps \bbD\to\bbE$ be a lax
  framed functor.  Then $F$ preserves monoids, monoid homomorphisms,
  bimodules, and equivariant maps, for the same reasons that lax
  monoidal functors do.  We define the unit constraint for $\Mod(F)$
  to be the identity on $FA$, and the composition constraint to be the
  result of factoring the composite
  \begin{equation}
    FM \odot FN \too[F_\odot] F(M\odot N) \too F(M\odot_B N)\label{eq:mod-fr-1}
  \end{equation}
  through the coequalizer
  \begin{equation}
    FM\odot FN \to FM\odot_{FB} FN\label{eq:mod-fr-2}
  \end{equation}
  It is straightforward to check that this makes $\Mod(F)$ into a
  normal lax double functor.  Similarly, the components of a framed
  transformation $F\to G$ define a framed transformation
  $\Mod(F)\to\Mod(G)$.

  Finally, if $F$ is strong and preserves local coequalizers,
  then~(\ref{eq:mod-fr-1}) is a coequalizer of the same maps
  that~(\ref{eq:mod-fr-2}) is.  Hence the induced composition
  constraint is an isomorphism, so $\Mod(F)$ is strong.  It is easy to
  see that $\Mod(F)$ also preserves local coequalizers, so that it
  lies in \FrBiq.
\end{proof}

\begin{prop}\xlabel{mod-pres-monoidal}
  If \bbD\ is a monoidal framed bicategory with local coequalizers
  preserved by $\ten$, then so is $\Mod(\bbD)$.  If \bbD\ is
  symmetric, so is $\Mod(\bbD)$.
\end{prop}
\begin{proof}
  It is easy to check that the 2-functor $\Mod\maps \FrBiq\to\FrBiq$
  preserves products, so it must preserve pseudo-monoids and symmetric
  pseudo-monoids.
\end{proof}

\begin{prop}
  Suppose that \bbD\ has local coequalizers and each category
  $\calD(A,B)$ has equalizers.  If \bbD\ is closed, then $\Mod(\bbD)$
  is closed.  If \bbD\ is monoidal and externally closed with local
  coequalizers preserved by $\ten$, then $\Mod(\bbD)$ is externally
  closed.
\end{prop}
\begin{proof}
  Just as for monoidal categories.
\end{proof}

\begin{prop}
  If \bbD\ has local coequalizers and is equipped with an involution, so
  is $\Mod(\bbD)$.  If \bbD, $\Mod(\bbD)$, and the involution on \bbD\
  are monoidal or symmetric monoidal, so is the involution on
  $\Mod(\bbD)$.
\end{prop}
\begin{proof}
  It is easy to see that $\Mod(\bbD\hop)\eqv\Mod(\bbD)\hop$, so we can
  simply apply the 2-functor $\Mod$ to $(-)\op$ and $\xi$.
\end{proof}

Note, however, that since the vertical arrow components of $\xi$ in
$\Mod(\bbD)$ are defined from the \emph{2-cell} components of $\xi$ in
\bbD, the involution of $\Mod(\bbD)$ may not be vertically strict even
if the involution of \bbD\ is so.

\section{Monoidal fibrations}
\label{sec:monoidal-fibrations}

The generalized \Mod\ construction from \S\ref{sec:modules} defines a
horizontal composition from an external product via a coequalizer.  In
\S\ref{sec:monfib-to-framed} we will explain how in a cartesian
situation, horizontal compositions can be constructed using a pullback
or equalizer-type construction instead.  The basic input for this
construction is a structure called a `monoidal fibration', which
includes base change operations and an external product, but \emph{a
  priori} no horizontal composition.

\begin{defn}
  A \textbf{monoidal fibration} is a functor $\Phi\maps \sA\to\sB$ such
  that
  \begin{enumerate}
  \item \sA\ and \sB\ are monoidal categories;
  \item $\Phi$ is a fibration and a strict monoidal functor; and
  \item The tensor product $\ten$ of \sA\ preserves cartesian arrows.
  \end{enumerate}
  If $\Phi$ is also an opfibration and $\ten$ preserves opcartesian
  arrows, we say that $\Phi$ is a \textbf{monoidal bifibration}.  We say
  that $\Phi$ is \textbf{braided} (resp.\ \textbf{symmetric}) if \sA, \sB,
  and the functor $\Phi$ are braided (resp.\ symmetric).
\end{defn}

We will also speak of `monoidal $*$-fibrations' and `monoidal
$*$-bifibrations', but without implying any compatibility between the
monoidal structure and the right adjoints $f_*$.  This is because in
most cases there is no such compatibility.

\begin{eg}
  Let \sC\ be a category with finite limits.  Recall that if $\sC^\dn$
  denotes the category of arrows in \sC, the codomain functor gives a
  bifibration $\ttArr_\sC\maps \sC^\dn\to\sC$ called the
  `self-indexing' of \sC.  It is easy to see that $\ttArr_\sC$ is a
  monoidal bifibration when \sC\ and $\sC^\dn$ are equipped with their
  cartesian products.
\end{eg}

\begin{eg}\xlabel{eg:LR-monfib}
  If \bbD\ is a monoidal framed bicategory, then $(L,R)\maps
  \bbDo\to\bbDz\times\bbDz$ is a monoidal bifibration.  If \bbD\ is
  braided or symmetric, so is $(L,R)$.
\end{eg}

\begin{eg}\xlabel{eg:mod-monfib}
  The fibration $\ttMod\maps \mathbf{Mod}\to\mathbf{Ring}$ is a
  monoidal $*$-bifibration under the tensor product of rings and the
  `external' tensor product of modules.
\end{eg}

For most of our applications, such as \xref{locmon->monfib} below and
the construction of framed bicategories in
\S\ref{sec:monfib-to-framed}, we will require the base category \sB\
to be cartesian or cocartesian monoidal.  However, we see from
Examples~\ref{eg:LR-monfib} and~\ref{eg:mod-monfib} that this is not
always the case, and the general notion of monoidal fibration is
interesting in its own right.

Recall from \xref{fib-psfr-eqv} that the 2-category of fibrations
$\Phi\maps \sA\to\sB$ is equivalent to the 2-category of
pseudofunctors $\sB\op\to\Cat$.  We intend to prove an analogous
result for monoidal fibrations over cartesian base categories, but
first we must define the 2-category of monoidal fibrations.

\begin{defn}\xlabel{def:monoidal-mor-fib}
  Let $\Phi\maps \sA\to\sB$ and $\Phi'\maps \sA'\to\sB'$ be monoidal
  fibrations.
  \begin{itemize}
  \item An \textbf{oplax monoidal morphism of fibrations} is a commuting
    square
    \begin{equation}
      \xymatrix{\sA' \ar[r]^{F_1}\ar[d]_{\Phi'} & \sA \ar[d]^\Phi\\
        \sB' \ar[r]_{F_0} & \sB}\label{eq:mon-mor-fib}
    \end{equation}
    (that is, an oplax morphism of fibrations) together with the data
    of oplax monoidal functors on $F_0$ and $F_1$ such that the
    identity $\Phi F_1=F_0\Phi'$ is a monoidal natural transformation.
  \item An oplax morphism is \textbf{strong} if $F_0$ and $F_1$ are
    strong monoidal functors and $F_1$ preserves cartesian arrows.
  \item A \textbf{lax} morphism is a square~(\ref{eq:mon-mor-fib}) such
    that $F_0$ and $F_1$ are lax monoidal functors, $F_1$ preserves
    cartesian arrows, and the equality $\Phi F_1=F_0\Phi'$ is a
    monoidal transformation.
  \end{itemize}
  Any sort of morphism is \textbf{over \sB} if $F_0$ is an identity
  $\sB'=\sB$.  If $\Phi$ and $\Phi'$ are braided (resp.\ symmetric),
  then any sort of monoidal morphism is braided (resp.\ symmetric) if
  the functors $F_0$ and $F_1$ and the equality $\Phi F_1=F_0\Phi'$
  are braided (resp.\ symmetric).

  If $\Phi$ and $\Phi'$ are monoidal bifibrations, then a \textbf{lax
    monoidal morphism of bifibrations} is just a lax monoidal morphism
  of fibrations, while an \textbf{oplax} (resp.\ \textbf{strong})
  \textbf{monoidal morphism of bifibrations} is an oplax (resp.\ strong)
  monoidal morphism of fibrations which also preserves opcartesian
  arrows.

  A \textbf{monoidal transformation of fibrations}, or of bifibrations,
  is a transformation of fibrations whose components are monoidal
  natural transformations.  If the two morphisms are over \sB, then
  the transformation is \textbf{over \sB} if its downstairs component is
  an identity.
\end{defn}

\begin{notns}
  Let \MFc\ (resp.\ \MF, \MFl) be the 2-category of monoidal
  fibrations, oplax (resp.\ strong, lax) monoidal morphisms of
  fibrations, and monoidal transformations of fibrations.  We write
  \BMF\ and \SMF\ for the braided and symmetric versions.  Let $\MFB$
  denote the sub-2-category of \MF\ consisting of fibrations,
  morphisms, and transformations over \sB, and so on.  Finally, we
  write \MonCat\ for the 2-category of monoidal categories, strong
  monoidal functors, and monoidal natural transformations, and
  similarly \BrMonCat\ and \SymMonCat.
\end{notns}

\begin{thm}\xlabel{locmon->monfib}
  If \sB\ is cartesian monoidal, the equivalence of
  \xref{fib-psfr-eqv} lifts to equivalences of 2-categories
  \begin{align*}
    \MFB &\eqv [\sB\op,\MonCat]\\
    \BMFB &\eqv [\sB\op,\BrMonCat]\\
    \SMFB &\eqv [\sB\op,\SymMonCat].
  \end{align*}
\end{thm}

This means that, in particular, in a monoidal fibration with cartesian
base, each fiber is monoidal and each transition functor $f^*$ is
strong monoidal.  We call the monoidal structure on \sA\ the
\textbf{external} monoidal structure, and the monoidal structures on
fibers the \textbf{internal} monoidal structures.

In many cases, the internal monoidal structures on the fibers are more
familiar and predate the external monoidal structure.  For example, in
$\ttArr_\sC$, the fiber over $B$ is the slice category $\sC/B$, and
the internal monoidal structure is the fiber product over $B$.

It is crucial that \sB\ be \emph{cartesian} monoidal for
\xref{locmon->monfib} to be true.  For example, the fiber of \ttMod\
over a noncommutative ring $R$ is the category $\mathbf{Mod}_R$ of
$R$-modules, which does not in general have an internal tensor
product.  But if we restrict to the monoidal fibration $\ttCMod$ of
modules over \emph{commutative} rings, the tensor product in
$\mathbf{CRing}$ becomes the coproduct, so we can apply the dual
result, obtaining the familiar tensor product on $\mathbf{Mod}_R$ in
the commutative case.

\begin{notn}
  In a cartesian monoidal category \sB, we write $\pr_B$ for any map
  which projects $B$ out of a product; thus we have $\pr_B\maps B\to
  1$, but also $\pr_B\maps A\times B\times C \to A\times C$.  We also
  write $\Delta_B\maps B\to B\times B$ for the diagonal, and other maps
  constructed from it such as $A\times B\times C \to A\times B\times
  B\times C$.
\end{notn}

\begin{proof}[\xref{locmon->monfib}]
  Let $\Phi\maps \sA\to\sB$ be a monoidal fibration with a chosen
  cleavage, and let $B\in\sB$.  We define a monoidal structure on the
  fiber $\sA_B$ as follows.  The unit object is $I_B=\pr_B^*I$, and the
  product is given by
  \begin{equation}
    M\boxtimes N = \Delta_B^*(M\ten N)\label{eq:external->internal}
  \end{equation}
  where $M,N\in\sA_B$ and $\ten$ is the monoidal structure of \sA.  To
  obtain the associativity isomorphism, we tensor the cartesian arrow
  \[M\boxtimes N \too M\ten N\]
  (which lives over $\Delta_B$) with $Q$ to get an arrow
  \[(M\boxtimes N)\ten Q \too (M\ten N)\ten Q
  \]
  which is cartesian since $\ten$ preserves cartesian arrows.  We then
  compose with another cartesian arrow over $\Delta_B$ to obtain a
  composite cartesian arrow
  \[(M\boxtimes N)\boxtimes Q \too (M\ten N)\ten Q.
  \]
  We do the same on the other side to get a cartesian arrow
  \[M\boxtimes(N\boxtimes Q) \too M\ten(N\ten Q)\]
  and the unique factorization of
  \[\fa\maps (M\ten N)\ten Q\iso M\ten(N\ten Q)\]
  through these cartesian arrows gives an associativity isomorphism
  \[(M\boxtimes N)\boxtimes Q \iso M\boxtimes(N\boxtimes Q).
  \]
  for $\sA_B$.  The pentagon axiom follows from unique factorization
  through cartesian arrows and the pentagon axiom for \sA.  The unit
  constraints and axioms are analogous, using the fact that
  $\pr_B\Delta_B=1_B$, as is the braiding when $\Phi$ is braided or
  symmetric.

  Now consider a map $f\maps A\to B$; we show that $f^*$ is strong
  monoidal.  We have the composite cartesian arrows
  \[f^*M\boxtimes f^*N \too f^*M\ten f^*N \too M\ten N\]
  and
  \[f^*(M\boxtimes N) \too M\boxtimes N \too M\ten N,\]
  both lying over $\Delta_B f = (f\times f)\Delta_A$; hence we obtain
  a canonical isomorphism 
  \[f^*M\boxtimes f^*N \iso f^*(M\boxtimes N).
  \]
  The unit constraint is similar and, as before, the coherence of
  these constraints follows from the uniqueness of factorization
  through cartesian arrows, as does the fact that the isomorphisms
  $(fg)^*\iso f^*g^*$ and $(1_B)^*\iso \Id$ are monoidal.  Therefore,
  we have constructed a pseudofunctor $\sB\op\to \MonCat$ from a
  monoidal fibration.  It is straightforward to extend this
  construction to give 2-functors
  \begin{align*}
    \MFB &\too{} [\sB\op,\MonCat]\\
    \BMFB &\too{} [\sB\op,\BrMonCat]\\
    \SMFB &\too{} [\sB\op,\SymMonCat].
  \end{align*}
  Uniqueness of factorization again gives the coherence to show that
  the resulting pseudonatural transformations are pointwise monoidal.
  
  Conversely, given a pseudofunctor $\sB\op\to\MonCat$, we define a
  fibration over \sB\ in the usual way, and define an external
  product as follows: given $M,N$ over $A,B$ respectively, let
  \begin{equation}
    M\ten N = \pr_B^*M\boxtimes \pr_A^*N.\label{eq:internal->external}
  \end{equation}
  The external unit is $I_1$, the internal unit in the fiber over
  $1$.  For an associativity isomorphism we use
  \begin{align*}
    (M\ten N)\ten Q
    &= \pr_C^*(\pr_B^*M\boxtimes \pr_A^*N)\boxtimes\pr_{AB}^*Q\\
    &\iso (\pr_{BC}^* M \boxtimes \pr_{AC}^*N)\boxtimes \pr_{AB}^*Q\\
    &\iso \pr_{BC}^* M \boxtimes (\pr_{AC}^*N\boxtimes \pr_{AB}^*Q)\\
    &\iso \pr_{BC}^*M \boxtimes \pr_A^*(\pr_C^*N\boxtimes \pr_B^*Q)\\
    &= M\ten (N\ten Q)
  \end{align*}
  using the monoidal constraints for the strong monoidal functors
  $\pr^*$, the composition constraints for the pseudofunctor, and the
  associativity for the internal products.  It is straightforward, if
  tedious, to check that this isomorphism satisfies the pentagon
  axiom.  Similarly, we have a unit constraint
  \begin{align*}
    M\ten I_1
    &= \pr_1^* M \boxtimes \pr_A^* I_1\\
    &\iso M \boxtimes I_A\\
    &\iso M
  \end{align*}
  which can be checked to be coherent; thus \sA\ is monoidal, and
  $\Phi$ is strict monoidal by definition.  It is obvious how to
  define a braiding in the braided or symmetric case making \sA\ and
  $\Phi$ braided or symmetric.  Finally, using the composition
  constraints and monoidal constraints, we have:
  \begin{align*}
    f^*M \ten g^*N
    &= \pr^* f^* M \boxtimes \pr^* g^* N\\
    &\iso (f\times g)^* \pr^* M\boxtimes (f\times g)^* \pr^* N\\
    &\iso (f\times g)^* \big(\pr^*M\boxtimes \pr^*N\big)\\
    &= (f\times g)^*(M\ten N),
  \end{align*}
  which we can then use to verify that $\ten$ preserves cartesian
  arrows.  Thus we have constructed a monoidal fibration of the
  desired type.  It is straightforward to extend this to a 2-functor
  and verify that these constructions are inverse equivalences.
\end{proof}

\begin{rmk}
  Under the above equivalence, pseudofunctors which land in
  \emph{cartesian} monoidal categories correspond to fibrations where
  the total category \sA\ is cartesian monoidal.
\end{rmk}

We end this section by introducing a few new examples of monoidal
fibrations.

\begin{eg}\xlabel{eg:param-obj}
  Let \sC\ be a category with finite limits and colimits, and assume
  that pullbacks in \sC\ preserve finite colimits.  (For example, \sC\
  could be locally cartesian closed.)  Let $\mathrm{Retr}(\sC)$ be the
  category of \emph{retractions} in \sC.  That is, an object of
  $\mathrm{Retr}(\sC)$ is a pair of maps $A\too[s] X \too[r] A$
  such that $rs=1_A$.  This is also known as an object $X$
  `parametrized' over $A$, in which case $s$ is called the `section'.
  We define $\ttRetr_\sC\maps \mathrm{Retr}(\sC)\to\sC$ to take the
  above retraction to $A$.  It is easy to check that pullback and
  pushout make $\Phi$ into a bifibration, which is a $*$-bifibration
  if \sC\ is locally cartesian closed.

  The fiber over $B\in \sC$ is the category $\sC_B$ of objects
  parametrized over $B$.  It has finite products, given by pullback
  over $B$, but usually the relevant monoidal structure is not the
  cartesian product but the \emph{fiberwise smash product}, defined as
  the pushout
  \[\xymatrix{X\sqcup_B Y \ar[r]\ar[d] & X\times_B Y \ar[d]\\
    B \ar[r] & X \sm_B Y.}\]
  The unit is $B\sqcup B \to B$ with section given by one of the
  coprojections.  Under the assumption that pullbacks preserve finite
  colimits, this defines a symmetric monoidal structure on $\sC_B$,
  all the functors $f^*$ are strong symmetric monoidal, and the
  coherence isomorphisms are also monoidal.  Thus by
  \xref{locmon->monfib}, $\ttRetr_\sC$ is a symmetric monoidal
  fibration, and it is easy to check that it is actually a monoidal
  bifibration.  The external monoidal structure on
  $\mathrm{Retr}(\sC)$ is called the \emph{external smash product}
  $\exsm$.
\end{eg}

\begin{eg}\xlabel{eg:param-spaces}
  Suppose that \sC\ has finite limits and colimits, and \emph{not} all
  pullbacks preserve finite colimits, but there is some full
  subcategory \sB\ of \sC\ such that pullbacks along morphisms in \sB\
  do preserve finite colimits.  Then we can repeat the construction of
  \xref{eg:param-obj} using parametrized objects whose base objects
  are restricted to lie in \sB.  This is what is done
  in~\cite[\S2.5]{pht}, with $\sC=\sK$ the category of $k$-spaces and
  $\sB=\sU$ the category of compactly generated spaces.  By a slight
  abuse of notation, we call the resulting monoidal $*$-bifibration
  $\ttRetr_\Top$, since we have only been prevented from considering
  all retractions in \Top\ by point-set technicalities.  The objects
  of $\mathrm{Retr}(\Top)$ are called \emph{ex-spaces}.
\end{eg}

\begin{eg}\xlabel{eg:spectra-ps}
  For each space $B\in\sU$, a category $\sS_B$ of orthogonal spectra
  parametrized over $B$ is defined in~\cite[Ch.~11]{pht}.  A map
  $f\maps A\to B$ of spaces gives rise to a string of adjoints
  $f_!\adj f^*\adj f_*$ which are pseudofunctorial in $f$.  Each
  category $\sS_B$ is closed symmetric monoidal under an internal
  smash product $\sm_B$, each functor $f^*$ is closed symmetric
  monoidal, and so are the composition constraints.  Thus, by
  \xref{locmon->monfib}, we obtain a symmetric monoidal fibration
  which we denote \ttSp.  The external smash product $\exsm$ is
  defined in~\cite[11.4.10]{pht} just as we have done
  in~(\ref{eq:internal->external}).

  To show that \ttSp\ is in fact a monoidal $*$-bifibration, one can
  check directly that $\exsm$ preserves opcartesian arrows.  However,
  this will also follow from \xref{extclosed-mates} below.
\end{eg}

\begin{eg}\xlabel{eg:spectra-derived}
  Let $\sB=\sU$ as in \xref{eg:spectra-ps}, but instead of $\sS_B$ we
  use its homotopy category $\Ho\sS_B$.  It is proven
  in~\cite[12.6.7]{pht} that $f_!\adj f^*$ is a Quillen adjunction,
  for a suitable choice of model structures on $\sS_B$, hence it
  descends to an adjunction on homotopy categories which is still
  pseudofunctorial; thus we obtain another functor $\Ho(\ttSp)\maps
  \sA\to\sB$ which is a bifibration.

  The external smash product $\exsm$ is proven to be a Quillen left
  adjoint in~\cite[12.6.6]{pht}; thus it descends to homotopy
  categories to make $\sA$ symmetric monoidal.
  Since~\cite[13.7.2]{pht} shows that $\exsm$ preserves cartesian
  arrows, $\Ho(\ttSp)$ is a symmetric monoidal fibration, and the same
  methods as in \xref{eg:spectra-ps} show that it is a monoidal
  bifibration.  The derived functors $f^*$ also have right adjoints,
  although these are constructed in~\cite[13.1.18]{pht} using Brown
  representability rather than by deriving the point-set level right
  adjoints; thus $\Ho(\ttSp)$ is a monoidal $*$-bifibration.
\end{eg}

\section{Closed monoidal fibrations}
\label{sec:closed-monfib}

We now consider two different notions of when a monoidal fibration is
`closed'.  To fix terminology and notation, we say an ordinary
monoidal category \sC\ with product $\boxtimes$ is \emph{closed} if
the functors $(M\boxtimes -)$ and $(-\boxtimes N)$ have right adjoints
$(-\LHD M)$ and $(N\RHD -)$, respectively, for all $M,N$.  Of course, if
\sC\ is symmetric, then $P\LHD M\iso M\RHD P$.  If \sC\ and \sD\ are
closed monoidal categories and $f^*\maps \sC\to\sD$ is a strong
monoidal functor, then there are canonical natural transformations
\begin{align}
  f^*(N\RHD P) &\too f^*N \RHD f^* P\label{eq:intcl-map-1}\\
  f^*(P\LHD N) &\too f^*P \LHD f^* N.\label{eq:intcl-map-2}
\end{align}
When these transformations are isomorphisms, we say that $f^*$ is
\emph{closed monoidal}.  Of course, in the symmetric
case,~(\ref{eq:intcl-map-1}) is an isomorphism if and only
if~(\ref{eq:intcl-map-2}) is.

\begin{defn}\xlabel{def:closed-mf}
  Let $\Phi\maps \sA\to\sB$ be a monoidal fibration where \sB\ is
  cartesian monoidal (so that each fiber is a monoidal category).  We
  say $\Phi$ is \textbf{internally closed} if each fiber $\sA_B$ is
  closed monoidal and each functor $f^*$ is closed monoidal.
\end{defn}

However, in any monoidal fibration, we can also ask whether the
external product
\[\ten\maps \sA_A\times \sA_B \to \sA_{A\ten B}\]
has adjoints $\xlhd,\xrhd$, with defining isomorphisms
\[\sA_{A\ten B}(M\ten N, P) \iso \sA_{A}(M, N\xrhd P) \iso \sA_B (N, P\xlhd M).\]
If so, then for any $f\maps C\to A$ and $g\maps D\to B$ there are
canonical transformations
\begin{align}
  f^*(N\xrhd P) &\too N \xrhd (f\ten 1)^*P\label{eq:extcl-map-1}\\
  g^*(P\xlhd M) &\too (1\ten g)^*P\xlhd M\label{eq:extcl-map-2}
\end{align}
defined analogously to~(\ref{eq:intcl-map-1})
and~(\ref{eq:intcl-map-2}).  For example,~\eqref{eq:extcl-map-1} is
the adjunct of the composite
\[f^*(N\xrhd P) \ten N \too[\iso] (f\ten 1)^*\big((N\xrhd P)\ten N\big) \too (f\ten 1)^*P.
\]

\begin{defn}
  Let $\Phi\maps \sA\to\sB$ be a monoidal fibration.  We say that
  $\Phi$ is \textbf{externally closed} if the adjoints $\xlhd,\xrhd$
  exist and the maps~(\ref{eq:extcl-map-1}) and~(\ref{eq:extcl-map-2})
  are isomorphisms for all $f,g$.
\end{defn}

\begin{egs}
  If \sC\ is locally cartesian closed, then $\ttArr_\sC$ is internally
  and externally closed.  If \sC\ also has finite colimits, then
  $\ttRetr_\sC$ is internally and externally closed.
\end{egs}

\begin{eg}
  The fibration \ttSp\ of parametrized orthogonal spectra over spaces
  is internally and externally closed; its internal homs are defined
  in~\cite[11.2.5]{pht} and the base change functors are shown to be
  closed in~\cite[11.4.1]{pht}.  We postpone consideration of
  $\Ho(\ttSp)$ until later.
\end{eg}

\begin{eg}
  The fibration $\ttMod\maps \mathbf{Mod}\to\mathbf{Ring}$ is
  externally closed.  If $N$ is a $B$-module and $P$ is an $A\ten
  B$-module, the external-hom $N\xrhd P$ is $\Hom_B(N,P)$, which
  retains the $A$-module structure from $P$.  In this case, internal
  closure makes no sense because the fibers are not even monoidal.
\end{eg}

\begin{eg}
  The monoidal $*$-bifibration \ttCMod\ of modules over
  \emph{commutative} rings is also externally closed.  In this case
  the fibers $\mathbf{Mod}_R$ are closed monoidal, but neither $f_!$
  nor $f^*$ is a closed monoidal functor.
\end{eg}

\begin{eg}
  If \bbD\ is a monoidal framed bicategory, then the monoidal
  bifibration $(L,R)$ is externally closed just when \bbD\ is
  externally closed in the sense of \S\ref{sec:monoidal-structures}.
  The fact that~(\ref{eq:extcl-map-1}) and~(\ref{eq:extcl-map-2}) are
  isomorphisms in this case will follow from \xref{extclosed-mates},
  below.
\end{eg}

\begin{rmk}
  Contrary to what one might expect (see, for
  example,~\cite[\S2.4]{pht}), external closedness does not imply that
  the monoidal category \sA\ is closed in its own right.  For one
  thing, $N\xrhd P$ is only defined when $N\in\sA_B$ and $P\in
  \sA_{A\times B}$.  But even when defined, $N\xrhd P$ is not an
  internal-hom for \sA: if $M\in\sA_C$, then the morphisms $M\to
  N\xrhd P$ in \sA\ are bijective not to all morphisms $M\ten N\to P$,
  but only those lying over $f\times 1$ for some $f\maps C\to A$.
\end{rmk}

In the rest of this section, we will prove that under mild hypotheses,
internal and external closedness are equivalent, and give useful dual
versions of the maps (\ref{eq:intcl-map-1}), (\ref{eq:intcl-map-2}),
(\ref{eq:extcl-map-1}), and (\ref{eq:extcl-map-2}).  We begin by
comparing the internal and external homs.

\begin{prop}\xlabel{clfib<->extadj}
  Let $\Phi$ be either
  \begin{enumerate}
  \item a monoidal $*$-fibration in which \sB\ is cartesian monoidal, or
  \item a monoidal bifibration in which \sB\ is cocartesian monoidal.
  \end{enumerate}
  Then the right adjoints $\xlhd,\xrhd$ exist if and only if $\Phi$
  has closed fibers (i.e.\ the right adjoints $\LHD,\RHD$ exist).
\end{prop}
\begin{proof}
  Suppose first that \sB\ is cartesian and each fiber is closed.  Then
  for $N\in\sA_B$ and $Q\in \sA_{A\times B}$ we define
  \begin{equation}
    N\xrhd Q = (\pr_B)_*\big(\pr_A^*N \RHD Q\big).\label{eq:extcl->loccl}
  \end{equation}
  and similarly for $\xlhd$.  Conversely, if $\xlhd,\xrhd$ exist, then
  for $N,Q\in\sA_A$ we define
  \begin{equation}
    N\RHD Q = N \xrhd ((\Delta_A)_* Q)\label{eq:loccl->extcl}
  \end{equation}
  and similarly for $\LHD$.  It is easy to check, using the
  relationships between $\ten$ and $\boxtimes$ established in
  \xref{locmon->monfib}, that these definitions suffice.

  In the cocartesian case, these relationships become
  \begin{align*}
    M\boxtimes N &\iso \nabla_!(M\ten N)\\
    M\ten N &\iso \eta_!M\boxtimes \eta_!N,
  \end{align*}
  where $\nabla_A\maps A\sqcup A\to A$ denotes the `fold' or
  codiagonal, and $\eta\maps \emptyset\to A$ is the unique map from
  the initial object.  Therefore, the analogous definitions:
  \begin{align*}
    M\RHD N &= M \xrhd (\nabla^* N)\\
    M\xrhd N &= \eta^*\big(\eta_!M \RHD N\big).
  \end{align*}
  allow us to pass back and forth between internal and external
  closedness.
\end{proof}

This equivalence is valuable because sometimes one of the two types of
right adjoints is much easier to construct than the other.

\begin{eg}
  The homotopy-level fibration $\Ho(\ttSp)$ has the adjoints
  $\xlhd,\xrhd$, since the adjunction between $\exsm$ and $\xrhd$ in
  \ttSp\ is Quillen (see~\cite[12.6.6]{pht}).  This then implies, by
  \xref{clfib<->extadj}, that the fibers of $\Ho(\ttSp)$ are all
  closed monoidal.  This would be difficult to prove directly, since
  we have no homotopical control over the internal monoidal structures
  in \ttSp.
\end{eg}

In order to prove a full equivalence of local and external closedness,
we need to assume an extra condition on the commutativity of right and
left adjoints.  Suppose that $\Phi\maps \sA\to\sB$ is a fibration and
that the square
\begin{equation}
  \xymatrix{A \ar[r]^h\ar[d]_k & B \ar[d]^g\\ C \ar[r]_f & D}\label{eq:bchv-downstairs}
\end{equation}
commutes in \sB.  Thus we obtain a square
\begin{equation}
  \xymatrix{\sA_A \ar@{}[dr]|{\iso} & \sA_B \ar[l]_{h^*} \\
    \sA_C\ar[u]^{k^*} & \sA_D \ar[l]^{f^*} \ar[u]_{g^*}
  }\label{eq:bchv-upstairs}
\end{equation}
which commutes up to canonical isomorphism.  If $\Phi$ is a
bifibration, there is a canonical natural transformation
\begin{equation}
  k_!h^*\too f^*g_!,\label{eq:bchv-trans}
\end{equation}
namely the `mate' of the isomorphism~(\ref{eq:bchv-upstairs}).
Explicitly, it is the composite
\[k_!h^* \too[\eta] k_!h^*g^*g_! \iso k_!k^*f^*g_! \too[\ep] f^*g_!.
\]
Similarly, if $\Phi$ is a $*$-fibration, there is a canonical
transformation
\begin{equation}
  g^*f_* \too h_*k^*.\label{eq:bchv-trans-other}
\end{equation}

\begin{defn}
  If $\Phi$ is a bifibration (resp.\ a $*$-fibration), we say that the
  square~(\ref{eq:bchv-downstairs}) satisfies the \textbf{Beck-Chevalley
    condition} if the natural transformation~(\ref{eq:bchv-trans})
  (resp.~(\ref{eq:bchv-trans-other})) is an isomorphism.  We say that
  $\Phi$ is \textbf{strongly BC} if this condition is satisfied by every
  pullback square, and \textbf{weakly BC} if it is satisfied by every
  pullback square in which one of the legs ($f$ or $g$, above) is a
  product projection.  If instead all \emph{pushout} squares satisfy
  the Beck-Chevalley condition, we say that $\Phi$ is \textbf{strongly
    co-BC}.
\end{defn}

If $\Phi$ is a $*$-bifibration, then~(\ref{eq:bchv-trans})
and~(\ref{eq:bchv-trans-other}) are mates under the composite
adjunctions $f^*g_!\adj g^*f_*$ and $k_!h^*\adj h_*k^*$, so that one
is an isomorphism if and only if the other is.  Thus, a
$*$-bifibration is strongly or weakly BC as a bifibration if and only
if it is so as a $*$-fibration.

\begin{egs}
  The monoidal bifibrations $\ttArr_\sC$ and $\ttRetr_\sC$ are always
  strongly BC, as is the monoidal $*$-bifibration \ttSp\
  (see~\cite[11.4.8]{pht}).
\end{egs}

\begin{eg}
  The monoidal $*$-bifibration \ttCMod, whose base is cocartesian
  monoidal, is strongly co-BC.
\end{eg}

\begin{eg}
  The homotopy-level monoidal $*$-bifibration $\Ho(\ttSp)$ is only
  weakly BC; it is proven in~\cite[13.7.7]{pht} that the
  Beck-Chevalley condition is satisfied for pullback squares one of
  whose legs is a fibration in the topological sense (which includes
  product projections, of course).  It does \emph{not} satisfy the
  Beck-Chevalley condition for arbitrary pullback squares; a concrete
  counterexample is given in~\cite[0.0.1]{pht}.  One intuitive reason
  for this is that since \ttSp\ also incorporates `homotopical'
  information about the base spaces, we should only expect the derived
  operations to be well-behaved on \emph{homotopy} pullback squares.
  This is our main motivation for introducing the notion of `weakly
  BC'.
\end{eg}

Of course, the idea of commuting adjoints is older than the term
`Beck-Chevalley condition'.  In the theory of fibered categories, what
we call a `strongly BC bifibration' is referred to as a `fibration
with indexed coproducts'.

We will eventually use Beck-Chevalley conditions in our construction
of a framed bicategory from a monoidal fibration
(\xref{monfib->framed}), but we mention them in this section for the
purposes of the following result.

\begin{prop}\xlabel{closed<->extclosed}
  Let $\Phi\maps \sA\to\sB$ be a monoidal $*$-fibration in which \sB\
  is cartesian monoidal.  Then
  \begin{enumerate}
  \item if $\Phi$ is internally closed and weakly BC, then it is
    externally closed, and\label{item:cl->extcl}
  \item if $\Phi$ is externally closed and strongly BC, then it is
    internally closed.\label{item:extcl->cl}
  \end{enumerate}
  In particular, a strongly BC monoidal $*$-fibration is internally
  closed if and only if it is externally closed.
\end{prop}
\begin{proofsk}
  Under the equivalences~(\ref{eq:extcl->loccl})
  and~(\ref{eq:loccl->extcl}), each of the maps~(\ref{eq:intcl-map-1})
  and~(\ref{eq:extcl-map-1}) is equal to the composite of the other
  with a Beck-Chevalley transformation, and similarly
  for~(\ref{eq:intcl-map-2}) and~(\ref{eq:extcl-map-2}).  It turns out
  that both of these transformations come from pullback squares; thus
  since $\pr$ appears in~(\ref{eq:extcl->loccl}) but $\Delta$ appears
  in~(\ref{eq:loccl->extcl}), the weak condition is good enough in one
  case but not the other.
\end{proofsk}

Now, if $f^*$ is strong monoidal and has a left adjoint $f_!$, there
is a canonical map
\begin{equation}
  f_!(M\boxtimes f^*N) \too f_!M\boxtimes N.\label{eq:intcl-mate-1}
\end{equation}
When the monoidal categories in question are closed, this is the mate
of~(\ref{eq:intcl-map-1}), so one is an isomorphism if and only if the
other is.  In particular, if $\Phi$ is a monoidal bifibration with
cartesian monoidal base and closed fibers, then $\Phi$ is internally
closed if and only if the maps~(\ref{eq:intcl-mate-1}), together with
the analogous maps
\begin{equation}
  f_!(f^*N\boxtimes M) \too N\boxtimes f_!M,\label{eq:intcl-mate-2}
\end{equation}
are all isomorphisms.  This dual condition is sometimes easier to
check.

\begin{eg}
  Topological arguments involving excellent prespectra are used in
  \cite[13.7.6]{pht} to show that the derived
  maps~(\ref{eq:intcl-mate-1}) are isomorphisms, and therefore
  $\Ho(\ttSp)$ is internally closed.  Since it is weakly BC, we can
  then conclude, by \xref{closed<->extclosed}\ref{item:cl->extcl},
  that it is externally closed as well.
\end{eg}

In a similar way, if $\Phi$ is any monoidal bifibration with right
adjoints $\xlhd,\xrhd$, then~(\ref{eq:extcl-map-1}) has a mate
\begin{equation}
  (f\times 1)_!(M\ten N) \too f_!M\ten N\label{eq:extcl-mate-1}
\end{equation}
which is an isomorphism if and only if~(\ref{eq:extcl-map-1}) is.
But~(\ref{eq:extcl-mate-1}) is an isomorphism just when $-\ten N$
preserves the opcartesian arrow $M\to f_!M$, so we have the following.

\begin{prop}\xlabel{extclosed-mates}
  Let $\Phi$ be a monoidal fibration which is also an opfibration and
  such that the right adjoints $\xlhd,\xrhd$ exist.  Then $\ten$
  preserves opcartesian arrows (that is, $\Phi$ is a monoidal
  bifibration) if and only if $\Phi$ is externally closed.
\end{prop}

\begin{eg}
  As remarked earlier, this implies that a monoidal framed bicategory
  \bbD\ is externally closed in the sense of
  \S\ref{sec:monoidal-structures} if and only if the monoidal
  bifibration $(L,R)$ is externally closed in the sense of this
  section.
\end{eg}

\begin{eg}
  In the converse direction, \xref{extclosed-mates} can be used to
  show that $\Ho(\ttSp)$ and \ttSp\ are monoidal bifibrations, since
  we know that they are externally closed.  This could also be shown
  directly.
\end{eg}

\begin{cor}\xlabel{monfib-closed->loc-and-ext}
  Let $\Phi$ be a strongly BC monoidal $*$-bifibration over a
  cartesian base and having closed fibers.  Then $\Phi$ is internally
  and externally closed.
\end{cor}
\begin{proof}
  Since $\Phi$ is a $*$-fibration, by \xref{clfib<->extadj} it also
  has right adjoints $\xlhd,\xrhd$.  Then, since it is a monoidal
  bifibration, it is externally closed by \xref{extclosed-mates}.  But
  since it is strongly BC,
  \xref{closed<->extclosed}\ref{item:extcl->cl} then implies that it
  is also internally closed.
\end{proof}

\section{From fibrations to framed bicategories}
\label{sec:monfib-to-framed}

We now prove that any well-behaved monoidal bifibration gives rise to
a framed bicategory.  The reader may not be too surprised that there
is some relationship, since many of our examples of monoidal
bifibrations look very similar to our examples of framed bicategories.
In this section we state our results; the proofs will be given in
\S\S\ref{sec:beck-chev}--\ref{sec:monfib-to-framed-ii} after we
consider an important class of examples in \S\ref{sec:monoids-monfib}.

To motivate the precise construction, consider the relationship
between the framed bicategory $\Span(\sC)$ and the monoidal
bifibration $\ttArr_\sC\maps \sC^\dn\to\sC$.  A horizontal 1-cell
$M\maps A\hto B$ in $\Span(\sC)$ is a span $A\leftarrow M\to B$, which
can also be considered as an arrow $M\to A\times B$, and hence an
object of $\sC^\dn$ over $A\times B$.  The horizontal composition of
$M\maps A\hto B$ and $N\maps B\hto C$ is given by pulling back along
the maps to $B$, then remembering only the maps to $A$ and $C$:
\[\xymatrix@R=1pc@C=1pc{&& M\times_B N \ar[dl]\ar[dr]\\
  & M  \ar[dl]\ar[dr] && N  \ar[dl]\ar[dr]\\
  A && B && C}
\]
But this can also be phrased in terms of the maps $M\to A\times B$
and $N\to B\times C$ by taking the product map
\[M\times N \to A\times B\times B\times C,\]
pulling back along the diagonal $\Delta_B$:
\[\xymatrix{M\times_B N \ar[r]\ar[d] & M\times N \ar[d]\\
  A\times B\times C\ar[r] & A\times B\times B\times C}
\]
and then composing with the projection $\pr_B\maps A\times B\times
C\to A\times C$.  In terms of the monoidal bifibration $\ttArr_\sC$,
this can be written as
\begin{equation}
  M\times_B N = (\pr_B)_! \Delta_B^*(M\times N).\label{eq:span-mf->anch-comp}
\end{equation}

Similarly, the unit object $U_A$ in $\Span(\sC)$ is the span $A
\leftarrow A \to A$, alternatively viewed as the diagonal map $A\to
A\times A$.  This can be obtained (in a somewhat perverse way) by
pulling back the terminal object $1$ along the map $\pr_A\maps A\to
1$, then composing with the diagonal $\Delta_A\maps A\to A\times A$.
In the language of $\ttArr_\sC$, we have
\begin{equation}
  U_A = (\Delta_A)_!\pr_A^* 1.\label{eq:span-mf->anch-unit}
\end{equation}
We now observe that the expressions~(\ref{eq:span-mf->anch-comp})
and~(\ref{eq:span-mf->anch-unit}) can easily be generalized to any
monoidal bifibration in which the base is cartesian monoidal, so that
we have diagonals and projections.  This may help to motivate the
following result.

\begin{defn}
  We say that a monoidal bifibration $\Phi\maps \sA\to\sB$ is
  \textbf{frameable} if \sB\ is cartesian monoidal and $\Phi$ is either
  \begin{enumerate}
  \item strongly BC or\label{item:frameable-strong}
  \item weakly BC and internally closed.\label{item:frameable-closed}
  \end{enumerate}
\end{defn}

\begin{thm}\xlabel{monfib->framed}
  Let $\Phi\maps \sA\to\sB$ be a frameable monoidal bifibration.  Then
  there is a framed bicategory $\bbFr(\Phi)$ with a vertically strict
  involution, defined as follows.
  \begin{enumerate}
  \item $\bbFr(\Phi)_0=\sB$.
  \item $\bbFr(\Phi)_1$, $L$, and $R$ are defined by the following
    pullback square.
    \[\xymatrix{\bbFr(\Phi)_1 \ar[r]\ar[d]_{(L,R)} & \sA \ar[d]^\Phi\\
      \sB\times \sB\ar[r]_<>(.5)\times & \sB.}\]
    Thus the horizontal 1-cells $A\hto B$ are the objects of \sA\ over
    $A\times B$, and the 2-cells $M\sto{f}{g} N$ are the arrows of \sA\ over
    $f\times g$.
  \item The horizontal composition of $M\maps A\hto B$ and $N\maps
    B\hto C$ is 
    \[M\odot N = (\pr_B)_!\Delta_B^*(M\ten N),
    \]
    and similarly for 2-cells.
  \item The horizontal unit of $A$ is
    \[U_A = (\Delta_A)_!\pr_A^* I.\]
  \item The involution is the identity on objects and we have $M\op$ =
    $\fs^*M$, where $\fs$ is the symmetry isomorphism.
  \end{enumerate}
  If $\Phi$ is externally closed and a $*$-bifibration, then
  $\bbFr(\Phi)$ is closed in the sense of \S\ref{sec:duality-theory}.
  If $\Phi$ is symmetric, then $\bbFr(\Phi)$ is symmetric monoidal in
  the sense of \S\ref{sec:monoidal-structures} and its involution is
  also symmetric monoidal.
\end{thm}

\begin{egs}
  As alluded to above, if \sC\ has finite limits, the symmetric
  monoidal bifibration $\ttArr_\sC$ gives rise to the symmetric
  monoidal framed bicategory $\Span(\sC)$, which is closed if \sC\ is
  locally cartesian closed.
\end{egs}

\begin{eg}
  If \sC\ has finite limits and colimits preserved by pullback, then
  the monoidal bifibration $\ttRetr_\sC$ gives rise to a symmetric
  monoidal framed bicategory of \emph{parametrized objects}, which we
  denote $\Ex(\sC)$.  It is also closed if \sC\ is locally cartesian
  closed.

  Applied to the monoidal $*$-bifibration $\ttRetr_\Top$ of ex-spaces
  from \xref{eg:param-spaces}, we obtain a framed bicategory
  $\Ex(\Top)$ of parametrized spaces which is both symmetric monoidal
  and closed.
\end{eg}

\begin{eg}
  The monoidal $*$-bifibration \ttSp\ of parametrized orthogonal
  spectra gives rise to a point-set level framed bicategory of
  parametrized spectra, which we may denote \Sp.  It is symmetric
  monoidal and closed.
\end{eg}

\begin{eg}\xlabel{eg:ho-spectra}
  The homotopy-category monoidal $*$-bifibration $\Ho(\ttSp)$, which
  is weakly BC and internally closed, gives rise to a framed
  bicategory $\Ho(\Sp)$.  This is the same as the framed bicategory we
  have been calling \Ex\ ever since \S\ref{sec:double-categories}.
  Similarly, $\Ho(\ttRetr_\Top)$ gives rise to a homotopy-level framed
  bicategory $\Ho(\Ex(\Top))$ of parametrized spaces.  Both of these
  framed bicategories are symmetric monoidal and closed.  These are
  the only ones of our examples which are weakly rather than strongly
  BC.
\end{eg}

The dual version of \xref{monfib->framed} says the following.

\begin{thm}\xlabel{dual-monfib->framed}
  If $\Phi\maps \sA\to\sB$ is a strongly co-BC monoidal bifibration
  where \sB\ is cocartesian monoidal, then there is a framed
  bicategory $\bbFr(\Phi)$ with a vertically strict involution,
  defined as in \xref{monfib->framed}, except that composition is
  given by
  \[M\odot N = \eta^*\nabla_!(M\ten N)
  \]
  and units are given by
  \[U_A = \nabla^* \eta_! I.
  \]
  If $\Phi$ is externally closed and is a $*$-bifibration, then
  $\bbFr(\Phi)$ is closed.  If $\Phi$ is symmetric, then $\bbFr(\Phi)$
  and its involution are symmetric monoidal.
\end{thm}

\begin{eg}
  The monoidal $*$-bifibration \ttCMod\ gives rise to the framed
  bicategory \CMod, which is symmetric monoidal and closed.  However,
  \Mod\ cannot be constructed in this way, because the category of
  noncommutative rings is not cocartesian monoidal.
\end{eg}

Like our other ways of constructing framed bicategories, these results
are 2-functorial.  To state this precisely, we need to define the
right 2-categories.  In \S\ref{sec:monoidal-fibrations} we defined
lax, strong, and oplax monoidal morphisms of bifibrations to be those
that preserve both the monoidal structure and the base change in
appropriate ways.  These morphisms, together with the monoidal
transformations defined there, give us 2-categories $\MbF$, \MbF, and
$\MbFl$.  Let $\MFfr$ denote the full sub-2-category of \MbF\ spanned
by the frameable monoidal bifibrations, and similarly for $\MFfrc$ and
$\MFfrl$.

\begin{thm}\xlabel{monfib->framed-2functor}
  The construction of \xref{monfib->framed} extends to a 2-functor
  \[\bbFr\maps \MFfr \too \FrBi\]
  and similarly for oplax and lax morphisms.
\end{thm}

\begin{eg}
  The 2-functor $\Span\maps \Cart\to\FrBi$ clearly factors through
  \bbFr\ via a 2-functor $\ttArr\maps \Cart\to\MFfr$.
\end{eg}

\begin{eg}
  Let \biCartpsc\ denote the 2-category of categories with finite
  limits and finite colimits preserved by pullback, functors which
  preserve finite limits and colimits, and natural transformations.
  Then we have a 2-functor $\ttRetr\maps \biCartpsc\to\MFfr$.
  Composing this with $\bbFr$ defines a 2-functor $\Ex\maps
  \biCartpsc\to \FrBi$.
\end{eg}

As with our other 2-categories, we automatically obtain notions of
equivalence and adjunction between monoidal bifibrations, and these
are preserved by 2-functors such as \bbFr.  As usual, we can also
characterize these more explicitly; we omit the proof of the
following.

\begin{prop}\xlabel{monfib-adjn}
  An adjunction $F\adj G$ in $\MFfrl$ between $\Phi\maps \sA\to\sB$
  and $\Phi'\maps \sA'\to\sB'$ consists of the following properties
  and structure.
  \begin{enumerate}
  \item $F$ is a strong monoidal morphism of bifibrations and $G$ is a
    lax monoidal morphism of bifibrations;
  \item We have monoidal adjunctions $F_0\maps \sB\toot \sB'\spam G_0$
    and $F_1\maps \sA\toot \sA'\spam G_1$;
  \item We have equalities $\Phi'F_1 = F_0\Phi$ and $\Phi G_1 =
    G_0\Phi'$ which are monoidal transformations; and
  \item The adjunction $F_1\adj G_1$ `lies over' $F_0\adj G_0$ in the
    sense that the following square commutes:
    \[\xymatrix{\sA(M,G_1N) \ar[r]^\iso\ar[d]_\Phi &
      \sA'(F_1M,N) \ar[d]^{\Phi'}\\
      \sB(A,G_0 B)\ar[r]_\iso & \sB'(F_0A,B).}
    \]
  \end{enumerate}
\end{prop}

\begin{rmk}
  In fact, these conditions are somewhat redundant.  For example, left
  adjoints automatically preserve opcartesian arrows and right
  adjoints automatically preserve cartesian ones, and the right
  adjoint of a strong monoidal functor is always lax monoidal.  These
  are consequences of `doctrinal adjunction' (see \xref{ladj-strong}
  and~\cite{kelly:doc-adjn}) and a property called `lax-idempotence'
  (see~\cite{property-like}).
\end{rmk}

In many cases, $F_0$ and $G_0$ are the identity, and the entire
adjunction is `over \sB' in the sense introduced in
\xref{def:monoidal-mor-fib}.

\begin{eg}\xlabel{disjoint-sections}
  Let \sC\ have finite limits and finite colimits preserved by
  pullback.  Then there is a forgetful lax monoidal morphism
  $\ttRetr_\sC\to \ttArr_\sC$ lying over \sC.  Cartesian arrows are
  given by pullback in both cases, and hence are preserved strongly,
  but opcartesian arrows are given by pushout in $\ttRetr_\sC$ and
  mere composition in $\ttArr_\sC$, hence are preserved only laxly.
  The lax monoidal constraint is given by the quotient map $M\times N
  \to M\exsm N$.

  This forgetful morphism has a left adjoint
  \begin{equation}
    (-)_+\maps \ttArr_\sC\too \ttRetr_\sC.\label{eq:djsec-fibr}
  \end{equation}
  which takes an object $X\to A$ over $A$ to the retraction
  \[A\too X_+=X\sqcup A \too A.
  \]
  We say that the functor $(-)_+$ \emph{adjoins a disjoint section}.
  It is straightforward to check that $(-)_+$ is a strong monoidal
  morphism of bifibrations and the pair satisfies \xref{monfib-adjn}.

  As \sC\ varies, the forgetful morphisms define a 2-natural
  transformation from the 2-functor $\ttRetr$ to the 2-functor
  $\ttArr$, while the morphisms $(-)_+$ form an oplax natural
  transformation $\ttArr\to\ttRetr$.  This remains true upon composing
  with the 2-functor \bbFr, so we obtain framed adjunctions
  \[(-)_+\maps \Span(\sC) \toot \Ex(\sC) \spam U
  \]
  where the right adjoint is 2-natural in \sC\ and the left adjoint is
  oplax natural in \sC.
\end{eg}

\begin{eg}\xlabel{sigma-infty}
  It is essentially shown in~\cite[Ch.~11]{pht} that
  we have an adjunction 
  \begin{equation}
    \Sigma^\infty\maps \ttRetr_\Top \toot \ttSp\spam \Omega^\infty.
    \label{eq:siginf-monfib-adjn}
  \end{equation}
  of monoidal bifibrations lying over \sU.  Applying \bbFr, we obtain
  a framed adjunction
  \begin{equation}
    \Sigma^\infty\maps \Ex(\Top) \toot \Sp\spam \Omega^\infty.\label{eq:siginf-framed-adjn}
  \end{equation}
  The fiber adjunctions are shown to be Quillen in~\cite[12.6.2]{pht},
  so passing to homotopy categories of horizontal 1-cells, we obtain a
  framed adjunction
  \begin{equation}
    \Sigma^\infty\maps \Ho(\Ex(\Top)) \toot \Ex\spam \Omega^\infty.\label{eq:siginf-ho-framed-adjn}
  \end{equation}
\end{eg}

\section{Monoids in monoidal fibrations and examples}
\label{sec:monoids-monfib}

We have seen that well-behaved monoidal bifibrations give rise to
framed bicategories, and that framed bicategories with local
coequalizers admit the \Mod\ construction, so it is natural to ask
what conditions on a monoidal fibration ensure that the resulting
framed bicategory has local coequalizers.

\begin{defn}
  Let $\Phi\maps \sA\to\sB$ be a fibration.
  \begin{itemize}
  \item We say that $\Phi$ has \textbf{fiberwise coequalizers} if each
    fiber $\sA_B$ has coequalizers and the functors $f^*$ preserve
    coequalizers.  Note that this latter condition is automatic if
    $\Phi$ is a $*$-fibration, since then $f^*$ is a left adjoint.
  \item Similarly, we say that $\Phi$ has \textbf{fiberwise equalizers}
    if each fiber $\sA_B$ has equalizers and $f^*$ preserves
    equalizers, the second condition being automatic if $\Phi$ is a
    bifibration.
  \item If $\Phi$ is a monoidal fibration with fiberwise coequalizers,
    we say these coequalizers are \textbf{preserved by $\ten$} if the
    functors $\ten\maps \sA_A\times \sA_B\too \sA_{A\ten B}$ all
    preserve coequalizers in each variable.  This is automatic if the
    right adjoints $\xlhd,\xrhd$ exist.
  \end{itemize}
\end{defn}

\begin{prop}\xlabel{monfib->mod}
  Let $\Phi$ be a frameable monoidal bifibration with fiberwise
  coequalizers preserved by $\ten$.  Then:
  \begin{enumerate}
  \item The framed bicategory $\bbFr(\Phi)$ has local coequalizers, so
    there is a framed bicategory $\Mod(\bbFr(\Phi))$.\label{item:monfib->monoids-i}
  \item If $\Phi$ is symmetric, then $\bbFr(\Phi)$ is a monoidal
    framed bicategory with local coequalizers preserved by $\ten$;
    hence $\Mod(\bbFr(\Phi))$ is also
    monoidal.\label{item:monfib->monoids-ii}
  \item If $\Phi$ is externally closed, a $*$-fibration, and has
    fiberwise equalizers, then $\bbFr(\Phi)$ is closed and its
    hom-categories have equalizers; hence $\Mod(\bbFr(\Phi))$ is also
    closed.\label{item:monfib->monoids-iii}
  \end{enumerate}
\end{prop}
\begin{proof}
  Since the hom-category $\mathcal{F}\mathit{r}(\Phi)(A,B)$ is just
  the fiber of $\Phi$ over $A\times B$, it has coequalizers.  And
  since $M\odot N = \pr_!\Delta^*(M\ten N)$, where $\ten$ and
  $\Delta^*$ preserve coequalizers by assumption and $\pr_!$ preserves
  all colimits as it is a left adjoint, these coequalizers are
  preserved by $\odot$; thus $\bbFr(\Phi)$ has local coequalizers.
  Item~\ref{item:monfib->monoids-i} then follows from
  \xref{frbi->mod}.

  For~\ref{item:monfib->monoids-ii}, we need to know what the external
  monoidal structure of $\bbFr(\Phi)$ is.  When we prove in
  \xref{monfib->framed-i} that $\bbFr(\Phi)$ is monoidal, we will
  define this external product to be essentially that of $\Phi$, but
  with a slight twist.  Namely, if $\Phi(M)=A\times B$ and
  $\Phi(N)=C\times D$, then we have $\Phi(M\ten N)=(A\times
  B)\times(C\times D)$, whereas their product in $\bbFr(\Phi)$ should
  lie over $(A\times C)\times (B\times D)$.  Thus we define the
  external product $M\ten' N$ of $\bbFr(\Phi)$ to be the base change
  of $M\ten N$ along the constraint isomorphism
  \[(A\times C)\times (B\times D) \iso (A\times B)\times(C\times D).\]
  In particular, we have $M\ten' N\iso M\ten N$; thus $\ten'$
  preserves coequalizers because $\ten$ does.  It then follows from
  \xref{frbi->mod} that $\Mod(\bbFr(\Phi))$ is monoidal.

  Finally,~\ref{item:monfib->monoids-iii} follows directly from
  Theorems~\ref{monfib->framed} and~\ref{frbi->mod}.
\end{proof}

\begin{eg}
  If \sC\ has finite limits and coequalizers preserved by pullback,
  then its self-indexing $\ttArr_\sC$ satisfies the conditions of
  \xref{monfib->mod}, and we obtain the framed bicategory
  $\Mod(\Span(\sC))$ of internal categories and distributors which we
  mentioned in \xref{eg:internal-dist}.
\end{eg}

However, we can also obtain \emph{enriched} categories and
distributors, by starting with a different monoidal bifibration.

\begin{eg}\xlabel{eg:fib-can}
  Given any ordinary category \sV, let $\Fam(\sV)$ be the category of
  \emph{families} of objects of \sV.  That is, an object of
  $\Fam(\sV)$ is a set $X$ together with an $X$-indexed family
  $\{A_x\}_{x\in X}$ of objects in \sV.  Then there is a fibration
  $\ttFam_\sV\maps \Fam(\sV)\to\Set$ which is sometimes called the
  \emph{naive indexing} of \sV; its fiber over a set $X$ is the
  category $\sV^X$.  The reader may check the following.
  \begin{itemize}
  \item If \sV\ is a monoidal category, then $\ttFam_\sV$ is a
    monoidal fibration; the external product of $\{A_x\}_{x\in X}$ and
    $\{B_y\}_{y\in Y}$ is $\{A_x\ten B_y\}_{(x,y)\in X\times Y}$.  The
    fiberwise monoidal structure is the obvious one.  If \sV\ is
    braided or symmetric, then so is $\ttFam_\sV$.
  \item If \sV\ has small coproducts (resp.\ products), then
    $\ttFam_\sV$ is a strongly BC bifibration (resp.\ $*$-fibration).
    If \sV\ is also monoidal and $\ten$ preserves coproducts, then
    $\ttFam_\sV$ is a monoidal bifibration.
  \item If \sV\ has coequalizers preserved by $\ten$, then
    $\ttFam_\sV$ has fiberwise coequalizers preserved by $\ten$.  If
    \sV\ has equalizers, then $\ttFam_\sV$ has fiberwise equalizers.
  \item If \sV\ is closed, then $\ttFam_\sV$ is internally closed.
    Thus, by \xref{monfib-closed->loc-and-ext}, if \sV\ also has small
    products, then $\ttFam_\sV$ is externally closed.
  \end{itemize}
  In particular, when \sV\ is monoidal and has colimits preserved by
  $\ten$, $\ttFam_\sV$ is frameable and $\bbFr(\ttFam_\sV)$ has local
  coequalizers, so we can define $\Mod(\bbFr(\ttFam_\sV))$.  It is
  easy to see that $\bbFr(\ttFam_\sV)$ is equivalent to the framed
  bicategory $\Mat(\sV)$ defined in \xref{eg:mat->dist}, where we
  observed that $\Mod(\Mat(\sV))\eqv \Dist(\sV)$.
\end{eg}

This suggests that we can view a frameable closed symmetric monoidal
$*$-bifibration $\Phi$ with fiberwise equalizers and coequalizers as a
`parametrized' version of a complete and cocomplete closed symmetric
monoidal category \sV, and that we can view the associated framed
bicategory $\Mod(\bbFr(\Phi))$ as a parametrized version of
$\Dist(\sV)$.  The other example of $\ttArr_\sC$ seems to bear this
out.

In fact, we can view the monoidal bifibrations $\ttFam_\sV$ and
$\ttArr_\sC$ as living at opposite ends of a continuum.  In
$\ttFam_\sV$, the base category \Set\ is fairly uninteresting, while
all the interesting things happen in the fibers.  On the other hand,
in $\ttArr_\sC$, the base category \sC\ can be interesting, but the
fibers carry essentially no new information, being determined by the
base.  Other monoidal bifibrations will fall somewhere in between the
two, and monoids in the resulting framed bicategories can be thought
of as `categories which are both internal and enriched'.

\begin{eg}\xlabel{pointed-internal-cats}
  If $\sC$ has finite limits and finite colimits preserved by
  pullback, then \xref{monfib->mod} applies to the monoidal
  bifibration $\ttRetr_\sC$, so that $\Mod(\Ex(\sC))$ is a symmetric
  monoidal and closed framed bicategory.  A monoid in $\Ex(\sC)$ may
  be thought of as a `pointed internal category' in \sC.  For example,
  a monoid in $\Ex(\Set)$ is a small category enriched over the
  category $\Set_*$ of pointed sets with smash product, meaning that
  each hom-set has a chosen basepoint and composition preserves
  basepoints.  Similarly, a monoid in $\Ex(\Top)$ is a `based
  topological category'.  If its space of objects is discrete, then it
  is just a small category enriched over based topological spaces, but
  in general it will be `both internal and enriched'.

  Applying \Mod\ to the disjoint-sections functor $(-)_+$ from
  \xref{disjoint-sections}, we obtain a framed functor
  $\Mod(\Span(\sC)) \to \Mod(\Ex(\sC))$.  Thus, any internal category
  can be made into a pointed internal category by `adjoining disjoint
  basepoints to hom-objects'.
\end{eg}

\begin{eg}\xlabel{eg:mod-ex}
  \xref{monfib->mod} applies to the point-set fibration \ttSp\ of
  parametrized spectra, so $\Mod(\Sp)$ is a symmetric monoidal and
  closed framed bicategory.  A monoid in $\Sp$ can be viewed as a
  category `internal to spaces and enriched over spectra'; if its
  space of objects is discrete, then it is just a small category
  enriched over orthogonal spectra.

  To obtain other examples, we can apply $\Mod(\Sigma^\infty)$ to any
  based topological category as in \xref{pointed-internal-cats}, and
  thereby to any internal category in \Top\ with a disjoint section
  adjoined.  Certain monoids in \Sp\ arising in this way from the
  topologized fundamental groupoid $\Pi M$ or path-groupoid
  $\mathcal{P}M$ of a space $M$ play an important role
  in~\cite{kate:traces}.

  A good case can be made (see~\cite{pht}) that a monoid in \Sp\ is
  the right parametrized analogue of a classical ring spectrum, since
  when its space of objects is a point, it reduces to an orthogonal
  ring spectrum.  The more naive notion of a monoid in $\mathbf{Sp}_B$
  with respect to the internal smash product $\sm_B$ is poorly behaved
  because, unlike the situation for the external smash product
  $\exsm$, we have no homotopical control over $\sm_B$.
\end{eg}

The above two examples give framed bicategories with involutions which
are not vertically strict, since the 2-cell components of $\xi$ in
$\bbFr(\Phi)$ are not identities.

\begin{eg}
  \xref{monfib->mod} does \emph{not} apply to the homotopy-level
  monoidal fibration $\Ho(\ttSp)$, since the stable homotopy
  categories of parametrized spectra do not in general admit
  coequalizers.  Rather than $\Mod(\Ex) = \Mod(\Ho(\Sp))$, the correct
  thing to consider is `$\Ho(\Mod(\Sp))$'.  Here the objects are
  honest monoids in \Sp, whose multiplication is associative and
  unital on the point-set level, just like in $\Mod(\Sp)$, but we pass
  to homotopy categories of horizontal 1-cells.  We then need to use a
  `homotopy tensor product' to define the derived horizontal
  composition, as was done in~\cite{kate:traces}.  We hope to
  investigate the homotopy theory of framed bicategories more fully in
  a later paper.
\end{eg}

\begin{eg}
  We can construct various monoidal fibrations over any `\Set-like'
  category \sE\ by mimicking the constructions of classical \Set-based
  monoidal categories.  For example, we have a fibration
  $\mathtt{Ab}_\sE$ over \sE\ whose fiber over $B$ is the category
  $\Ab(\sE/B)$ of abelian group objects in $\sE/B$.  If \sE\ is
  locally cartesian closed, has finite colimits, and the forgetful
  functors $\Ab(\sE/B) \to \sE/B$ have left adjoints, then this is a
  strongly BC $*$-bifibration (see~\cite[D5.3.2]{elephant2}).  If \sE\
  is cocomplete, the tensor product of abelian group objects can be
  defined internally and makes $\mathtt{Ab}_\sE$ a monoidal
  bifibration; monoids in the corresponding framed bicategory
  $\mathbb{A}\mathbf{b}(\sE)$ give a notion of `\Ab-category in \sE'.

  For example, if \sE\ is a category of topological spaces, then any
  vector bundle over a space $B$ gives an object of $\Ab(\sE/B)$.  One
  might argue, analogously to \xref{eg:mod-ex}, that monoids in
  $\mathbb{A}\mathbf{b}(\sE)$ give a good notion of a `bundle of
  rings'.
\end{eg}

The theory of such \emph{relative enriched categories} appears to be
fairly unexplored; the only references we know are~\cite{gg:fib-rel}
and~\cite{przybylek:enriched-internal}.  We will explore this theory
more extensively in a later paper; in many ways, it is very similar to
classical enriched category theory.  We end with one further example
of this phenomenon.

If \sV\ is an ordinary monoidal category with coproducts preserved by
$\ten$, then any small unenriched category $C$ gives rise to a `free'
\sV-category $\sV[C]$ whose hom-objects are given by copowers of the
unit object:
\begin{equation}
  \sV[C](x,y) = \coprod_{C(x,y)} I.\label{eq:free-vcat}
\end{equation}
For a monoidal fibration $\Phi\maps \sA\to\sB$, the analogue of an
unenriched category is an internal category in \sB.  The following is
an analogue of this construction in our general context.

\begin{prop}
  Let \sB\ have finite limits and let $\Phi\maps \sA\to\sB$ be a
  strongly BC monoidal bifibration.  Then there is a canonical strong
  monoidal morphism of bifibrations
  \begin{equation}
    \ttArr_\sB\too\Phi\label{eq:discrete-monfib}
  \end{equation}
  which takes an object $X\too[f] B$ of $\sB/B$ to the object
  $f_!\pr_X^*I$ of $\sA_B$.  Consequently, there is a canonical framed
  functor
  \begin{equation}
    \Span(\sB)\too \bbFr(\Phi)\label{eq:discrete-frbi}
  \end{equation}
  and thus, if $\Phi$ has fiberwise coequalizers preserved by $\ten$
  and \sB\ has coequalizers preserved by pullback, a framed functor
  \begin{equation}
    \Mod(\Span(\sB))\too\Mod(\bbFr(\Phi)).\label{eq:discrete-mod}
  \end{equation}
\end{prop}

\begin{eg}
  When $\Phi = \ttFam_\sV$ for an ordinary monoidal category \sV, the
  morphism $\ttArr_\sB\to\ttFam_\sV$ takes a set $A$ to $\coprod_A I$.
  Thus the induced framed functor $\Dist(\Set) \to \Dist(\sV)$ is
  exactly the `free \sV-category' operation~(\ref{eq:free-vcat})
  described above.
\end{eg}

\begin{egs}
  More interestingly, when $\Phi$ is $\ttRetr_\sC$, the morphism
  $\ttArr_\sC\to\ttRetr_\sC$ is the disjoint-section operation
  described in~\xref{disjoint-sections}.  And when $\Phi$ is the
  monoidal fibration \ttSp\ of orthogonal spectra, the morphism
  $\ttArr_\sU \to \ttSp$ first adjoins a disjoint section, then
  applies the parametrized suspension-spectrum functor $\Sigma^\infty$
  from~\xref{sigma-infty}.  Therefore, if $C$ is an internal category
  in topological spaces, the `topologically internal and spectrally
  enriched category' $\Sigma^\infty C_+$ considered in
  \xref{eg:mod-ex} is in fact `freely generated by $C$' in this
  canonical way.
\end{egs}

\section{Two technical lemmas}
\label{sec:beck-chev}

In preparation for our proof of \xref{monfib->framed} in
\S\ref{sec:monfib-to-framed-ii}, in this section we reformulate the
Beck-Chevalley condition and internal closedness in terms of cartesian
arrows.

\begin{lem}\xlabel{bchv-via-cart}
  Let $\Phi\maps \sA\to\sB$ be a bifibration.  Then a commuting square
  \begin{equation}
    \xymatrix{A \ar[r]^h\ar[d]_k & B \ar[d]^g\\ C \ar[r]_f & D}\label{eq:bchv-downstairs-2}
  \end{equation}
  in \sB\ satisfies the Beck-Chevalley condition if and only if for
  every $M\in \sA_B$, the square~(\ref{eq:bchv-downstairs-2}) lifts to
  some commutative square
  \begin{equation}
    \xymatrix{M' \ar[r]^\phi\ar[d]_\chi & M \ar[d]^\xi\\
      M'''\ar[r]_\psi & M''}\label{eq:bchv-cartsq}
  \end{equation}
  in \sA\ in which $\phi$ and $\psi$ are cartesian and $\chi$ and
  $\xi$ are opcartesian.
\end{lem}
Note that given $\phi,\chi,\xi$ lifting $h,k,g$ with $\chi$
opcartesian, there is exactly one $\psi$ lifting $f$ which
makes~(\ref{eq:bchv-cartsq}) commute.  Thus the condition can also be
stated as ``Given any cartesian $\phi$ and opcartesian $\chi,\xi$, the
unique morphism $\psi$ over $f$ making~(\ref{eq:bchv-cartsq}) commute
is cartesian''.
\begin{proof}
  Choose a cleavage.  Then by the universal properties of cartesian
  and opcartesian arrows, there is a unique dotted arrow living over
  $1_C$ which makes the following pentagon commute:
  \[\xymatrix{& h^*M \ar[r]^{\mathrm{cart}}\ar[dl]_{\mathrm{opcart}} &
    M \ar[d]^{\mathrm{opcart}}\\
    k_!h^*M \ar@{.>}[dr] & & g_!M\\
    & f^*g_!M \ar[ur]_{\mathrm{cart}} .}
  \]
  We claim that in fact this dotted arrow is the component of the
  Beck-Chevalley natural transformation~(\ref{eq:bchv-trans}) at
  $M$.  To see this, we fill out the diagram as follows.
  \[\xymatrix{h^*M \ar[rrrr]^{\mathrm{cart}}\ar[ddd]_{\mathrm{opcart}}
    \ar@{-->}[dr]_{h^*\eta} &&&&
    M \ar[dddd]^{\mathrm{opcart}} \ar@{-->}[dl]^\eta\\
    & h^*g^*g_!M \ar[rr]_{\mathrm{cart}}\ar[dd]^{\mathrm{opcart}}
    \ar@{-->}[dr]^\iso
    && g^*g_!M \ar[dddr]_{\mathrm{cart}}\\
    && k^*f^*g_!M \ar@(r,u)[ddr]^{\mathrm{cart}} \ar[d]_{\mathrm{opcart}}\\
    k_!h^*M \ar@{-->}[r]^(0.4){k_!h^*\eta} \ar@{.>}@(dr,l)[drrr] &
    k_!h^*g^*g_!M \ar@{-->}[r]^\iso &
    k_!k^*f^*g_!M \ar@{-->}[dr]^{\ep f^*g_!}\\
    &&& f^*g_!M \ar[r]_{\mathrm{cart}} &
    g_!M.}
  \]
  Here the dashed arrows are unique factorizations through
  (op)cartesian arrows.
  This exhibits the dotted arrow as the composite of a unit, canonical
  isomorphism, and counit, which is the definition of the
  transformation~(\ref{eq:bchv-trans}).  This proves our claim.

  Therefore, if~(\ref{eq:bchv-trans}) is an isomorphism, the
  composite $k_!h^*M\iso f^*g_!M\to g_!M$ is cartesian, and hence we
  have a commuting square of cartesian and opcartesian arrows as
  desired.  Conversely, if we have such a commuting square, then
  clearly for some choice of cleavage, the dotted arrow is the
  identity; hence it is an isomorphism for all cleavages.
\end{proof}

As always, it simplifies our life greatly to work with cartesian
arrows rather than chosen cleavages.  For example, we can now easily
show the following.

\begin{cor}\xlabel{monfrbi->bchv}
  If \bbD\ is a monoidal framed bicategory, then the square
  \begin{equation}
    \xymatrix{A\ten C \ar[r]^{f\ten 1}\ar[d]_{1\ten g} &
      B\ten C \ar[d]^{1\ten g}\\
      A\ten D\ar[r]_{f\ten 1} & B\ten D}\label{eq:bchv-frbi}
  \end{equation}
  in \bbDz\ satisfies the Beck-Chevalley condition for the
  bifibrations $L$ and $R$.
\end{cor}
\begin{proof}
  Let $M\maps B\ten C \hto E$ be a horizontal 1-cell and consider the
  following diagram in $R^{-1}(E)$.
  \[\xymatrix{({}_fB\ten C) \odot M
    \ar[d]_{\mathrm{opcart}} \ar[r]^{\mathrm{cart}} &
    (B\ten C) \odot M
    \ar[d]^{\mathrm{opcart}} \ar[r]^(0.7)\iso & M\\
    ({}_fB\ten D_g)\odot M \ar[r]_{\mathrm{cart}} &
    (B\ten D_g) \odot M
  }\]
  The arrows labeled cartesian or opcartesian are obtained from the
  cartesian ${}_fB\to B$ and opcartesian $C\to D_g$ via $\ten$ and
  $\odot$.  The square commutes by functoriality of $\ten$, so the
  result follows from \xref{bchv-via-cart}.
\end{proof}

\begin{cor}\xlabel{frbi->weak-bchv}
  If \bbD\ is a monoidal framed bicategory in which \bbDz\ is
  cartesian monoidal, then the monoidal bifibration $(L,R)$ is weakly
  BC.
\end{cor}
\begin{proof}
  Taking $D=1$ in the square~(\ref{eq:bchv-frbi}) shows immediately
  that $L$ and $R$ are weakly BC; an analogous square in
  $\bbDz\times\bbDz$ applies to $(L,R)$.
\end{proof}

To deal with the `weakly BC and internally closed' case of
\xref{monfib->framed}, we also need a statement about cartesian arrows
that makes use of the closed structure.  Recall that if $f^*$ is
strong monoidal and has a left adjoint $f_!$, then $f^*$ is closed
monoidal if and only if the dual maps
\begin{align}
  f_!(M\boxtimes f^*N) \too f_!M\boxtimes N\label{eq:intcl-mate-repeat-1}\\
  f_!(f^*N\boxtimes M) \too N\boxtimes f_!M\label{eq:intcl-mate-repeat-2}
\end{align}
are isomorphisms.  This latter condition is amenable to restatement in
fibrational terms, using the characterization of $\boxtimes$ in terms
of $\ten$ and $\Delta^*$.

\begin{lem}\xlabel{closed->some-bchev}
  Let $\Phi\maps \sE\to\sB$ be an internally closed monoidal bifibration,
  where \sB\ is cartesian monoidal.  Then for any $f\maps A\to B$ in
  \sB, and any $M,N$ in \sE\ with $\Phi(M)=A,\Phi(N)=B$, the square
  \begin{equation}
    \xymatrix{A \ar[r]^f\ar[d]_{\Delta_A} & B \ar[d]^{\Delta_B}\\
      A\times A \ar[r]_{f\times f} & B\times B}\label{eq:closed-bchev-downstairs}
  \end{equation}
  in \sB\ lifts to a square
  \begin{equation}
    \xymatrix{\Delta_A^*(M\ten f^*N) \ar[r]^{\mathrm{opcart}}\ar[d]_{\mathrm{cart}} &
      \Delta_B^*(f_!M\ten N)\ar[d]^{\mathrm{cart}}\\
      M\ten f^*N \ar[r]_{\mathrm{opcart}\ten\mathrm{cart}} & f_!M\ten
      N}\label{eq:closed-bchev-ok}
  \end{equation}
  in \sA, and dually.
\end{lem}
\begin{proof}
  A cleavage gives us an opcartesian $M\to f_!M$, a cartesian $f^*N\to
  N$, and cartesian arrows on the left and right, inducing a unique
  arrow on the top which lifts $f$.  We then observe that
  $\Delta_A^*(M\ten f^*N) \iso M\boxtimes f^*N$ and
  $\Delta_B^*(f_!M\ten N) \iso f_!M\boxtimes N$, so factoring the top
  arrow through an opcartesian arrow gives us
  precisely~(\ref{eq:intcl-mate-repeat-1}).  Since $\Phi$ is
  internally closed, this is an isomorphism, so the top arrow must be
  opcartesian.
\end{proof}

\section{Proofs of Theorems \ref{monfib->framed} and \ref{monfib->framed-2functor}}
\label{sec:monfib-to-framed-ii}

This section is devoted to the proofs of Theorems~\ref{monfib->framed}
and~\ref{monfib->framed-2functor} (and the dual version
\xref{dual-monfib->framed}).  To make the proofs more manageable, we
split them up into several propositions.

\begin{prop}\xlabel{monfib->framed-i}
  Let $\Phi\maps \sA\to\sB$ be a strongly BC monoidal bifibration,
  where \sB\ is cartesian monoidal.  Then there is a framed bicategory
  $\bbFr(\Phi)$ defined as follows.
  \begin{enumerate}
  \item $\bbFr(\Phi)_0=\sB$.
  \item $\bbFr(\Phi)_1$, $L$, and $R$  are defined by the following pullback square.
    \[\xymatrix{\bbFr(\Phi)_1 \ar[r]\ar[d]_{(L,R)} & \sA \ar[d]^\Phi\\
      \sB\times \sB\ar[r]_<>(.5)\times & \sB.}\]
    Thus the horizontal 1-cells $A\hto B$ are the objects of \sA\ over
    $A\times B$, and the 2-cells $M\sto{f}{g} N$ are the arrows of
    \sA\ over $f\times g$.
  \item The horizontal composition of $M\maps A\hto B$ and $N\maps
    B\hto C$ is 
    \[M\odot N = (\pr_B)_!\Delta_B^*(M\ten N),
    \]
    and similarly for 2-cells.
  \item The horizontal unit of $A$ is
    \[U_A = (\Delta_A)_!\pr_A^* I.\]
  \end{enumerate}
  If $\Phi$ is symmetric, then $\bbFr(\Phi)$ is symmetric monoidal.
\end{prop}
\begin{proof}
  Throughout the proof, we will write $\bbP=\bbFr(\Phi)$ for brevity.
  Since we intend to construct an algebraic structure (a framed
  bicategory), we choose once and for all a cleavage (and opcleavage)
  on $\Phi$, and reserve the notations $f^*,f_!$ and so on for the
  functors given by this cleavage.  However, we will still use
  (op)cartesian arrows which are not in this cleavage in order to
  construct the constraints and coherence.

  We have the structure and operations, at least, of a double category
  essentially already defined, except for the functoriality of $\odot$
  and $U$.  It is easy to see that $\odot$ is a functor, since $\ten$,
  $\Delta^*$, and $\pr_!$ are functors.  The functoriality of $U$ is
  similar, but perhaps not as obvious since it is a functor
  $\bbPz\to\bbPo$.  Its action on an arrow $f\maps A\to B$ is given by
  the unique factorization $U_f$ as follows.
  \[\xymatrix{I &
    \pr_A^* I \ar[l]_{\mathrm{cart}}\ar[r]^{\mathrm{opcart}}
    \ar@{-->}[d] & U_A \ar@{-->}[d]^{U_f}\\
    & \pr_B^* \ar[ul]^{\mathrm{cart}}\ar[r]_{\mathrm{opcart}} & U_B}\]
  Thus, to show that \bbP\ is a double category, it suffices to
  construct coherent associativity and unit constraints.  The
  following arguments should remind the reader of the proof of
  \xref{locmon->monfib}, although they are more complicated.

  Note first that for horizontal 1-cells $M\maps A\hto B$ and $N\maps
  B\hto C$ in \bbP, we have $\Phi(M)=A\times B$ and $\Phi(N)=B\times
  C$, and the chosen cleavage gives us canonical morphisms
  \begin{equation}
    M\ten N \oot[\mathrm{cart}] \Delta_B^*(M\ten N)
    \too[\mathrm{opcart}] (\pr_B)_!\Delta_B^*(M\ten N) = M\odot N.
    \label{eq:odot-cart-arr}
  \end{equation}
  We begin with the associativity isomorphism.  So suppose in addition
  to $M,N$ we have $Q\maps C\hto D$.  Then since $\ten$ preserves
  (op)cartesian arrows, we can construct the following diagram.
  \[\small\xymatrix{
    && {} \save[] *{(M\odot N) \ten Q}
    \ar@{<-}[dl]_{\mathrm{opcart}}
    \ar@{<-}[dr]^{\mathrm{cart}}
    \restore \\
    (M\ten N)\ten Q &
    \Delta_B^*(M\ten N) \ten Q \ar[l]_{\mathrm{cart}} &&
    \Delta_C^*\big((M\odot N) \ten Q\big)
    \ar[r]^<>(.5){\mathrm{opcart}} &
    (M\odot N)\odot Q.\\
    && {} \save[] *{\Delta_{BC}^* \big((M\ten N) \ten Q\big)}
    \ar@{-->}[ul]^{\mathrm{cart}} \ar@(l,dr)[ull]^<>(.5){\mathrm{cart}}
    \ar@{.>}[ur]_{\mathrm{opcart}}}  \restore
  \]
  Here the solid arrows are part of the chosen cleavage.  The dashed
  arrow is a unique factorization, which is cartesian by
  \xref{fib-facts}\ref{item:cart-fact-cart}.  The dotted arrow, also a
  unique factorization, is opcartesian by the Beck-Chevalley condition
  (\xref{bchv-via-cart}), because the square in question lifts the
  pullback square
  \[\xymatrix{A \times B\times C \times D \ar[r]^<>(.5){\Delta_C}\ar[d]_{\pr_B} &
    A\times B\times C\times C\times D\ar[d]^{\pr_B}\\
    A\times C\times D\ar[r]_<>(.5){\Delta_C} & A\times C\times C\times D.}\]
  Composing the two opcartesian arrows on the right, we obtain a span
  \begin{equation}
    \xymatrix{(M\ten N)\ten Q &
      \Delta_{BC}^* \big((M\ten N) \ten Q\big)
      \ar[l]_<>(.5){\mathrm{cart}} \ar[r]^<>(.5){\mathrm{opcart}} &
      (M\odot N)\odot Q.}\label{eq:odot-defn-span}
  \end{equation}
  We perform an analogous construction for $M\odot(N\odot Q)$, then
  factor the associativity isomorphism for $\ten$ through these
  cartesian and opcartesian arrows to obtain an associativity
  isomorphism for $\odot$:
  \begin{equation}
    \xymatrix{(M\ten N)\ten Q \ar[d]_{\iso} &
      \Delta_{BC}^* \big((M\ten N) \ten Q\big)
      \ar[l]_<>(.5){\mathrm{cart}} \ar[r]^<>(.5){\mathrm{opcart}} \ar@{.>}[d]_{\iso} &
      (M\odot N)\odot Q \ar@{.>}[d]^{\iso}\\
      M\ten (N\ten Q) &
      \Delta_{BC}^* \big(M\ten (N\ten Q)\big)
      \ar[l]^<>(.5){\mathrm{cart}} \ar[r]_<>(.5){\mathrm{opcart}} &
      M\odot (N\odot Q).}\label{eq:associator}
  \end{equation}
  This isomorphism is natural because it is defined by unique
  factorization.  The proof that it satisfies the pentagon axiom is
  similar to its construction: we tensor~(\ref{eq:odot-defn-span})
  with $R\maps D\hto E$, then use the Beck-Chevalley condition again
  for the square
  \[\xymatrix@R=1pc@C=1pc{ABCDE \ar[r]\ar[d] & ABCDDE \ar[d]\\
    ADE \ar[r] & ADDE}\]
  (where we omit the symbol $\times$) to obtain a span
  \begin{equation}
    \small\xymatrix{((M\ten N)\ten Q)\ten R &
      \Delta_{BCD}^* \big(((M\ten N) \ten Q)\ten R\big)
      \ar[l]_<>(.5){\mathrm{cart}} \ar[r]^<>(.5){\mathrm{opcart}} &
      ((M\odot N)\odot Q)\odot R.}\label{eq:pentagon-span}
  \end{equation}
  By uniqueness of factorizations, the isomorphism
  \[((M\odot N)\odot Q)\odot R \iso (M\odot (N\odot Q))\odot R,\]
  obtained by applying the functor $-\odot R$
  to~(\ref{eq:associator}), is the same as the isomorphism obtained by
  factoring the isomorphism
  \[((M\ten N)\ten Q)\ten R \iso (M\ten (N\ten Q))\ten R\] through the
  span~(\ref{eq:pentagon-span}).  Therefore, by inspecting the
  following diagram and using unique factorization again, we see that
  the pentagon axiom for $\ten$ implies the pentagon axiom for
  $\odot$.
  \[\small\xymatrix@R=1pc@C=1pc{
    &&& \save[] *{((M\odot N)\odot Q)\odot R}="top" \restore \\
    \save[] *{(M\odot N) \odot (Q\odot R)}="topleft" \ar[ddddddd]|\iso \restore &&&
    \ar[u]\ar[d] \ar[drr]|\iso \ar[dll]|\iso &&&
    \save[] *{(M\odot (N\odot Q))\odot R}="topright" \ar[ddddddd]|\iso  \restore\\
    & \ar[ddr] \ar[ul] \ar[ddddd]|\iso &&
    \save[] *{((M\ten N)\ten Q)\ten R} \ar[ddl]|\iso \ar[ddr]|\iso \restore &&
    \ar[ddl] \ar[ur]  \ar[ddddd]|\iso \\
    \\
    && (M\ten N) \ten (Q\ten R) \ar[dd]|\iso &&
    (M\ten (N\ten Q))\ten R \ar[dd]|\iso\\
    \\
    && M\ten (N\ten (Q\ten R)) &&
    M\ten ((N\ten Q)\ten R) \ar[ll]|\iso\\
    & \ar[ur]\ar[dl] &&&& \ar[ul] \ar[dr] \ar[llll]|\iso\\
    \save[] *{M\odot (N\odot (Q\odot R))}="botleft" \restore &&&&&&
    \save[] *{M\odot ((N\odot Q)\odot R)} \ar"botleft"|\iso \restore
    \ar"top";"topleft"|\iso
    \ar"top";"topright"|\iso
  } 
  \]

  Now we consider the left unit transformation.  Let $M\maps A\hto B$
  and recall that $U_A=(\Delta_A)_!\pr_A^*I$, so that we have
  \[\xymatrix{I & \pr_A^*I \ar[l]_{\mathrm{cart}}
    \ar[r]^{\mathrm{opcart}} & U_A}.\]
  Tensoring this with $M$ and adding the arrows
  from~(\ref{eq:odot-cart-arr}) for $U_A\odot M$, we have
  \begin{equation}
    \small\xymatrix{
      && U_A\ten M\\
      I\ten M &
      \pr_A^*I\ten M \ar[l]_{\mathrm{cart}} \ar[ur]^{\mathrm{opcart}} &&
      \Delta_A^*(U_A\ten M) \ar[ul]_{\mathrm{cart}}
      \ar[r]^<>(.5){\mathrm{opcart}} &
      U_A\odot M.\\
      && M
      \ar@{-->}[ul]_{\mathrm{cart}} \ar@(l,dr)[ull]^{\fl}_{\iso}
      \ar@{.>}[ur]_{\mathrm{opcart}}}\label{eq:unit-trans}
  \end{equation}
  The solid arrows marked cartesian or opcartesian are part of the
  chosen cleavage.  The other solid arrow is the left unit constraint
  for $\ten$, which is an isomorphism, hence also cartesian.  The
  dashed arrow is cartesian by
  \xref{fib-facts}\ref{item:cart-fact-cart}, and the dotted arrow is
  opcartesian by the Beck-Chevalley condition for the pullback square
  \begin{equation}
    \xymatrix{A\times B \ar[r]^<>(.5)\Delta\ar[d]_\Delta &
      A\times A\times B \ar[d]^{\Delta\times 1}\\
      A\times A\times B\ar[r]_<>(.5){1\times \Delta} & A \times A\times
      A\times B.
    }\label{eq:bchv-bad-1}
  \end{equation}
  Since the composite of the two opcartesian arrows on the right is
  opcartesian and lies over $1_{A\times B}$, it is an isomorphism
  \[M\iso U_A\odot M
  \]
  which we take as the left unit isomorphism for $\odot$.  Its
  naturality follows, as before, from unique factorization.  The right
  unit isomorphism is analogous.

  We now show the unit axiom.  We tensor the diagram
  \[\xymatrix{I\ten M & M \ar[l]_<>(.5)\iso \ar[r]^<>(.5)\iso & U_A\odot M}\]
  with $N$ and compose with the defining cartesian and opcartesian
  arrows for $\odot$ to obtain the following diagram.
  \begin{equation}
    \small\xymatrix{
      && \save[] *{N\ten (U_A\odot M)}="top" \restore\\
      N\ten (I\ten M) &
      N\ten M \ar[l]_<>(.5){\iso} \ar"top"^{\iso} &&
      \Delta_A^*(N\ten (U_A\odot M)) \ar"top"_{\mathrm{cart}}
      \ar[r]^<>(.5){\mathrm{opcart}} &
      N\odot (U_A\odot M).\\
      && \save[] *{\Delta_A^*(N\ten (I\ten M))}
      \ar@{-->}[ul]^{\mathrm{cart}}
      \ar@{.>}[ur]^<>(.5)\iso_{\mathrm{opcart}} \restore
    }\label{eq:unit-axiom}
  \end{equation}
  We do the same for $(N\odot U_A)\odot M$.  By universal
  factorization, the two unit isomorphisms
  \begin{align*}
    N\odot (U_A\odot M) &\iso N\odot M\\
    (N\odot U_A)\odot M &\iso N\odot M
  \end{align*}
  are given by factorization through these (op)cartesian arrows, as
  is the associativity isomorphism
  \[(N\odot U_A)\odot M \iso N\odot (U_A\odot M).
  \]
  Thus, as for the pentagon, the unit axiom for $\ten$ implies the
  unit axiom for $\odot$.

  This shows that \bbP\ is a double category.  Since the pullback of a
  bifibration is a bifibration, $(L,R)$ is a bifibration.  Thus, by
  \xref{thm:framed}, \bbP\ is a framed bicategory.

  We now assume that $\Phi$ is symmetric and show that \bbP\ is a
  symmetric monoidal framed bicategory.  Since $\bbPz=\sB$, it is
  already (cartesian) symmetric monoidal.  The monoidal structure of
  \bbPo\ is almost the same as that of \sA, but with a slight twist.
  If $M\maps A\hto B$ and $N\maps C\hto D$, so that $\Phi(M)=A\times
  B$ and $\Phi(N)=C\times D$, then we have
  \[\Phi(M\ten N)=(A\times B)\times(C\times D)
  \]
  whereas the product of $M$ and $N$ in \bbP\ should be an object of
  $\sA$ lying over $(A\times C)\times (B\times D)$.  But the chosen
  cleavage gives us a cartesian arrow ending at $M\ten N$ lying over
  the unique constraint
  \[(A\times C)\times (B\times D) \iso (A\times B)\times(C\times D),\]
  and we call its domain $M\ten' N$.  Since cartesian arrows over
  isomorphisms are isomorphisms, we have $M\ten' N \iso M\ten N$.
  Similarly, the unit for \sA\ should be $U_1=(\Delta_1)_!(\pr_1)^*
  I$, and since $\pr_1=1_1$ and $\Delta_1$ is the unique isomorphism
  $1\iso 1\times 1$ we have $U_1\iso I$; we define $I'=U_1$.  The
  constraints and coherence axioms for $\ten$ and $I$ pass across
  these isomorphisms to make \bbPo\ a symmetric monoidal category
  under $\ten'$, with $(L,R)$ a strict symmetric monoidal functor.

  Thus, to make \bbP\ a symmetric monoidal framed bicategory, it
  remains to construct coherent interchange and unit isomorphisms and
  show that the monoidal associativity and unit constraints are framed
  transformations.  Our by-now familiar procedure gives the following
  diagram for the interchange isomorphism.
  \[\xymatrix{(M\ten' P)\ten (N\ten' Q) \ar[d]_{\iso} &
    \Delta^*\big((M\ten' P)\ten (N\ten' Q)\big)
    \ar[l]_<>(.5){\mathrm{cart}} \ar[r]^<>(.5){\mathrm{opcart}}
    \ar[d]_\iso &
    (M\ten' P)\odot (N\ten' Q)\ar[d]^\iso\\
    (M\ten N)\ten' (P\ten Q) &
    \Delta^*(M\ten N)\ten' \Delta^*(P\ten Q)
    \ar[l]^<>(.5){\mathrm{cart}} \ar[r]_<>(.5){\mathrm{opcart}} &
    (M\odot N)\ten(P\odot Q).}
  \]
  For the the unit isomorphism we have
  \[\xymatrix{I\ar[d]_\iso &
    \pr_{A\times B}^*I \ar[l]_{\mathrm{cart}}
    \ar[r]^{\mathrm{opcart}}\ar[d]_\iso &
    U_{A\times B} \ar[d]^\iso\\
    I\ten I &
    \pr_A^*I \ten \pr_B^*I \ar[l]^{\mathrm{cart}} \ar[r]_{\mathrm{opcart}} &
    U_A\ten' U_B.}
  \]
  As before, by factoring known commuting diagrams through cartesian
  and opcartesian arrows, we can show that these constraints are
  framed transformations and satisfy the monoidal coherence axioms.
\end{proof}

\begin{cor}\xlabel{dual-monfib->framed-i}
  If $\Phi\maps \sA\to\sB$ is a strongly co-BC monoidal bifibration
  where \sB\ is cocartesian monoidal, then there is a framed
  bicategory $\bbFr(\Phi)$ defined as in \xref{monfib->framed-i},
  except that composition is given by
  \[M\odot N = \eta^*\nabla_!(M\ten N),
  \]
  units are given by
  \[U_A = \nabla^* \eta_! I,
  \]
  and similarly for the other data.  If $\Phi$ is symmetric, then
  $\bbFr(\Phi)$ is symmetric monoidal.
\end{cor}
\begin{proof}
  Simply apply \xref{monfib->framed-i} to the strongly BC monoidal
  bifibration $\Phi\op\maps \sA\op\to\sB\op$, since $\sB\op$ is
  cartesian monoidal.
\end{proof}

We now consider the case when $\Phi$ is only weakly BC.  Most of the
pullback squares for which we used the Beck-Chevalley condition in
\xref{monfib->framed-i} had one of their legs a product projection, so
those parts of the proof carry over with no problem.  However, there
was one which involved only diagonal maps, and this is the problem
that \xref{closed->some-bchev} was designed to solve.  This is
essentially the same method as that used in~\cite[Ch.~17]{pht} for the
case of $\Ho(\ttSp)$.

\begin{prop}\xlabel{monfib->framed-ii}
  Let $\Phi\maps \sA\to\sB$ be a weakly BC and internally closed monoidal
  bifibration, where \sB\ is cartesian monoidal.  Then the same
  definitions as in \xref{monfib->framed-i} give a framed bicategory,
  which is symmetric monoidal if $\Phi$ is.
\end{prop}
\begin{proof}
  There is only one place in the proof of \xref{monfib->framed} where
  we used a Beck-Chevalley property for a `bad' square: in proving
  that the unit transformation is an isomorphism, using the
  square~(\ref{eq:bchv-bad-1}).  In this case, the dotted arrow
  in~(\ref{eq:unit-trans}) which we want to be opcartesian is defined
  by unique factorization from a square of the form
  \begin{equation}
    \xymatrix{
      M \ar[d]_{\mathrm{cart}} \ar@{.>}[r] &
      \Delta_A^*(U_A\ten M) \ar[d]^{\mathrm{cart}}\\
      \pr_A^*I \ten M \ar[r]_{\mathrm{opcart}\ten 1} &
      U_A\ten M.}\label{eq:unit-def-square}
  \end{equation}
  This is almost of the form~(\ref{eq:closed-bchev-ok}), but not
  quite, since it lies over the square
  \[\xymatrix{AB \ar[r]^{\Delta_A B}\ar[d]_{\Delta_A B} &
    AAB \ar[d]^{A\Delta_A B}\\
    AAB\ar[r]_{\Delta_A AB} & AAAB}
  \]
  which is not of the form~(\ref{eq:closed-bchev-downstairs}).  But we
  can decompose it into another pair of squares:
  \[\xymatrix@C=3pc{AB \ar[rr]^{\Delta_A B}
    \ar[d]^{\Delta_{AB}} &&
    AAB
    \ar[d]_{\Delta_{AAB}}\\
    ABAB \ar[d]^{\pr_B} \ar[rr]|<>(.5){\Delta_A B\Delta_A B} &&
    AABAAB \ar[d]_{\pr_{A}\pr_B}\\
    AAB \ar[rr]_{\Delta_A AB} &&
    AAAB}
  \]
  Here the top square is of the
  form~(\ref{eq:closed-bchev-downstairs}), where $f=\Delta_A B$.  If
  we then construct~(\ref{eq:unit-def-square}) by lifting in stages,
  we obtain
  \begin{equation}
    \xymatrix{
      M \ar[d]^{\mathrm{cart}}
      \ar@{.>}[r] &
      \Delta_A^*(U_A\ten M) \ar[d]_{\mathrm{cart}}
      \\
      \pr_{AB}^*I \ten M \ar@{-->}[r] \ar[d]^{\mathrm{cart}} &
      \pr_B^*U_A\ten \pr_A^*M  \ar[d]_{\mathrm{cart}}\\
      \pr_A^*I \ten M \ar[r]_{\mathrm{opcart}\ten 1} &
      U_A\ten M.}\label{eq:unit-def-rectangle}
  \end{equation}
  where the outer rectangle is the same as~(\ref{eq:unit-def-square}).
  We can then obtain the bottom square as the product of a square
  \[\xymatrix{\pr_{AB}^*I \ar@{-->}[r]\ar[d]_{\mathrm{cart}} &
    \pr^* U_A \ar[d]^{\mathrm{cart}}\\
    \pr_A^*I \ar[r]_{\mathrm{opcart}} & U_A,}\]
  where the dashed arrow is opcartesian by the Beck-Chevalley
  condition, and a square
  \[\xymatrix{M \ar@{-->}[r]\ar@{=}[d] &
    \pr_A^* M \ar[d]^{\mathrm{cart}}\\
    M\ar@{=}[r] & M.}\]
  where the dashed arrow is cartesian by
  \xref{fib-facts}\ref{item:cart-fact-cart}.  Thus the dashed arrow
  in~(\ref{eq:unit-def-rectangle}) is of the form
  $\mathrm{opcart}\ten\mathrm{cart}$, so by \xref{closed->some-bchev},
  the dotted arrow is opcartesian as desired.  Since the unit
  transformation that we have just shown to be an isomorphism is the
  same as the transformation defined in the proof of
  \xref{monfib->framed}, the same proof of the coherence axioms
  applies.
\end{proof}

Note that \xref{frbi->weak-bchv} shows that any monoidal framed
bicategory with cartesian base is weakly BC, so being weakly BC is a
\emph{necessary} condition for the construction of
\xref{monfib->framed} to give a framed bicategory.  We do not know
whether being weakly BC is sufficient for frameability without
closedness, but we suspect not.

\begin{prop}
  If $\Phi$ is a frameable monoidal bifibration, then $\bbFr(\Phi)$
  has a vertically strict involution given by the identity on objects
  and $M\op = \fs^* M$ on 1-cells.  If $\Phi$ is symmetric, this
  involution is symmetric monoidal.
\end{prop}
\begin{proof}
  Left to the reader.
\end{proof}

\begin{prop}\xlabel{monfib->closed-i}
  Let $\Phi\maps \sA\to\sB$ be a frameable monoidal $*$-bifibration
  which is externally closed.  Then the resulting framed bicategory
  $\bbFr(\Phi)$ is closed.
\end{prop}
\begin{proof}
  Define $N\rhd P = N\xrhd\Delta_*\pr^*P$.  Writing \calD\ for the
  horizontal bicategory of $\bbFr(\Phi)$, we have
  \begin{align*}
    \calD(M\odot N,P)
    &= \calD(\pr_!\Delta^*(M\ten N),P)\\
    &\iso \calD(\Delta^*(M\ten N),\pr^* P)\\
    &\iso \calD(M\ten N,\Delta_*\pr^*P)\\
    &\iso \calD(M, N \xrhd \Delta_*\pr^*P)\\
    &= \calD(M,N\rhd P).
  \end{align*}
  The construction of $\lhd$ is similar.
\end{proof}

\begin{prop}\xlabel{dual-monfib->closed-i}
  Let $\Phi$ be an externally closed and strongly BC monoidal
  $*$-bifibration in which \sB\ is cocartesian monoidal.  Then the
  resulting framed bicategory $\bbFr(\Phi)$ is closed.
\end{prop}
\begin{proof}
  Define $N\rhd P = N \xrhd \nabla^*\eta_*P$.  Again writing \calD\
  for the horizontal bicategory of $\bbFr(\Phi)$, we have
  \begin{align*}
    \calD(M\odot N,P)
    &= \calD(\eta^*\nabla_!(M\ten N),P)\\
    &\iso \calD(\nabla_!(M\ten N),\eta_* P)\\
    &\iso \calD(M\ten N,\nabla^*\eta_*P)\\
    &\iso \calD(M, N \xrhd \nabla^*\eta_*P)\\
    &= \calD(M,N\rhd P).
  \end{align*}
  The construction of $\lhd$ is similar.
\end{proof}

Finally, we sketch the proof of \xref{monfib->framed-2functor}.

\begin{prop}\xlabel{monfib->framed-2functor-i}
  The construction of \xref{monfib->framed} extends to a 2-functor
  \[\bbFr\maps \MFfr \too \FrBi\]
  and similarly for oplax and lax morphisms.
\end{prop}
\begin{proofsk}
  Let $F\maps \Phi\to\Psi$ be a morphism in the appropriate domain
  category.  We define $\bbFr(F)$ to be $F_0$ on vertical categories.
  If $M\maps A\hto B$, so that $\Phi(M)=A\times B$ and thus
  $\Psi(F_1(M))=F_0(A\times B)$, we let $\bbFr(F)(M) = (F_\times)_!
  F_1(M)$, where $F_\times\maps F_0(A\times B)\to F_0A\times F_0B$ is
  the unique oplax constraint downstairs (which is an isomorphism if
  $F$ is strong or lax).

  The horizontal composition and units are built out of the monoidal
  structure and the functors $f^*$ and $f_!$, so the lax or oplax
  constraints for these induce lax or oplax constraints for a strong
  double functor.  For example, suppose $F\maps \Phi\to\Psi$ is a lax
  monoidal morphism of fibrations and that $M\maps A\hto B$ and
  $N\maps B\hto C$ are horizontal 1-cells in $\bbFr(\Phi)$.  Then
  $M\odot N$ comes with a diagram
  \begin{equation}
    \xymatrix{M\ten N &
      \Delta^*(M\ten N)\ar[l]_<>(.5){\mathrm{cart}} \ar[r]^<>(.5){\mathrm{opcart}} &
      M\odot N}\label{eq:mfanch-funct-upstairs}
  \end{equation}
  lying over
  \begin{equation}
    \xymatrix{A\times B\times B\times B &
      A\times B \times C \ar[l]\ar[r] &
      A\times C.
    }\label{eq:mfanch-funct-downstairs}
  \end{equation}
  Applying $F$ to~(\ref{eq:mfanch-funct-downstairs}), we obtain the
  following diagram (omitting the symbol $\times$).
  \begin{equation}
    \xymatrix{F(ABBC) & F(ABC) \ar[l]\ar[r] & F(AC)\\
      (FA)(FB)(FB)(FC) \ar[u]_\iso &
      (FA)(FB)(FC) \ar[l]\ar[r] \ar[u]_\iso &
      (FA)(FC) \ar[u]_\iso.}\label{eq:mfanch-funct-downstairs-image}
  \end{equation}
  Applying $F$ to~(\ref{eq:mfanch-funct-upstairs}), and adding the
  defining arrows for $FM\odot FN$, we obtain
  \[\xymatrix{F(M\ten N) &
    F(\Delta^*(M\ten N))\ar[l]_<>(.5){\mathrm{cart}}\ar[r] &
    F(M\odot N)\\
    FM\ten FN \ar[u] &
    \Delta^*(FM\ten FN)\ar[l]^<>(.5){\mathrm{cart}}\ar[r]_<>(.5){\mathrm{opcart}} \ar@{-->}[u] &
    FM\odot FN \ar@{.>}[u].}
  \]
  The dashed and dotted arrows follow by factoring the lax constraint
  of $F$ through the given cartesian and opcartesian arrows.  Since
  $F$ does not preserve opcartesian arrows, the top-right solid arrow
  is not necessarily opcartesian, but this does not matter.  The unit
  constraint is similar.

  Finally, the oplax case is dual to this; the only difference is that
  all the vertical arrows go the other way, and
  in~(\ref{eq:mfanch-funct-downstairs-image}) they are no longer
  isomorphisms.
\end{proofsk}

\appendix

\section{Connection pairs}
\label{sec:folding}

As mentioned in \S\ref{sec:introduction}, the questions which led us
to framed bicategories have been addressed by others in several ways.
In this section we explain how framed bicategories are related to
connection pairs on a double category; in the other appendices we
consider their relationship to various parts of bicategory theory.

For further detail on connection pairs, we refer the reader
to~\cite{dblgpd-xedmod,dbl-thin-conn} and also to~\cite{fiore:pscat},
which proved that connection pairs are equivalent to `foldings'.  Our
presentation of the theory differs from the usual one because we focus
on the pseudo case, which turns out to \emph{simplify} the definition
greatly.  The following terminology is
from~\cite{double-adjoints,pp:spans}.

\begin{defn}\xlabel{def:companion}
  Let \bbD\ be a double category and $f\maps A\to B$ a vertical arrow.
  A \textbf{companion} for $f$ is a horizontal 1-cell ${}_fB\maps A\hto
  B$ together with 2-cells
  \begin{equation*}
    \begin{array}{c}
      \xymatrix{
        \ar[r]|<>(.5){|}^<>(.5){{}_fB} \ar[d]_f \ar@{}[dr]|\Downarrow
        & \ar@{=}[d]\\
        \ar[r]|<>(.5){|}_<>(.5){U_B} & }
    \end{array}\quad\text{and}\quad
    \begin{array}{c}
      \xymatrix{
        \ar[r]|<>(.5){|}^<>(.5){U_A} \ar@{=}[d] \ar@{}[dr]|\Downarrow
        & \ar[d]^f\\
        \ar[r]|<>(.5){|}_<>(.5){{}_fB} & }
    \end{array}
  \end{equation*}
  such that the following equations hold.
  \begin{align*}
    \begin{array}{c}
      \xymatrix{
        \ar[r]|<>(.5){|}^<>(.5){U_A} \ar@{=}[d] \ar@{}[dr]|\Downarrow
        & \ar[d]^f\\
        \ar[r]|<>(.5){{}_fB} \ar[d]_f \ar@{}[dr]|\Downarrow
        & \ar@{=}[d]\\
        \ar[r]|<>(.5){|}_<>(.5){U_B} & }
    \end{array} &= 
    \begin{array}{c}
      \xymatrix{ \ar[r]|<>(.5){|}^<>(.5){U_A} \ar[d]_f
        \ar@{}[dr]|{\Downarrow U_f} &  \ar[d]^f\\
        \ar[r]|<>(.5){|}_<>(.5){U_B} & }
    \end{array}
    &
    \begin{array}{c}
      \xymatrix{
        \ar[rr]|<>(.5){|}^<>(.5){{}_fB} \ar@{}[drr]|\iso \ar@{=}[d] &&
        \ar@{=}[d] \\
        \ar[r]|<>(.5){|}^<>(.5){U_A} \ar@{=}[d] \ar@{}[dr]|\Downarrow &
        \ar[r]|<>(.5){|}^<>(.5){{}_fB} \ar[d]_f \ar@{}[dr]|\Downarrow
        & \ar@{=}[d]\\
        \ar[r]|<>(.5){|}_<>(.5){{}_fB} &
        \ar[r]|<>(.5){|}_<>(.5){U_B} &\\
        \ar[rr]|<>(.5){|}_<>(.5){{}_fB} \ar@{}[urr]|\iso \ar@{=}[u] &&
        \ar@{=}[u]}
    \end{array} &=
    \begin{array}{c}
      \xymatrix{
        \ar[r]|<>(.5){|}^<>(.5){{}_fB} \ar@{=}[d]
        & \ar@{=}[d]\\
        \ar[r]|<>(.5){|}_<>(.5){{}_fB} & }
    \end{array}
  \end{align*}
  A \textbf{conjoint} for $f$ is a horizontal 1-cell $B_f\maps B\hto
  A$ together with 2-cells
  \[\begin{array}{c}
    \xymatrix{
      \ar[r]|<>(.5){|}^<>(.5){B_f} \ar@{=}[d] \ar@{}[dr]|\Downarrow
      & \ar[d]^f\\
      \ar[r]|<>(.5){|}_<>(.5){U_B} & }
  \end{array}\quad\text{and}\quad
  \begin{array}{c}
    \xymatrix{
      \ar[r]|<>(.5){|}^<>(.5){U_A} \ar[d]_f \ar@{}[dr]|\Downarrow
      & \ar@{=}[d]\\
      \ar[r]|<>(.5){|}_<>(.5){B_f} & }
  \end{array}\]
  such that the following equations hold.
  \begin{align*}
    \begin{array}{c}
      \xymatrix{
        \ar[r]|<>(.5){|}^<>(.5){U_A} \ar[d]_f \ar@{}[dr]|\Downarrow
        & \ar@{=}[d]\\
        \ar[r]|<>(.5){B_f} \ar@{=}[d] \ar@{}[dr]|\Downarrow
        & \ar[d]^f\\
        \ar[r]|<>(.5){|}_<>(.5){U_B} & }
    \end{array} &= 
    \begin{array}{c}
      \xymatrix{ \ar[r]|<>(.5){|}^<>(.5){U_A} \ar[d]_f
        \ar@{}[dr]|{\Downarrow U_f} &  \ar[d]^f\\
        \ar[r]|<>(.5){|}_<>(.5){U_B} & }
    \end{array}
    &
    \begin{array}{c}
      \xymatrix{
        \ar[rr]|<>(.5){|}^<>(.5){B_f} \ar@{}[drr]|\iso \ar@{=}[d] &&
        \ar@{=}[d] \\
        \ar[r]|<>(.5){|}^<>(.5){{}_fB} \ar@{=}[d] \ar@{}[dr]|\Downarrow &
        \ar[r]|<>(.5){|}^<>(.5){U_A} \ar[d]_f \ar@{}[dr]|\Downarrow
        & \ar@{=}[d]\\
        \ar[r]|<>(.5){|}_<>(.5){U_B} &
        \ar[r]|<>(.5){|}_<>(.5){B_f} &\\
        \ar[rr]|<>(.5){|}_<>(.5){B_f} \ar@{}[urr]|\iso \ar@{=}[u] &&
        \ar@{=}[u]}
    \end{array} &=
    \begin{array}{c}
      \xymatrix{
        \ar[r]|<>(.5){|}^<>(.5){B_f} \ar@{=}[d]
        & \ar@{=}[d]\\
        \ar[r]|<>(.5){|}_<>(.5){B_f} & }
    \end{array}
  \end{align*}
\end{defn}

Comparing this definition with
\xref{thm:framed}\ref{item:frbi-compconj}, the following becomes
evident.

\begin{thm}\xlabel{frbi<->connpair}
  A double category is a framed bicategory exactly when every vertical
  arrow has both a companion and a conjoint.
\end{thm}

One can prove in general that companions and conjoints are unique up
to canonical isomorphism, that $({}_fB,B_f)$ is a dual pair if both
are defined, and that the operations $f\mapsto {}_fB$ and $f\mapsto
B_f$ are pseudofunctorial insofar as they are defined.

The following definition is then easily seen to be equivalent to those
given in~\cite{dblgpd-xedmod,dbl-thin-conn,fiore:pscat}.  Because it
was originally motivated by double categories like the `quintets' of a
2-category, it includes only companions and not conjoints.

\begin{defn}
  Let \bbD\ be a strict double category.  A \textbf{connection pair}
  on \bbD\ is a choice of a companion ${}_fB$ each vertical arrow $f$
  such that the pseudofunctor $f\mapsto {}_fB$ is a strict 2-functor.
\end{defn}

Thus, an arbitrary choice of companions on a non-strict double
category may be called a `pseudo connection pair', and a choice of
conjoints may be called a `pseudo op-connection pair'.
\xref{frbi<->connpair} then states that a double category is a framed
bicategory precisely when it admits both a pseudo connection pair and
a pseudo op-connection pair.

\section{Biequivalences, biadjunctions, and monoidal bicategories}
\label{sec:frbi-vs-bicat}

We now consider the question of how much of the structure of a framed
bicategory \bbD\ is reflected in its underlying bicategory \calD.
Note that any bicategory may be considered as a framed bicategory with
only identity vertical arrows; we call such framed bicategories
\textbf{vertically discrete}.  If $\FrBi_0$ denotes the underlying
1-category of \FrBi\ and $\bfBicat$ denotes the 1-category of
bicategories and pseudo 2-functors, we have an
adjunction
\begin{equation}\label{eq:dbl->bi}
  \bfBicat \toot \FrBi_0.
\end{equation}
in which the left adjoint considers a bicategory as a vertically
discrete framed bicategory, while the right adjoint takes a framed
bicategory to its underlying horizontal bicategory.

The left adjoint $\bfBicat\to\FrBi_0$ does not extend to a 2-functor
or 3-functor, but in the other direction, any framed transformation
$\alpha\maps F\to G\maps \bbD\to\bbE$ can be `lifted' to an oplax
transformation between the underlying pseudofunctors as follows.  For
an object $A\in\bbD$, we define
\[\altil_A = (\alpha_A)^*(U_{GA})\maps FA \hto GA.\]
For a horizontal 1-cell $M\maps A\hto B$, we define
\[\altil_M\maps FM\odot \altil_B \iso (FM)(\alpha_B)_! \too
(\alpha_A)^*(GM) \iso \altil_A\odot GM
\]
to be the globular 2-cell corresponding to $\alpha_M\maps
FM\sto{\alpha_A}{\alpha_B}GM$.  It is easy to check that this defines
an oplax natural transformation between the pseudo 2-functors induced
by $F$ and $G$.

If \calD\ and \calE\ are bicategories, we write $\Bicatc(\calD,\calE)$
for the bicategory of pseudo 2-functors, oplax natural transformations,
and modifications from \calD\ to \calE.  By the pseudofunctoriality of
base change, \xref{bco-psfr}, the above construction defines a
pseudofunctor
\begin{equation}
  \FrBi(\bbD,\bbE)\too \Bicatc(\calD,\calE).\label{eq:fibr->bicat}
\end{equation}
We can also allow lax or oplax functors on both sides.  Note, however,
that framed transformations always give rise to \emph{oplax} natural
transformations.

We would like to say that this construction extends to a functor from
\FrBi\ to the tricategory of bicategories, but unfortunately there is
no tricategory of bicategories which includes oplax natural
transformations, since the composition operation
\[\Bicatc(\calF,\calE)\times\Bicatc(\calD,\calE)\to\Bicatc(\calD,\calF)\]
would be only an oplax 2-functor.  We could allow the codomain to be a
sort of `oplax tricategory', such as the `bicategory op-enriched
categories' of~\cite[1.3]{verity:base-change}, but this would take us
too far afield.  Instead, we merely observe that if $\alpha$ happens
to be a framed natural isomorphism, then $\altil$ is a pseudo natural
equivalence.  This suffices to prove the following.

\begin{prop}
  An equivalence of framed bicategories induces a biequivalence of
  horizontal bicategories.
\end{prop}
\begin{proof}
  If $F,G$ are inverse equivalences in \FrBi, then they give rise to
  pseudofunctors, and by the above observation, the framed natural
  isomorphisms $FG\iso \Id$ and $\Id\iso GF$ give rise to pseudo
  natural equivalences.
\end{proof}

For example, this implies that in \xref{eg:ex-not-equiv}, the
horizontal bicategories $G\calEx_{G/H}$ and $H\calEx$ actually
\emph{are} biequivalent.  However, we believe the equivalence is more
naturally stated, and easier to work with, in \FrBi.

In a similar way, we can lift a monoidal structure on a framed
bicategory to a monoidal structure on its horizontal bicategory.  Many
examples of monoidal bicategories actually arise from monoidal framed
bicategories.  This is useful, because monoidal bicategories are
complicated `tricategorical' objects, whereas monoidal framed
bicategories are much easier to get a handle on.
See~\cite{tricats,nick:tricats,cg:degeneracy-ii} for a definition of monoidal
bicategory.

\begin{thm}\xlabel{monfrbi->monbi}
  If \bbD\ is a monoidal framed bicategory, then any cleavage for
  \bbD\ makes \calD\ into a monoidal bicategory in a canonical way.
\end{thm}
\begin{proofsk}
  \calD\ already has a product and a unit object induced from \bbD, so
  it suffices to construct the constraints and coherence.  We consider
  the associativity constraints, leaving the unit constraints to the
  reader.  Since \bbD\ is a monoidal double category, it has a
  \emph{vertical} associativity constraint
  \[\fa\maps (A\ten B)\ten C\iso A\ten(B \ten C).\]
  But since \bbD\ is a framed bicategory, this vertical isomorphism
  can be `lifted' to an equivalence in \calD:
  \[\fatil = \fa^*\big((A\ten B)\ten C\big)\maps (A\ten B)\ten C\hto A\ten(B \ten C)
  \]
  with adjoint inverse $\big((A\ten B)\ten C\big)\fa^*$; this will be
  the associativity equivalence for the monoidal bicategory \calD.  We
  need further a `pentagonator' 2-isomorphism
  \[\fatil\circ\fatil\circ\fatil \iso \fatil\circ\fatil.\]
  But the coherence pentagon for the vertical isomorphism $\fa$ tells
  us that 
  \[\fa\circ\fa\circ \fa = \fa \circ \fa\]
  and since base change objects are pseudofunctorial by
  \xref{bco-psfr}, this \emph{equality} in \bbDz\ becomes a canonical
  \emph{isomorphism} in \calD, which we take as the pentagonator.

  It remains to check that this pentagonator satisfies the `cocycle
  equation' for relations between quintuple products.  However, since
  all the pentagonators are defined by universal properties (being
  canonical isomorphisms between two cartesian arrows), both sides of
  the cocycle equation are also characterized by the same universal
  property, and therefore must be equal.
\end{proofsk}

In contrast to these well-behaved cases, a framed adjunction $F\maps
\bbD\toot\bbE\spam G$ does not generally give rise to a biadjunction
$\calD\toot\calE$.  It does, however, give rise to a \emph{local
  adjunction} in the sense of~\cite{bp:local-adjointness}; this
consists of an oplax 2-functor $F\maps \calD\to\calE$, a lax 2-functor
$G\maps \calE\to\calD$, and an adjunction
\begin{equation}
  \calD(A,GB) \toot \calE(FA,B).\label{eq:local-adjunction}
\end{equation}
In a biadjunction, $F$ and $G$ would be pseudo 2-functors
and~\eqref{eq:local-adjunction} would be an equivalence.

When $F$ and $G$ arise from a framed or op-framed adjunction, a local
adjunction~\eqref{eq:local-adjunction} is given by
\[\xymatrix@C=3pc{\calD(A,GB) \ar@<1mm>[r]^{(F-)\ep_!} & \calE(FA,B) \ar@<1mm>[l]^{\eta^*(G-)}}.\]
Of course, in a framed adjunction $F$ is strong, while in an op-framed
adjunction $G$ is strong.  A bit more 2-category theory than we have
discussed here (see~\cite{kelly:doc-adjn}) gives us a notion of
`lax/oplax' framed adjunction, in which the left adjoint is oplax and
the right adjoint is lax; these also give rise to local adjunctions
between horizontal bicategories.

In this way, practically any framed-bicategorical notion gives rise to
a counterpart on the purely bicategorical level.  For example, by a
process similar to that in \xref{monfrbi->monbi}, any involution on
\bbD\ gives rise to a `bicategorical involution' on \calD.

Of course, we can also define monoids and bimodules in any bicategory;
in this context monoids are often called \emph{monads}, since in \Cat\
they reduce to the usual notion of monad.  The fact that both internal
and enriched categories are monoids in appropriate bicategories is
well-known, and bicategory theorists have studied categories enriched
in a bicategory as a generalization of categories enriched in a
monoidal category;
see~\cite{sheaves-cauchy-1,street:cauchy-enr,street:absolute,street:enr-cohom,modulated-bicats,2sided-enrichment}.

However, pure bicategory theory usually starts to break down whenever
we need to use vertical arrows.  For example, it is harder to get a
handle on internal or enriched \emph{functors} purely bicategorically.
In the next appendix we introduce a structure called an
\emph{equipment} which is sometimes used for this purpose, for example
in~\cite{ftm2}.

\section{Equipments}
\label{sec:equipments}

For the theory of equipments we refer the reader
to~\cite{proarrows-i,proarrows-ii,modulated-bicats,verity:base-change}.
From our point of view, it is natural to introduce them by asking how
the vertical arrows of \bbD\ are reflected in \calD.  We know that
there is a pseudofunctor $\bbDz\to\calD$ sending $f\maps A\to B$ to
the base change object ${}_fB$; this pseudofunctor is bijective on
objects and each ${}_fB$ has a right adjoint in \calD.  Thus we almost
have an instance of the following structure.

\begin{defn}\xlabel{def:equipment}
  A \textbf{proarrow equipment} is a pseudo 2-functor $\overline{(-)}\maps
  \calK\to\calM$ between bicategories such that
  \begin{enumerate}
  \item \calK\ and \calM\ have the same objects and $\overline{(-)}$ is the
    identity on objects;
  \item For every arrow $f$ in \calK, $\fbar$ has a right adjoint $\ftil$
    in \calM; and
  \item $\overline{(-)}$ is locally full and
    faithful.\label{item:equip-lff}
  \end{enumerate}
\end{defn}

The only difference is that in an equipment, \calK\ is a bicategory
rather than the 1-category \bbDz, but condition~\ref{item:equip-lff}
means that the 2-cells in \calK\ are determined by those in \calM\
anyway.  Thus, given a framed bicategory \bbD, we can factor the
base-change object pseudofunctor $\bbDz\to\calD$ as
\[\bbDz\too[i] \calK \too[\overline{(-)}] \calD
\]
where $i$ is bijective on objects and 1-cells and $\overline{(-)}$ is
locally full and faithful.  The objects and morphisms of \calK\ are
those of \bbDz, and its 2-cells from $f\to g$ are the 2-cells
${}_fB\to {}_gB$ in \calD.  We have proven the following.

\begin{prop}\xlabel{frbi->equip}
  If \bbD\ is a framed bicategory, then the above pseudofunctor
  $\overline{(-)}$ is a proarrow equipment.
\end{prop}

Note that in the proarrow equipment arising from a framed bicategory,
the bicategory \calK\ is actually a strict 2-category.  However, this
is essentially the only restriction on the equipments which arise in
this way.

\begin{prop}\xlabel{equip->frbi}
  Let $\overline{(-)}\maps \calK\to\calM$ be a proarrow equipment such
  that \calK\ is a strict 2-category.  Define a double category \bbD\
  whose
  \begin{itemize}
  \item Objects are those of \calK\ (and \calM);
  \item Vertical arrows are the arrows of \calK;
  \item Horizontal 1-cells are the arrows of \calM; and
  \item 2-cells
    \[\xymatrix{A \ar[r]^M\ar[d]_f \ar@{}[dr]|{\Downarrow\alpha} & B \ar[d]^g\\
      C\ar[r]_N & D}\]
    are the 2-cells $\alpha\maps M \odot \gbar \too \fbar \odot N$ in
    \calM.
  \end{itemize}
  Then \bbD\ is a framed bicategory.
\end{prop}
\begin{proofsk}
  First we show that \bbD\ is a double category.  The vertical
  composite
  \[\xymatrix{\ar[r]^M\ar[d]_f \ar@{}[dr]|{\Downarrow\alpha} & \ar[d]^g\\
    \ar[r]|N\ar[d]_h \ar@{}[dr]|{\Downarrow\beta} & \ar[d]^k\\
    \ar[r]_P &}\]
  is defined to be the composite
  \begin{equation*}
    \begin{array}{rcl}
      M \odot \overline{kg} &\iso& M \odot (\gbar \odot \kbar)\\
      &\iso& (M\odot \gbar) \odot \kbar\\
      &\too[\alpha] & (\fbar \odot N) \odot \kbar\\
      &\iso& \fbar \odot (N\odot \kbar)\\
      &\too[\beta] & \fbar \odot (\overline{h} \odot P)\\
      &\iso& (\fbar \odot \overline{h}) \odot P\\
      &\iso& \overline{hf} \odot P.
    \end{array}
  \end{equation*}
  The coherence theorems for bicategories and pseudofunctors imply
  that this is vertically associative and unital.  Horizontal
  composition of 1-cells is defined as in \calM, and horizontal
  composition of 2-cells is defined analogously to their vertical
  composition.  The constraints come from those of \calM.

  Finally, for an arrow $f\maps A\to B$ in \calK, and \ftil\ the right
  adjoint of \fbar, it is easy to check that the 2-cells
  \begin{equation*}
    \begin{array}{c}
      \xymatrix{
        \ar[r]^<>(.5){\fbar} \ar[d]_f \ar@{}[dr]|{\Downarrow 1_{\fbar}}
        & \ar@{=}[d]\\
        \ar[r]_<>(.5){U_B} & }
    \end{array}\quad,\quad
    \begin{array}{c}
      \xymatrix{
        \ar[r]^<>(.5){\ftil} \ar@{=}[d] \ar@{}[dr]|{\Downarrow\ep}
        & \ar[d]^f\\
        \ar[r]_<>(.5){U_B} & }
    \end{array}\quad,\quad
    \begin{array}{c}
      \xymatrix{
        \ar[r]^<>(.5){U_A} \ar[d]_f \ar@{}[dr]|{\Downarrow \eta}
        & \ar@{=}[d]\\
        \ar[r]_<>(.5){\ftil} & }
    \end{array}, \quad\text{and}\quad
    \begin{array}{c}
      \xymatrix{
        \ar[r]^<>(.5){U_A} \ar@{=}[d] \ar@{}[dr]|{\Downarrow 1_{\fbar}}
        & \ar[d]^f\\
        \ar[r]_<>(.5){\fbar} & }
    \end{array},
  \end{equation*}
  defined by identities and by the unit and counit of the adjunction
  $\fbar\adj\ftil$, satisfy the equations of
  \xref{thm:framed}\ref{item:frbi-compconj}.  Thus \bbD\ is a framed
  bicategory.
\end{proofsk}

At the level of objects, it is easy to show that the two constructions
are inverses up to isomorphism.  In order to state this as an
equivalence of 2-categories, however, we would need morphisms and
especially transformations between equipments, and it is not
immediately obvious how to define these.

The approach to constructing a 2-category of equipments taken
in~\cite{verity:base-change} is essentially to first make equipments
into double categories, as we have done, and \emph{define} morphisms
and transformations of equipments to be morphisms between the
corresponding double categories.  This makes our desired equivalence
true by definition.  Actually,~\cite{verity:base-change} uses `doubly
weak' double categories to deal with equipments where \calK\ is not a
strict 2-category, and thus obtains a tricategory rather than a
2-category, but the idea is the same.

Thus, framed bicategories can be regarded as a characterization of the
double categories which arise from equipments.  However, since the
correct notions of morphism and transformation are apparent only from
the side of double categories, we believe it is more natural to work
directly with framed bicategories.

\begin{rmk}
  The authors of~\cite{basechange-i,basechange-ii} consider a related
  notion of `equipment' where \calK\ is replaced by a 1-category but
  the horizontal composition is forgotten.  If \bbD\ is a framed
  bicategory, then the span
  \begin{equation}
    \bbDz\oot[L]\bbDo\too[R]\bbDz\label{eq:starred-ptd-equip}
  \end{equation}
  has the property that $L$ is a fibration, $R$ is an opfibration, and
  the two types of base change commute, making it into a `two-sided
  fibration' from \bbDz\ to \bbDz\ in the sense of~\cite{street:fibi};
  these are essentially what~\cite{basechange-ii} studies under the
  name `equipment'.  The fact that $L$ is also an opfibration, and $R$
  a fibration, in a commuting way, make~(\ref{eq:starred-ptd-equip})
  into what they call a \emph{starred pointed equipment}.  This
  structure incorporates less of the structure of a framed bicategory,
  but it was sufficient in~\cite{basechange-i,basechange-ii} to obtain
  a 2-category or tricategory of equipments and a notion of equipment
  adjunction.  It is easy to check that any framed adjunction gives
  rise to an equipment adjunction in their sense.
\end{rmk}

\section{Epilogue: framed bicategories versus bicategories}
\label{sec:epilogue}

We end with some more philosophical remarks about the relationship of
framed bicategories to pure bicategory theory.  For any bicategory
\calB, there is a canonical proarrow equipment $\overline{(-)}\maps
\calK\to\calB$, where \calK\ is the bicategory of adjunctions
$\fbar\adj \ftil$ in \calB.  When \calB\ is a strict 2-category, so is
\calK, and the resulting framed bicategory is what we called
$\Adj(\calB)$ in \xref{eg:adj}.  In general, we obtain a `doubly weak'
framed bicategory which we also call $\Adj(\calB)$.

Thus, we can regard the theory of framed bicategories, or of
equipments, as a generalization of the theory of bicategories in which
we specify which adjunctions are the base change objects, rather than
using all of them.  A certain amount of pure bicategory theory can be
regarded as implicitly working with the framed bicategory
$\Adj(\calB)$; frequently 1-cells with right adjoints are called
\emph{maps} and take on a special role.  See, for
example,~\cite{street:cauchy-enr} and~\cite{axiom-mod}.

This purely bicategorical approach works well in bicategories like
$\calDist(\sV)$, because, as we mentioned in
Example~\ref{egs:dual-pairs}\ref{item:dual-pairs-dist}, the mild
condition of `Cauchy completeness' on the \sV-categories involved is
sufficient to ensure that any distributor with a right adjoint is
isomorphic to a base change object.  However, in other framed
bicategories, such as \Mod\ and \Ex, there will not be a good supply
of `Cauchy complete' objects, so framed bicategories or equipments are
necessary.  Moreover, even when working with $\calDist(\sV)$, framed
bicategories are implicit in some of the bicategorical literature,
such as the `calculus of modules' for enriched categories; see, for
example,~\cite{yoneda,proarrows-i,axiom-mod,street:absolute}.

Finally, framed bicategories are much easier to work with than
ordinary bicategories or equipments, because the vertical arrows form
a strict 1-category rather than a weak bicategory.  In situations
where this fails, we can still use `doubly weak' framed bicategories,
as in~\cite{verity:base-change}, but a good deal of simplicity is
lost.  However, in almost all examples, this strictness property does
hold, and the virtue of framed bicategories is that they take
advantage of this fact.  For example, this is what enables us to
define the strict 2-category \FrBi\ and apply the powerful methods of
2-category theory, rather than having to delve into the waters of
tricategories.

\bibliographystyle{halpha}
\bibliography{all,framed}

\end{document}